\documentclass[12pt,a4paper]{book}
\usepackage{amssymb}     
\DeclareFontShape{OMX}{cmex}{m}{n}{<-> sfixed * cmex10}{}
\usepackage{graphicx}
\usepackage{makeidx}
\usepackage[hidelinks]{hyperref}

\newenvironment{proof}%
    {\par\noindent{\it Proof. }\nopagebreak\normalsize}%
    {\linebreak[2]\hspace*{\fill}$\square$\ifdim\lastskip<12pt
    \removelastskip \penalty-200 \vskip12pt \fi}

\newtheorem{thm}{Theorem}[chapter]
\newtheorem{prop}[thm]{Proposition}
\newtheorem{cor}[thm]{Corollary}
\newtheorem{lem}[thm]{Lemma}

\newlength{\retraitexo}
\newlength{\largeurexo}
\newtheorem{exo}[thm]{Exercise}%
\newcommand{\debexo}{\begin{exo} \rm\parshape=2 0cm \textwidth
\retraitexo
\largeurexo }
\newcommand{\finexo}{\par \end{exo}}

\renewcommand{\Bbb}{\mathbb}
\newcommand{\frakt}{\mathfrak}

\newlength{\dessinw}\newlength{\dessinh}
\def\dessin #1 by #2 (#3){%
  \setlength{\dessinw}{#1}\setlength{\dessinh}{#2}\dessinimg{#3}}
\newlength{\dessinmargin}
\newcommand{\dessinimg}[1]{%
  \vbox to \dessinh{\vfill
    \hbox to \dessinw{\kern-\dessinmargin
      \raisebox{-\dessinmargin}[0pt][0pt]{\includegraphics{#1}}\hss}}}

\def\strictsubset{\mathrel{\subsetneq}}
\def\lc{\mathop{\mbox{\rm lc}}\nolimits}
\def\sign{\mathop{\mbox{\rm sign}}\nolimits}
\def\clos{\mathop{\mbox{\rm clos}}\nolimits}
\def\grad{\mathop{\overrightarrow{\mbox{\rm grad}}}\nolimits}
\def\et{\mathrel{\mbox{\rm\ and\ }}}
\def\ou{\mathrel{\mbox{\rm\ or\ }}}
\def\shitem{\setlength{\itemsep}{0pt plus1pt minus1pt}}

\newcommand{\test}[1]{%
\begin{picture}(0,0)
\put(0,6){\line(5,2){15}}
\put(0,6){\line(5,-2){15}}
\put(15,12){\line(5,-2){15}}
\put(15,0){\line(5,2){15}}
\put(15,6){\makebox(0,0){#1}}
\end{picture} }

\let\origtheindex\theindex
\renewcommand{\theindex}{\origtheindex\addcontentsline{toc}{chapter}{Index}}
\makeindex

\begin{document}
\title{AN INTRODUCTION\\
TO SEMIALGEBRAIC GEOMETRY}
\author{Michel COSTE\\ Institut de Recherche Math\'ematique de Rennes\\
michel.coste@univ-rennes1.fr}
\date{Octobre 2002}
\hypersetup{pageanchor=false}
\maketitle
\hypersetup{pageanchor=true}
\setcounter{page}{1}
\tableofcontents
\IfFileExists{SAGintro.tex}{\chapter*{Introduction}
\addcontentsline{toc}{chapter}{Introduction}
\markboth{INTRODUCTION}{INTRODUCTION}

Semialgebraic geometry is the study of sets of real solutions of systems of
polynomial equations and inequalities. These notes present the first results of
semialgebraic geometry and related algorithmic issues. Their content is by no
means original.

The first chapter explains algorithms for counting real roots (this is, in
some sense, $0$-dimensional semialgebraic geometry) and the Tarski-Seidenberg
theorem.

In the next two chapters we study semialgebraic subsets of ${\Bbb R}^n$, which are
defined by boolean combinations of polynomial equations and inequalities. The
main tool for this study is the cylindrical algebraic decomposition, which is
introduced in Chapter 2. The principal result of Chapter 3 is the triangulation
theorem. This theorem shows that semialgebraic sets have a simple topology,
which can be effectively computed from their definitions.

Chapter 4 contains examples of finiteness results and uniform bounds for
semialgebraic families. In particular, we give an explicit bound for the number
of connected components of a real algebraic set as a function of the degree of
the equations and the dimension of the ambient space.

Many results of semialgebraic geometry also hold true for o-minimal structures
(including for instance classes of sets definable with the exponential
functions). See \cite{D} (or the lecture notes \cite{Co}) for this theory. The
algorithmic aspects are specific to the semialgebraic case.

The bibliography is reduced to a minimum, containing mainly books and
recent surveys. Further references (in particular, references to original papers)
may be found there.

These notes have served as a basis for courses in Rennes, Paris and Pisa, for
mini-courses at MSRI (Berkeley), at a summer school in Laredo (Cantabria)
and at a CIMPA school in Niamey.
}{}
\chapter[Counting the real roots]{Counting the real roots of a polynomial}

\section{Sturm's theorem}

In this section, $P \in {\Bbb R}[X]$ is a nonconstant polynomial in one variable.

\subsection{$P$ without multiple root}
We assume that $P$ has no multiple root,
i.e., $\mbox{\rm gcd} (P, P') = 1$. We construct a sequence of polynomials in the following way:
$P_0 = P$, 
$P_1 = P'$ and, for $i > 0$,  $P_{i+1}$ is the negative of the remainder of the euclidean division
of
$P_{i-1}$ by
$P_i$ ($P_{i-1} = P_i Q_i - P_{i+1}$, with 
$\deg P_{i+1} < \deg P_i$). We stop just before we get  $0$.  The last polynomial $P_K$ is then a
nonzero constant (up to signs, this is just Euclide's algorithm for computing the gcd). The sequence
$P_0,\ldots,P_K$ is called the {\em Sturm sequence\index{Sturm sequence} of $P$ and $P'$}. Let
$a \in {\Bbb R}$, not a root of $P$. Denote by $v_P(a)$ the number of sign changes in the sequence
$P_0(a), P_1(a) ,\ldots, P_K(a)$. For instance, if $P=X^3-3X+1$, the Sturm sequence of $P$ and
$P'$ is
$$(X^3-3X+1,\ 3X^2-3,\ 2X-1,\ {9\over4})\;,$$
which gives $(-1,0,1,{9\over4})$ by evaluation at $a=1$. Here $v_P(1)=1$.
We drop the zeroes which occur when counting the sign changes.
The result is the following:

\begin{thm}[Sturm]\label{Sturm}
Let $a < b$ in ${\Bbb R}$, neither $a$ nor $b$ being a root of $P$. The number of roots of  $P$ in
the interval
$(a,b)$ is equal to $v_P(a) - v_P(b)$.
\end{thm}

 \begin{proof} We consider how $v_P(x)$ changes when $x$ passes through a root $c$ of a polynomial of
the Sturm sequence.
\begin{itemize}
\item If $c$ is a root of $P$, the signs of $P_0$ and $P_1$ behave as follows:
$$\mbox{\rm either}\qquad
\begin{array}[c]{c|ccc}
x &  & c & \\ \hline
P_0 & - & 0 & + \\ \hline
P_1 & + & + & + \\
\end{array}\qquad
\hbox{\rm or}\qquad
\begin{array}[c]{c|ccc}
x &  & c & \\ \hline
P_0 & + & 0 & - \\ \hline
P_1 & - & - & - \\
\end{array}\ .$$
In both cases, the contribution to $v_P(x)$ decreases by 1.
\item  If $c$ is a root of $P_i$, $0 < i < K$, we have $P_{i-1}(c) = -P_{i+1}(c)\not=0$. Hence the
contribution of the subsequence
$P_{i-1}(x), P_i(x), P_{i+1}(x)$ to $v_P(x)$ does not change and remains equal to 1.
\end{itemize}
The theorem follows from the preceding remarks.
\end{proof}

\subsection{$P$ with multiple roots}\label{Sturm-mult}

 The proof of Theorem \ref{Sturm} relies on the following properties of the sequence $P_0 ,\ldots,
P_K$:
\begin{enumerate}\shitem
\item $P = P_0$, and $P_K$ is a nonzero constant.
\item If $c$ is a root of $P_0$, the product $P_0 P_1$ is negative on some interval 
$(c-\varepsilon,c)$ and positive on some interval $(c,c+\varepsilon)$.
\item If $c$ is a root of $P_i$, $0 < i < K$, then $P_{i-1}(c)  P_{i+1}(c) < 0$.
\end{enumerate}
Assume now that $P$ has multiple roots. We construct, as above, the sequence $P_0 = P,
P_1 = P',\ldots,P_K$. Now $P_K$ is no longer a constant, but the gcd of $P$ and $P'$.
Consider the sequence  $$P_0/P_K ,\ P_1/P_K ,\ldots,\ P_{K-1}/P_K ,\ 1\;.$$ 
This sequence satisfies
 properties 1-2-3 above for the polynomial $P_0/P_K$, which has the same roots as $P$ (not counting
multiplicities). Moreover, if
$a$ is not a root of $P$, the number $v_P(a)$ of sign changes in the sequence $P_0(a),\ldots,P_K(a)$
is obviously the same as the number of sign changes in the sequence
$$P_0(a)/P_K(a),\ldots,\ P_{K-1}(a)/P_K(a),\ 1\;.$$ 
It follows:

\begin{thm}\label{multSturm}
Sturm's theorem still holds if $P$ has multiple roots. The difference $v_P(a) - v_P(b)$ is equal to
the number of {\em distinct} roots of $P$ in the interval $(a,b)$.
\end{thm}

\subsection{A bound for the roots}

\begin{prop}\label{borne}
Let $P = a_0 X^d +\cdots+ a_{d-1} X + a_d$, where $a_0\not = 0$. If $c\in{\Bbb C}$ is a root of $P$,
then
$${\vert}c{\vert} \leq \max_{i=1,\ldots,d}\left(d\,\left\vert {a_i\over a_0}\right\vert
\right)^{1/i}\;.$$
\end{prop}
\begin{proof}
Set $M=\max_{i=1,\ldots,d}\left(d\,\left\vert {a_i/a_0}\right\vert
\right)^{1/i}$ and let $z\in {\Bbb C}$ be such that $|z|>M$. Then $|a_i|<|a_0|\,|z|^i/d$, for
$i=1,\ldots,d$. Hence,
$$|a_1z^{d-1}+\cdots+a_d| \leq |a_1|\,|z|^{d-1}+\cdots+|a_d| < |a_0 z^d|$$
and $P(z)\not=0$.
\end{proof}

Set 
$$M = \max_{i=1,\ldots,d}\left(d\,\left\vert {a_i\over a_0}\right\vert
\right)^{1/i}\;.$$
Then $v_P(x)$ is constant on $(-\infty,-M)$ (resp.  $(M,+\infty)$) and equal to
$v_P(-\infty)$ (resp.
$v_P(+\infty)$), which is the number of sign changes in the sequence of leading coefficients of
$P_0(-X)$, $P_1(-X)$, \ldots, $ P_K(-X)$ (resp.
$P_0(X) $, \ldots, $ P_K(X)$).

\begin{prop}
The total number of distinct real roots of $P$ is $$v_P(-\infty) - v_P(+\infty).$$
\end{prop}

\section{Real roots satisfying inequalities}\label{racineg}

In this section, $P \in {\Bbb R}[X]$ is a non constant polynomial in one variable, and $Q$,
$Q_1$, \ldots, $Q_\ell$ are polynomials in ${\Bbb R}[X]$.

\subsection{One inequality}

We want to count the number of real roots $c$ of $P$ such that $Q(c)>0$.	We modify the construction
of the Sturm sequence by taking
$P_0 = P$,
$P_1 = P'Q$ and, as before, $P_{i+1}=$ the negative of the remainder of the euclidean division of
$P_{i-1}$ by
$P_i$, for
$i > 0$. We stop just before we obtain 
$0$, i.e.\ we stop with $P_K$ which is the gcd of $P$ and $P'Q$. The sequence of polynomials we
obtain in this way is called the Sturm sequence\index{Sturm sequence} of $P$ and
$P'Q$. If the real number $a$ is not a root of $P$, we denote by $v_{P,Q}(a)$ the number of sign
changes in the sequence
$P_0(a), P_1(a) ,\ldots, P_K(a)$.

\begin{thm}\label{Sturmineq}
Let $a < b$ be real numbers which are not roots of $P$. Then $v_{P,Q}(a) - v_{P,Q}(b)$
is equal to the number of distinct roots $c$ of $P$ in $(a,b)$ such that $Q(c) > 0$ minus the number
of those such that $Q(c)<0$.
\end{thm}

\begin{proof}
First consider the case where $P$ and $P'Q$ are relatively prime ($P_K$ is a nonzero constant).
This means that 
$P$ has no multiple root, and no common root with $Q$. The property 2 of Section
\ref{Sturm-mult} is replaced with:
\begin{enumerate}
\item[2'] If $c$ is a root of $P_0$, the product  $P_0 P_1 Q$  is negative on some interval
$(c-\varepsilon,c)$ and positive on some interval $(c,c+\varepsilon)$.
\end{enumerate}
The theorem follows from properties 1-2'-3.\par
If $P_K$ is not a constant, the sequence  $P_0/P_K,P_1/P_K,\ldots,P_{K-1}/P_K,1$ satisfies
properties 1-2'-3 for $P_0/P_K$. Hence, the difference between the numbers of sign changes in this
sequence evaluated at $a$ and $b$, respectively, is equal to the number we want to calculate, and it
coincides with
$v_{P,Q}(a) -v _{P,Q}(b)$.
\end{proof}

\noindent{\bf Remark.} $v_{P,Q^2}(a) - v_{P,Q^2}(b)$ counts the number of distinct roots of
$P$ in
$(a, b)$ which are not real roots of $Q$. Therefore the number of distinct roots $c$ of $P$ in
$(a,b)$ such that $Q(c)>0$ is equal to
$${1\over 2}\left(v_{P,Q}(a) + v_{P,Q^2}(a) - v_{P,Q}(b) - v_{P,Q^2}(b)\right)\;.$$
We can replace $v_{P,Q^2}$ with $v_P$ if $P$ and $Q$ are relatively prime.

\debexo
The {\em Cauchy index\index{Cauchy index}} of a rational fraction $F\in{\Bbb R}(X)$ between $a$ and
$b$ is the number of poles
$c$, with $a<c<b$, such that $\lim_{x\to c_-}F(x)=-\infty$ and $\lim_{x\to c_+}F(x)=+\infty$, minus
the number of those such that
$\lim_{x\to c_-}F(x)=+\infty$ and $\lim_{x\to c_+}F(x)=-\infty$.\\
Given two polynomials $P$ and $Q$ in ${\Bbb R}[X]$, define the Sturm sequence\index{Sturm sequence}
of
$P$ and
$Q$ by taking $P_0=P$, $P_1=Q$ and the rest constructed as above. If $a$ and $b$ are not roots of
$P$, show that the Cauchy index of $Q/P$ between $a$ and $b$ is equal to the difference $v(a)-v(b)$
between the numbers of sign changes in the Sturm sequence evaluated at $a$ and $b$, respectively.
Recover Theorems \ref{Sturm} and \ref{Sturmineq} from this result.
\finexo

\debexo
Show that we get the same result as Theorem \ref{Sturmineq} if we replace $P'Q$ with the
remainder of the euclidean division of $P'Q$ by $P$ in the construction of the Sturm sequence.
\finexo 

\subsection{Several inequalities}

We want to count the number of real roots of a system $P = 0,\  Q_1 > 0 ,...,\ Q_\ell
> 0$.	First assume that $P$ is relatively prime with all $Q_i$. Let
$\varepsilon = (\varepsilon_1 ,\ldots,
\varepsilon_\ell)\in \{0,1\}^\ell$ and $Q^\varepsilon =
Q_1^{\varepsilon_1}\cdots Q_\ell^{\varepsilon_\ell}$. By Theorem \ref{Sturmineq}, $s_\varepsilon = 
v_{P,Q^\varepsilon}(-\infty) - v_{P,Q^\varepsilon}(+\infty)$ is equal to the number of distinct real
roots $c$ of $P$ such that $Q^\varepsilon(c)>0$ minus the number of those such that
$Q^\varepsilon(c)<0$. If $\varphi=(\varphi_1 ,\ldots, \varphi_\ell)\in\{0,1\}^\ell$, we denote by
$c_\varphi$ the number of distinct real roots $c$ of $P$ such that the sign of $Q_i(c)$ is
$(-1)^{\varphi_i}$, for $i= 1,\ldots,\ell$. Let $s$ (resp. $c$) be the vector whose coordinates are
all $s_\varepsilon$ (resp. $c_\varphi$).

\begin{lem}
There is an invertible $2^\ell\times2^\ell$ matrix $A_\ell$, depending only on $\ell$, such that $s
= A_\ell . c$.
\end{lem}

\begin{proof}
We proceed by induction on $\ell$. For $\ell=0$, we have trivially
$s_{\emptyset} = c_{\emptyset}$. For $\ell = 1$, Sturm's theorem and Theorem \ref{Sturmineq} imply
$$\left(\begin{array}{c}s_0\\s_1\end{array}\right)=
\left(\begin{array}{cc}1&1\\1&-1\end{array}\right)\left(\begin{array}{c}c_0\\c_1\end{array}\right).$$
The induction step from $\ell$ to $\ell+1$ is as follows:
$$\left(\begin{array}{c}\vdots\\s_{\varepsilon,0}\\ \vdots\\s_{\varepsilon,1}\\ \vdots\end{array}\right)=
\left(\begin{array}{cc}A_\ell&A_\ell\\A_\ell&-A_\ell\end{array}\right)
\left(\begin{array}{c} \vdots\\c_{\varphi,0}\\ \vdots\\c_{\varphi,1}\\ \vdots\end{array}\right)\;.$$
(Exercise: check this equality). The matrix $ A_{\ell+1} =
\left(\begin{array}{cc}A_\ell&A_\ell\\A_\ell&-A_\ell\end{array}\right)$ has inverse  
$${1\over 2}\left(\begin{array}{cc}A_\ell^{-1}&A_\ell^{-1}\\A_\ell^{-1}&-A_\ell^{-1}\end{array}\right).$$
\end{proof}

We obtain from the lemma $c = A_\ell^{-1}.s$. Since we can compute $s$, we get $c$ and, in
particular, the number of solutions of
$P = 0, Q_1 > 0 ,..., Q_\ell > 0$. In the general case ($P$ can have multiple roots or common
roots with $Q_i$), we replace $Q^\varepsilon$ with 
$\left(\prod_{i=1}^\ell Q_i^2\right)/Q^\varepsilon$, in order to get rid of the roots of $P$ which
are also roots of some $Q_i$. \par
\vskip12pt

\noindent{\bf Remark.} The number of Sturm sequences\index{Sturm sequence} to be computed is
$2^\ell$, which grows exponentially with $\ell$. Actually, it is possible to avoid this exponential
growth:  the idea for this (due to Ben-Or, Kozen and Reif) is that, since the total number of real
roots of
$P$ is at most equal to
$d=\deg(P)$,  at most
$d$ among the $c_\varphi$ are nonzero (see \cite{Madrid}).

\subsection{Deciding the existence of a solution of a system of polynomial
equations and inequalities}\label{decision}
	
The preceding result allows us to decide the existence of a solution of a system
$P = 0$, $Q_1 > 0$, \ldots, $Q_\ell > 0$, where $P$ is nonconstant. If the system contains several
equations
$P_1 = 0$, \ldots, $P_m = 0$, we can replace them with one equation
$P_1^2 +\cdots+ P_m^2 = 0$. We can replace a nonstrict inequality $Q \geq 0$ with the disjunction 
$Q > 0$ or $Q = 0$. It remains the case where the system consists only of strict inequalities:  $Q_1
> 0 ,\ldots, Q_\ell > 0$. In this case:
\begin{itemize}
\item the system is satisfied on some unbounded interval of the form
$(a, +\infty)$ (resp.
$(-\infty, a)$) if and only if the leading coefficients of $Q_1$, \dots, $Q_\ell$ (resp. $Q_1(-X)$,
\dots, $Q_\ell(-X))$ are all positive;
\item the system is satisfied on an interval $(a,b)$, where $a$ and $b$ are real roots
of the product
$Q =\prod_{i=1}^\ell Q_j$, if and only if the system $Q' = 0$, $Q_1 > 0$, \ldots, $Q_\ell > 0$ has a
real solution. Note that this case happens only if $\deg(Q)\geq 2$, and then the derivative
$Q'$ is nonconstant. 
\end{itemize}

\section[Systems with parameters]{Systems of polynomial equations and inequalities with parameters}

We consider a system of polynomial equations and inequalities
$${\cal S}(T,
X)\ : \quad\left\{\begin{array}{l} S_1(T,X)\triangleright_1 0\\
S_2(T,X)\triangleright_2 0\\
\cdots\\
S_\ell(T,X)\triangleright_\ell 0
\end{array}\right.\;,$$
where the  $S_i$ are real polynomials in $T=(T_1,\ldots,T_p)$ and $X$, and $\triangleright_i$ are
either
$=$ or $\not=$ or
$>$ or
$\geq$. We consider $X$ as the variable and $T$ as parameters. In this section we explain how to
discuss the existence of a real solution of this system, depending on the real parameters $T$.

\subsection{Tarski-Seidenberg}

We shall prove in this section the following result.

\begin{thm}[Tarski-Seidenberg -- first form]\label{TS1}\index{Tarski-Seidenberg theorem} There exists
an algorithm which, given a system of polynomial equations and inequalities in the variables
$T=(T_1,\ldots,T_p)$ and $X$ with coefficients in ${\Bbb R}$
$${\cal S}(T, X)\ : \quad\left\{\begin{array}{l} S_1(T,X)\triangleright_1 0\\
S_2(T,X)\triangleright_2 0\\
\cdots\\ S_\ell(T,X)\triangleright_\ell 0
\end{array}\right.$$
 (where the $\triangleright_i$ are either $=$ or $\not=$ or $>$ or $\geq$), produces a finite list
${\cal C}_1(T)$, \dots, ${\cal C}_k(T)$ of systems of polynomial equations and inequalities in $T$
with coefficients in ${\Bbb R}$ such that, for every $t\in {\Bbb R}^p$, the system ${\cal S}(t,X)$
has a real solution if and only if one of the ${\cal C}_j(t)$ is satisfied.
\end{thm}

In other words, the formula ``$\exists X\ {\cal S}(T,X)$'' is equivalent to the disjunction ``${\cal
C}_1(X)$ or \dots\ or ${\cal C}_k(X)$''. The Tarski-Seidenberg theorem means that there is an
algorithm for eliminating the real variable $X$.  A well known example of elimination of a real
variable is
$$\displaylines
{\qquad\exists X\ AX^2+BX+C=0\quad \Leftrightarrow\hspace*{\fill}\cr
\hspace*{\fill} (A\not=0\ \mbox{\rm and}\
B^2-4AC\geq0)\ 
\mbox{\rm or}\ (A=0\ \mbox{\rm and}\ B\not=0)\ \mbox{\rm or}\ (A=B=C=0)\;.\cr}$$
This example shows that, in order to discuss a system depending on parameters, it is convenient to
fix the degrees (with respect to $X$) of the polynomials in the system. If $S\in {\Bbb R}[T,X]$,
denote by $\lc(S)$ the leading coefficient (in ${\Bbb R}[T]$) of
$S$, considered as a polynomial in $X$. We shall call {\em system with fixed degrees} a system of
the form $({\cal S}(T,X), {\cal D}(T))$, such that ${\cal D}(T)$ contains, or implies, the
inequations
$\lc(S_i)\not=0$ for all polynomials appearing in ${\cal S}(T,X)$.  Observe that any system is
equivalent to a finite disjunction of systems with fixed degrees. For instance, the inequality
$TX^3 + (U-1)X^2 +V >0$ is equivalent to
$$\displaylines{
\qquad(T\not=0 \et TX^3 + (U-1)X^2 +V >0) \ou \hspace*{\fill}\cr
\hspace*{\fill}(T=0\ \et U-1\not=0 \et (U-1)X^2+V>0) \ou \hspace*{\fill}\cr
 \hspace*{\fill}(T=0 \et U-1=0 \et V>0)\;.\qquad\cr}$$
Hence, it is sufficient to discuss systems with fixed degrees. 
Observe also that we can consider only equations $=0$ and strict inequalities $>0$, since the other
cases $\not=0$ and $\geq0$ are equivalent to disjunctions of the preceding ones.

\noindent{\bf Remark.} The proof of the Tarski-Seidenberg theorem\index{Tarski-Seidenberg theorem}
will show the following important fact: if all polynomials in ${\cal S}(T, X)$ have coefficients in
${\Bbb Q}$, the algorithm produces systems ${\cal C}_1(T),\ldots,{\cal C}_k(T)$ where all
polynomials have coefficients in
${\Bbb Q}$.

\subsection{Systems with one equation}\label{1eq}

First we consider the particular case of a system with fixed degrees containing one equation of
positive degree with respect to $X$. It is convenient to introduce the function  $\sign
:{\Bbb R}
\to
\{-1,0,1\}$ defined by
$$ \sign(r)=\cases{
1&if $r>0\;,$\cr
0&if $r=0\;,$\cr
-1&if $r<0\;.$\cr}
$$
Let $P, Q_1,\ldots, Q_\ell$ be real polynomials in
$T=(T_1,\ldots, T_p)$ and $X$, of positive degrees with respect to $X$. Let
${\cal D}(T)$ be the system  $\lc(P)\not=0$
and $\lc(Q_i)\not=0$, $i=1,\ldots,\ell$.

\begin{lem}
There is an algorithm which, given $(P,Q_1,\ldots,Q_\ell)$, produces a finite list $R_1$,\ldots,
$R_k$ of polynomials in $T$ and a function $c : \{-1,0,1\}^k\to {\Bbb N}$ such that, for every
$\varepsilon=(\varepsilon_1,\ldots,\varepsilon_k)\in \{-1,0,1\}^k$ and every $t\in{\Bbb R}^p$ which
satisfies
$${\cal D}(t)\et \sign(R_1(t))=\varepsilon_1\et\ldots\et \sign(R_k(t))=\varepsilon_k\;,$$
the system 
$$ P(t,X)=0\et Q_1(t,X)>0 \et\ldots\et Q_\ell(t,X)>0 $$
has exactly $c(\varepsilon)$ solutions.
\end{lem}

\begin{proof}
We perform the computations of subsection~\ref{racineg}, i.e.\
 we compute Sturm
sequences. For every new polynomial obtained in a Sturm sequence, we test
whether its leading coefficient is zero or nonzero. In the case where the leading coefficient is
zero, we replace the polynomial with its truncation. We do not test the leading coefficients of
polynomials
$P,Q_1,\ldots,Q_\ell$, since we assume they are all nonzero; this ensures, in particular, that all
Sturm sequences start with nonconstant polynomials.\par
 In this way we obtain a tree of computation of Sturm sequences; the
branching tests are polynomial equations ($= 0$) and inequations ($\not= 0$) in the parameters
$T$. Every branch of the computation tree gives a system of polynomial equations and inequations in
$T$ and the Sturm sequence corresponding to all parameters $t$ satisfying this system. The signs
($>0$ or
$<0$) of the leading coefficients of the polynomials in this Sturm sequence determine the
difference $v(-\infty) - v(+\infty)$ between the numbers of sign changes.\par
The leading coefficients are rational fractions $A(T)/B(T)$ in $T$, where $B$ is assumed to be
nonzero in the branch. Note that the sign of $A(t)/B(t)$ is the same as the sign of $A(t)B(t)$. We
take for
$R_1,\ldots,R_k$ the  $A(T)B(T)$, for all leading coefficients $A(T)/B(T)$ of
 polynomials appearing in all branches of trees of computation of Sturm sequences. If we fix the
sign ($-1$, $0$ or $1$) of each $R_1(t),\ldots,R_k(t)$ and assume ${\cal D}(t)$ holds, the results
of~\ref{racineg} give us the number of real solutions of the system
$$ P(t,X)=0\et Q_1(t,X)>0 \et\ldots\et Q_\ell(t,X)>0\;. $$ 
\end{proof}

\subsection{Example: polynomial of degree $4$}\label{quatre}

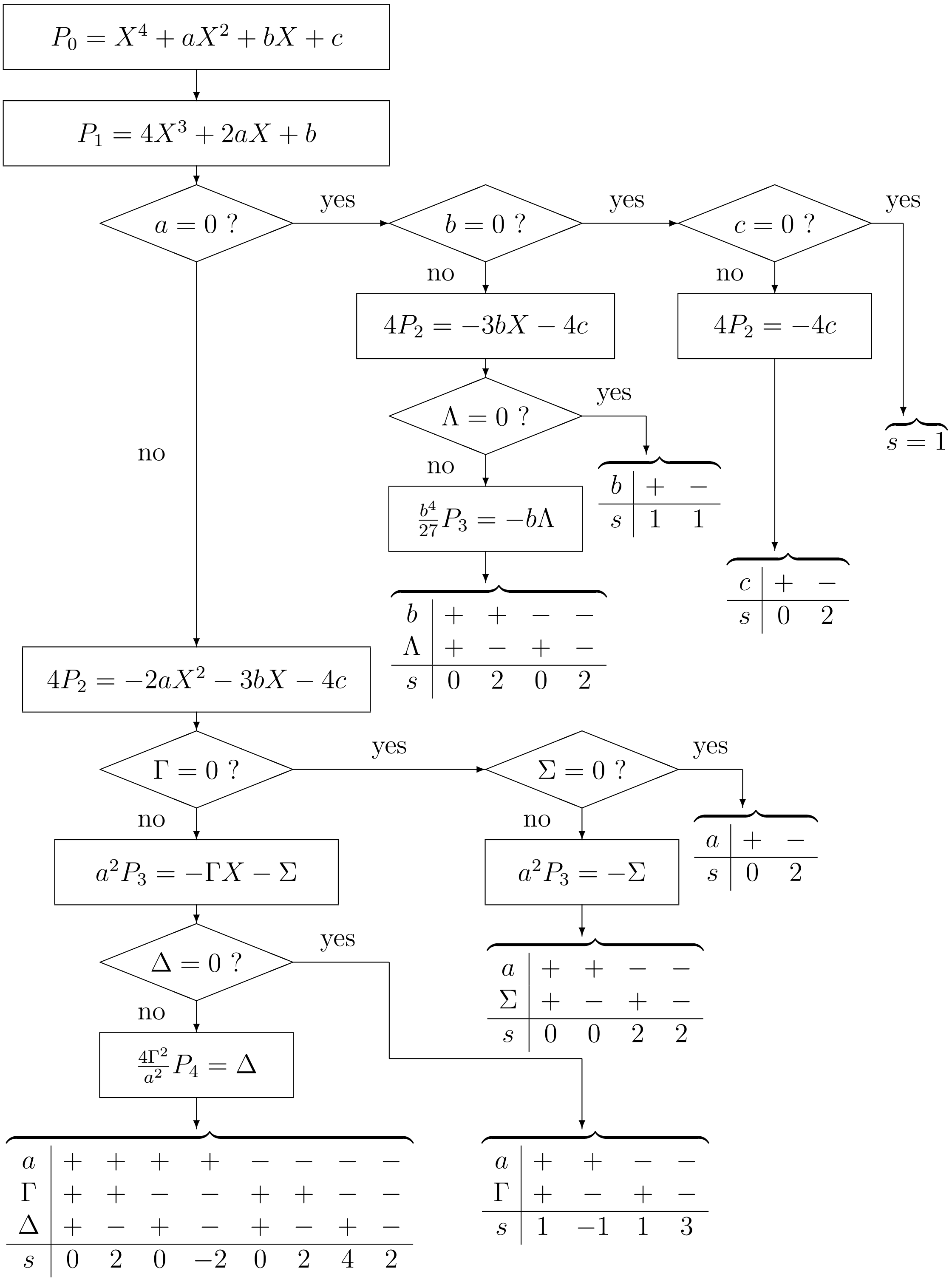
\begin{figure}
\setlength{\unitlength}{1mm}

\noindent
\begin{picture}(150,195)
\put(-5,0){
\begin{picture}(150,200)
\put(0,185){\framebox(60,10){$P_0=X^4+aX^2+bX+c$}}
				\put(30,185){\vector(0,-1){5}}
\put(0,170){\framebox(60,10){$P_1=4X^3+2aX+b$}}
				\put(30,170){\vector(0,-1){3}}
\put(15,155){\test{$a=0$ ?}}
				\put(30,155){\vector(0,-1){60}}
				\put(23,125){\makebox(0,0){no}}
\put(3,85){\framebox(54,10){$4 P_2=-2aX^2-3bX-4c$}}
				\put(30,85){\vector(0,-1){3}}
\put(15,70){\test{$\Gamma=0$ ?}}
				\put(30,70){\vector(0,-1){5}}
				\put(23,68){\makebox(0,0){no}}
\put(8,55){\framebox(44,10){$a^2 P_3=-\Gamma X-\Sigma$}}
				\put(30,55){\vector(0,-1){3}}
\put(15,40){\test{$\Delta=0$ ?}}
				\put(30,40){\vector(0,-1){5}}
				\put(23,38){\makebox(0,0){no}}
\put(15,25){\framebox(30,10){${4\Gamma^2\over a^2}P_4=\Delta$}}
				\put(30,25){\vector(0,-1){5}}
\put(30,20){
\makebox(0,0)[t]{
$\overbrace{
\begin{array}[t]{c|cccccccc}
a &+&+&+&+&-&-&-&-\\
\Gamma&+&+&-&-&+&+&-&-\\
\Delta&+&-&+&-&+&-&+&-\\
\hline
s&0&2&0&-2&0&2&4&2
\end{array}
}$}}


				\put(45,46){\line(1,0){15}}
				\put(52,49){\makebox(0,0){yes}}
				\put(60,46){\line(0,-1){15}}
				\put(60,31){\line(1,0){30}}
				\put(90,31){\vector(0,-1){11}}
\put(90,20){
\makebox(0,0)[t]{
$\overbrace{
\begin{array}[t]{c|cccc}
a &+&+&-&-\\
\Gamma&+&-&+&-\\
\hline
s&1&-1&1&3
\end{array}
}$}}


				\put(45,76){\vector(1,0){30}}
				\put(60,79){\makebox(0,0){yes}}
\put(75,70){\test{$\Sigma=0$ ?}}
				\put(90,70){\vector(0,-1){5}}
				\put(83,68){\makebox(0,0){no}}
\put(75,55){\framebox(30,10){$a^2 P_3=-\Sigma$}}
				\put(90,55){\vector(0,-1){5}}
\put(90,50){
\makebox(0,0)[t]{
$\overbrace{
\begin{array}[t]{c|cccc}
a &+&+&-&-\\
\Sigma&+&-&+&-\\
\hline
s&0&0&2&2
\end{array}
}$}}

				\put(105,76){\line(1,0){10}}
				\put(110,79){\makebox(0,0){yes}}
				\put(115,76){\vector(0,-1){6}}
\put(115,70){
\makebox(0,0)[t]{
$\overbrace{
\begin{array}[t]{c|cc}
a &+&-\\
\hline
s&0&2
\end{array}
}$}}


				\put(45,161){\vector(1,0){15}}
				\put(52,164){\makebox(0,0){yes}}
\put(60,155){\test{$b=0$ ?}}
				\put(75,155){\vector(0,-1){5}}
				\put(68,153){\makebox(0,0){no}}
\put(55,140){\framebox(40,10){$4P_2=-3bX-4c$}}
				\put(75,140){\vector(0,-1){3}}
\put(60,125){\test{$\Lambda=0$ ?}}
				\put(75,125){\vector(0,-1){5}}
				\put(68,123){\makebox(0,0){no}}
\put(60,110){\framebox(30,10){${b^4\over 27}P_3=-b\Lambda$}}
				\put(75,110){\vector(0,-1){5}}
\put(75,105){
\makebox(0,0)[t]{
$\overbrace{
\begin{array}[t]{c|cccc}
b &+&+&-&-\\
\Lambda&+&-&+&-\\
\hline
s&0&2&0&2
\end{array}
}$}}

				\put(90,131){\line(1,0){10}}
				\put(95,134){\makebox(0,0){yes}}
				\put(100,131){\vector(0,-1){6}}
\put(100,125){
\makebox(0,0)[t]{
$\overbrace{
\begin{array}[t]{c|cc}
b &+&-\\
\hline
s&1&1
\end{array}
}$}}

				\put(90,161){\vector(1,0){15}}
				\put(97,164){\makebox(0,0){yes}}
\put(105,155){\test{$c=0$ ?}}
				\put(120,155){\vector(0,-1){5}}
				\put(113,153){\makebox(0,0){no}}
\put(105,140){\framebox(30,10){$4P_2=-4c$}}
				\put(120,140){\vector(0,-1){30}}
\put(120,110){
\makebox(0,0)[t]{
$\overbrace{
\begin{array}[t]{c|cc}
c &+&-\\
\hline
s&0&2
\end{array}
}$}}

				\put(135,161){\line(1,0){5}}
				\put(140,164){\makebox(0,0){yes}}
				\put(140,161){\vector(0,-1){30}}
\put(140,131){
\makebox(0,0)[t]{
$\overbrace{
s=1
}$}}

\end{picture}
}
\end{picture}
\caption{The computation tree of the Sturm sequence \label{eq4}}
\end{figure}

In order to make the preceding lemma clearer, we treat the following example: the equation
$X^4 + aX^2 + bX + c = 0$, where $a, b, c$ are parameters. The leading coefficient of the equation
does not vanish. Since there is no inequality, we have only one Sturm sequence to compute. The
figure \ref{eq4} shows the computation tree of this Sturm sequence. In this tree the following
polynomial expressions in the parameters appear: 
\begin{eqnarray*}
\Gamma&=&2a^3 -8ac + 9b^2,\\
\Delta&=&16a^4c - 4a^3b^2 - 128a^2c^2 + 144ab^2c - 27b^4 + 256c^3,\\
\Sigma&=& b(a^2 + 12c),\\
\Lambda&=&-27b^4 + 256c^3.
\end{eqnarray*}

We have cleared the denominators in the presentation of the results of the computation. Observe
that the properties of Sturm sequences are not altered if we multiply a polynomial in this
sequence by a positive quantity (for instance, the square of a nonzero quantity). Remark that we
find cases where  $s = -2$ or $-1$; of course, these cases can never happen, since the number of
roots has to be nonnegative.

We can draw the list of cases where $s = 0$. This gives a necessary and sufficient condition for
the polynomial to have no real root. We can group some of the cases by noting that
$\Delta$ equals $\Lambda$ if $a=0$, and $256c^3$ (which has the same sign as $c$) if $a=b=0$.
Finally, we find that the polynomial
$X^4 + aX^2 + bX + c$ has no real root if and only if
$$\displaylines{
\qquad(a\geq 0\et\Delta>0)\ou (a>0 \et \Gamma=0)\hspace*{\fill}\cr
\hspace*{\fill} \ou (a<0\et\Gamma>0\et\Delta>0)\;.\qquad\cr}$$

The situation becomes clearer when we consider the picture in the space of parameters $(a, b, c)$
(see Figure \ref{disc}). The polynomial $P$ has a multiple root if and only if
$\Delta=0$ ($\Delta$ is the {\em discriminant\index{discriminant}} of $P$). If
$P$ has a multiple root, we denote by
$\alpha, \alpha, \beta, -(2\alpha+\beta)$ its four roots. It follows
\begin{eqnarray*}
a&=& -3\alpha^2 - 2\alpha\beta - \beta^2\;,\\
b&=&2\alpha(\alpha^2 + 2\alpha\beta + \beta^2) = 2\alpha(-2\alpha^2 - a)\; ,\\
 c&=&-\alpha^2(\beta^2 + 2\alpha\beta) =-\alpha^2(-3\alpha^2 - a)\;.
\end{eqnarray*}
 We obtain in this way parametrizations of the curves $\Delta=0$ in the planes $a = \mbox{\rm
constant}$. Note that, for $a > 0$, the two imaginary roots of
$\alpha^2 = - a/2$  give an isolated real point of the curve
$\Delta = 0$ in the corresponding plane (in this case, $P$ has two conjugate complex double roots).
These points form the half-branch of parabola that we see on Figure \ref{disc}; this parabola is
part of the surface $\Delta=0$.

\begin{figure}[ht]
\centering
\dessin 120mm by 96mm (qaeps)
\caption{The discriminant\index{discriminant} of the polynomial of degree four\label{disc}}
\end{figure}

When one moves continuously from one point of the space of parameters 
$(a,b,c)$ to another without meeting the surface $\Delta=0$, the implicit function theorem implies
that the number of real roots does not change. Hence, the number of real roots is constant in each
connected component of the complement of $\Delta=0$. 

\debexo
Identify on Figure \ref{disc} the set where $s = 0$ (no real root), and compare with the conditions
we have found.
\finexo

\debexo
Draw the semialgebraic subset
$$ S=\{ (a,c)\in{\Bbb R}^2\;;\ X^4+aX^2+c=0\ \mbox{\rm has no real root}\}\;.$$
We want to show that $S$ cannot be described by a conjunction of polynomial inequalities
$$P_1(a,c)>0 \et\ldots\et P_s(a,c)>0\;.$$
Suppose that $S$ admits such a description. Show that one of the $P_i$ vanishes on the parabola
$c=a^2/4$ and can be written $P_i=\linebreak (a^2-4c)^{2m+1}Q$, where $a^2-4c$ does not divide $Q$.
Deduce that $P_i$ should be negative on a part of $S$ (hint: the sign of $P_i$ changes along the
parabola). Conclude.\\
Show that the semialgebraic set
$$ \{ (a,b,c)\in{\Bbb R}^3\;;\ X^4+aX^2+bX+c=0\ \mbox{\rm has no real root}\}$$
cannot be described by a conjunction of polynomial inequalities.
\finexo

\subsection{General systems}

We now discuss general systems of equations and inequalities with fixed degrees. We repeat the
arguments of \ref{decision}
\par

Firt consider the case of a system with several equations (of positive degree with respect to the
variable $X$). We can replace the equations $P_1=\ldots=P_k=0$ with one equation
$P_1^2+\cdots+P_k^2=0$ and we can proceed as in susbsection \ref{1eq}.\par

Now consider the case where there is no equation of positive degree with respect to the
variable $X$. There are only inequalities $Q_1>0$, \ldots, $Q_\ell>0$ with $Q_j$ of positive
degree with respect to $X$. It can be decided, by looking at the signs of the leading coefficients of
the $Q_j$, whether the system is satisfied on an unbounded interval. The existence of an interval,
whose endpoints are roots of $Q=\prod_{j=1}^\ell Q_j$, in which the system holds can be decided
by discussing the system obtained by the adjunction of the equation $Q'=0$. We then proceed  as in 
susbsection \ref{1eq}.

\section{Another method for counting real roots}

\subsection{Hermite's method}\label{Hermite}\index{Hermite's method}

	Let $P\in{\Bbb R}[X]$ be a polynomial of degree $d$, $\alpha_1,\ldots, \alpha_d$ its roots in
${\Bbb C}$ (counted with multiplicities). If $i$ is a nonnegative integer, set $N_i = \sum_{i=1}^d
\alpha_j^i$. The $N_i$ are called {\em the Newton sums\index{Newton sums} of the roots of $P$}. They
are symmetric polynomials of the roots of $P$, with integer coefficients. Hence, if
$P=a_0X^d+a_1X^{d-1}+\cdots+a_d$, the Newton sums can be expressed as polynomials in  
$a_1/a_0$,\ldots,$a_d/a_0$ with coeficients in ${\Bbb Z}$.

\debexo\label{Newton}
Show that the expansion of the logarithmic derivative $P'/P$ as a series in $1/X$ is
$\sum_{i=0}^\infty N_i(1/X)^{1+i}$. Deduce from this fact the following relations:
\begin{eqnarray*}
N_0&=&d\\
a_0N_1&=&-a_1\\
a_0N_2&=&-(N_1 a_1 +2 a_2)\\
&&\cdots\\
a_0N_i&=&-(N_{i-1} a_1+N_{i-2} a_2+\cdots+N_1 a_{i-1} +  i a_i)\quad\mbox{\rm if } i\leq d\\
&&\cdots\\
a_0N_i&=&-(N_{i-1} a_1+N_{i-2} a_2+\cdots+N_{i-d} a_d)\quad\mbox{\rm if } i> d
\end{eqnarray*}
\finexo

 Let us consider the quadratic form whose matrix is
$${\cal H}(P)=\left(\begin{array}{cccc}
N_0&N_1&\ldots&N_{d-1}\\
N_1&N_2&\ldots&N_{d}\\
\vdots&\vdots&\vdots&\vdots\\
N_{d-1}&N_d&\ldots&N_{2d-2}
\end{array}\right)$$

\begin{thm}
The signature of the quadratic form with matrix ${\cal H}(P)$ is equal to the number of distinct
real roots of $P$. Its rank is equal to the number of distinct roots of $P$ (real and complex).
\end{thm}

Recall that the signature of a quadratic form $Q(U)$ in variables $U=(U_1,\ldots,U_d)$ with
coefficients in ${\Bbb R}$ is computed in the following way: we decompose $Q(U)$ as
$$Q(U)= \sum_{i=1}^p L_i(U)^2 - \sum_{i=p+1}^{p+q}L_i(U)^2\;,$$
where the $L_i$ are independant linear forms with coefficients in ${\Bbb R}$, and the signature of
$Q(U)$ is the difference $p-q$.

 \begin{proof} The quadratic form with matrix ${\cal H}(P)$ in  $d$ variables $U_1,\ldots,U_d$
can be decomposed over ${\Bbb C}$ as the sum
 $$(U_1 + \alpha_1 U_2 +\cdots+ \alpha_1^{d-1} U_d)^2
+\cdots+ (U_1 + \alpha_d U_2 +\cdots+ \alpha_d^{d-1}U_d)^2\; .$$
Indeed, if we denote by $V$ the Vandermonde matrix 
$$V=\left(\begin{array}{cccc} 1&1&\ldots&1\\ \alpha_1&\alpha_2&\ldots&\alpha_{d}\\
\vdots&\vdots&\vdots&\vdots\\ \alpha_1^{d-1}&\alpha_2^{d-1}&\ldots&\alpha_d^{d-1}
\end{array}\right)\;,$$
we have ${\cal H}(P) = V\,({}^tV)$. Remark that the linear forms $L_\alpha(U)=U_1 + \alpha U_2
+\cdots+ \alpha^{d-1} U_d $, for all {\em distinct} roots $\alpha$ of $P$, are linearly
independent. Hence, the rank of $Q(U)$ is equal to the number of distinct complex roots of $P$. If
$\alpha$ and $\overline\alpha$ are conjugate nonreal roots  of $P$, we have
$$\displaylines{
\qquad (L_\alpha)^2+(L_{\overline\alpha})^2=\hspace*{\fill}\cr
\hspace*{\fill} 2 \left((U_1 + \Re(\alpha) U_2 +\cdots+
\Re(\alpha^{d-1})U_d)^2 - (\Im(\alpha) U_2 +\cdots+
\Im(\alpha^{d-1})U_d)^2\right)\;.\qquad\cr} $$
The linear forms on the right-hand side of this equality have coefficients in ${\Bbb R}$. Hence,
the pairs of conjugate nonreal roots contribute for
$0$ to the signature. Each distinct real root contributes for $1$ to the signature. It follows that
the signature is equal to the number of distinct real roots of $P$.  
\end{proof}

\debexo\label{Hermiteineg}
Let $Q$ be a polynomial in ${\Bbb R}[X]$. Replace all $N_i$ with $N'_i=\sum_j Q(\alpha_j)\alpha_j^i$
in
${\cal H}(P)$. Show that the signature of the matrix obtained in this way is equal to the
 number of distinct real roots $c$ of $P$ such that $Q(c)>0$ minus the number of
those such that $Q(c)<0$. Explain how to compute the $N'_i$ by considering the expansion of
$P'Q/P$ as a series of $1/X$.
\finexo

\subsection{Determination of the signature}

Consider the {\em principal minors} of the matrix
${\cal H}(P)$, i.e. the determinants  $\delta_i$ consisting of the first $i$ rows and
columns of the matrix, for $i=1,\ldots,d$. Assume that none of the $\delta_i$ is zero. The 
signature of ${\cal H}(P)$ is equal to $d -{}$twice the number of sign changes in the sequence
$1,\delta_1,\ldots,\delta_d$. This result is known as Jacobi's theorem. We give indications on its
proof in the following exercise.

\debexo\label{Jacobi} Consider a symmetric bilinear form $b$ with matrix $M$ in the basis
$(e_1,\ldots,e_d)$. Assume that none of the principal minors  $\delta_1,\ldots,\delta_d$ of $M$ is
zero. We construct an orthogonal basis
$(\varepsilon_1,\ldots,\varepsilon_d)$ for $b$, of the form \begin{eqnarray*}
\varepsilon_1=e_1,\  
\varepsilon_i=\lambda_{i,1}e_1+\cdots+\lambda_{i,i-1}e_{i-1}+e_i \quad\mbox{\rm for } 2\leq i\leq d.
\end{eqnarray*}
Show that, for $i\geq 2$, the determinant of the linear system 
$b(e_1,\varepsilon_i)=0$,\ldots, $b(e_{i-1},\varepsilon_i)=0$ with unknowns $\lambda_{i,j}$,
$j=1,\ldots,i-1$, is equal to
$\delta_{i-1}$. Deduce from this fact that the orthogonal basis
$(\varepsilon_1,\ldots,\varepsilon_d)$ is uniquely determined. Show that
$b(\varepsilon_i,\varepsilon_i)=\delta_i/\delta_{i-1}$ for $i\geq 2$ (and, of course,
$b(\varepsilon_1,\varepsilon_1)=\delta_1$). Prove Jacobi's theorem.
\finexo

When one or several principal minors vanish, the situation is more complicated. For a general
quadratic form, the principal minors are not sufficient to determine the signature. But the
matrix ${\cal H}(P)$ has a special feature: it is a {\em Hankel matrix\index{Hankel matrix}}, which
means that all coefficients
$a_{ij}$ with $i+j =$ constant are equal. In this case, there is a method, due to Frobenius, to
determine the signature by using only the principal minors. For more details concerning
the theorem of Frobenius, see
\cite{G}, Chapter 10, Theorem 24. Using this result and the fact that the rank of ${\cal H}(P)$ is
equal to $r$ if and only if  $\delta_r\not=0$
and $\delta_i=0$ for $r<i\leq d$ (cf.\ \ref{racinesdist}), we obtain the following rule (see also
\cite{Madrid}):

\begin{prop}\label{comptesign} Let $P\in{\Bbb R}[X]$ be a polynomial of degree $d$, and let
$\delta_1,\ldots,\delta_d$ be the principal minors of the matrix ${\cal H}(P)$. Let $r$ be such that
$\delta_r\not=0$ and $\delta_{r+1}=\ldots=\delta_{d}=0$ (note that $r\geq1$ since $\delta_1=d$). For
$1\leq i\leq r$, we define the ``conventional signs'' $\widetilde{\sign}(\delta_i)$ as follows:
\begin{enumerate}
\item If $\delta_i\not=0$, $\widetilde{\sign}(\delta_i)=\sign(\delta_i)$.
\item If $\delta_i=\delta_{i-1}=\ldots=\delta_{i-j+1}=0$ and $\delta_{i-j}\not=0$, then
$$\widetilde{\sign}(\delta_{i})=(-1)^{j(j-1)/2}\sign(\delta_{i-j})\;.$$
\end{enumerate}
Then the number of distinct real roots of $P$ is equal to $r$ minus twice the number of changes in
the sequence $1,\widetilde{\sign}(\delta_1),\ldots,\widetilde{\sign}(\delta_r)$.
\end{prop}

Continuing the exercise  \ref{Hermiteineg}, we can obtain a similar result for counting the
difference between the number of distinct real roots $c$ of $P$ such that $Q(c)>0$ and the number
of those such that $Q(c)<0$.

Let us return to the example of a polynomial $P=X^4+aX^2+bX+c$ of degree 4. The Newton
sums\index{Newton sums}
$N_i$ are easily obtained by using the formulas of exercise
\ref{Newton}:
\begin{eqnarray*}
N_0&=&4\\
N_1&=&0\\
N_2&=&-2a\\
N_3&=&-3b\\
N_4&=&2a^2-4c\\
N_5&=&5ab\\
N_6&=&-2a^3+6ac+3b^2\;.
\end{eqnarray*}
From this we obtain the principal minors of the matrix ${\cal H}(P)$~:
\begin{eqnarray*}
\delta_1&=&4\\
\delta_2&=&-8a\\
\delta_3&=&-4(2a^3-8ac+9b^2)=-4\Gamma\\
\delta_4&=& 16 a^4  c - 4 a^3  b^2  - 128 a^2  c^2  + 144 a c b^2  - 27 b^4  + 256 c^3
=\Delta\;,
\end{eqnarray*}
where we use the notation of Section \ref{quatre}. \par

\debexo
Using the rule of Proposition \ref{comptesign}, recover a necessary and sufficient condition for
the polynomial $P$ to have no real root.
\finexo

The example of the polynomial of degree 4 shows a great advantage of Hermite's
method\index{Hermite's method} with respect to Sturm's method, in the case where there are
parameters: the computations in Hermite's method are made {\em without branching}. In other words,
the specialization to specific values of the parameters causes no trouble (if the degrees of the
given polynomials are fixed). Actually the specialization problems in the computation of the Sturm
sequence\index{Sturm sequence} can also be avoided by using {\em subresultant
polynomials\index{subresultant!polynomial}}. These specialization problems are treated extensively in
\cite{Madrid,Tapas}; see also \cite{Loo} for a survey on subresultant polynomials.\par

Another nice feature of Hermite's method is that it can be generalized for systems of
polynomial equations in several variables having finitely many solutions. The number of distinct
real solutions can be computed as the signature of a quadratic form over the quotient of the ring of
polynomials by the ideal generated by the equations. See for instance 
\cite{Madrid}.\par 

Hermite's and Sturm's method are actually strongly related. This is obvious in the example, if we
compare the principal minors with the leading coefficients of the polynomials in the main branch of
computation of the Sturm sequence (cf.\ \ref{quatre}). This is not mere coincidence, and it can be
explained by using the theory of subresultant polynomials. We refer to
\cite{Madrid,Tapas} for this explanation. We shall only consider principal subresultant
coefficients\index{subresultant!principal -- coefficient} in the next section. These are the leading
coefficients of the subresultant polynomials, except that we take them to be zero if the
subresultant polynomial has a degree smaller than its expected degree.

\subsection{Principal subresultant coefficients}

Consider two polynomials $P=a_0X^d+\cdots+a_d$ of degree $d$ and $Q=b_0X^e+\cdots+b_e$ of degree
$e$. The {\em resultant\index{resultant}} of $P$ and $Q$ is the determinant of the {\em Sylvester
matrix\index{Sylvester matrix}} of
$P$ and $Q$, which is the square matrix of size
$d+e$ whose rows are the coordinates of $X^{e-1}P,\ldots,XP,P,Q,XQ,\ldots,X^{d-1}Q$, respectively,
in the monomial basis $X^{d+e-1}$,\ldots, $X$, $1$. We draw this matrix in the case where $e=d-1$: 
$$\left(\begin{array}{cccccccccccc}
a_0   &a_1   &a_2   &\ldots&.     &.&.&a_d&0&.&\ldots&0\\
0     &a_0   &a_1   &\ldots&.     &.&.&a_{d-1}&a_d&0&\ldots&0\\
\vdots&\vdots&\vdots&\vdots&\vdots&\vdots&\vdots&\vdots&\vdots&\vdots&\vdots&\vdots\\
0     &.     &.     &\ldots&0     &a_0&a_1&a_2&.&.&\ldots&a_d\\
0     &.     &.     &\ldots&.&0&b_0&b_1&.&.&\ldots&b_e\\
\vdots&\vdots&\vdots&\vdots&\vdots&\vdots&\vdots&\vdots&\vdots&\vdots&\vdots&\vdots\\
0     &b_0   &.     &\ldots&.&.&b_{e-1}&b_e&0&.&\ldots&0\\
b_0   &b_1&.   &\ldots&.    &.&b_e&0&.&.&\ldots&0
\end{array}\right)\;.$$
Note that the Sylvester matrix is usually presented with another disposition of rows. The
nonclassical presentation given here is borrowed from \cite{Tapas} and has better properties with
respect to the signs of the principal subresultant coefficients.\par
 The resultant\index{resultant} is zero if and only
if
$P$ and
$Q$ have a common factor (i.e., their gcd is nonconstant). This is a well-known result (see for
instance
\cite{L} chapter V, section 10). We shall prove a generalization of this theorem.\par

For $0\leq j< \min(d,e)$, we call {\em principal subresultant
coefficient\index{subresultant!principal -- coefficient}} of order
$j$ of $P$ and $Q$ (denoted by $\mbox{\rm PSRC}_j(P,Q)$) the determinant of the square matrix of size
$d+e-2j$ which is obtained from the Sylvester matrix\index{Sylvester matrix} of $P$ and $Q$ by
deleting the first
$j$ rows, the last $j$ rows, the first $j$ columns and the last $j$ columns. The resultant of
$P$ and $Q$ is $\mbox{\rm PSRC}_0(P,Q)$

\begin{prop}\label{CSP} Let $\ell$ be an integer, $0\leq \ell< \min(d,e)$.
The gcd of $P$ et $Q$ has degree $>\ell$ if and only if $\mbox{\rm PSRC}_0(P,Q)=\ldots=\mbox{\rm
PSRC}_\ell(P,Q)=0$. \end{prop}

\begin{proof}
Consider the following problem:
\begin{itemize}
\item[(*)] Do there exist nonzero polynomials $U$ and $V$, with $\deg(U)<e-\ell$ and
$\deg(V)<d-\ell$, such that
$\deg(UP+VQ)<\ell$ ?
\end{itemize} 
The last inequality can be translated to a homogeneous linear system of $d+e-2\ell$ equations in
$d+e-2\ell$ unknowns (the coefficients of $U$ and
$V$). The determinant of this system is $\pm\mbox{\rm PSRC}_\ell(P,Q)$. Hence, problem (*) has
an affirmative answer if and only if
$\mbox{\rm PSRC}_\ell(P,Q)=0$.\par 
Note that
$P$ and $Q$ have a gcd of degree $>\ell$ if and only if there exist nonzero polynomials  $U$ and
$V$ as above, such that $UP+VQ=0$. The case 
$\ell=0$ of the proposition follows from the remarks already made.
Let $\ell>0$ and assume that we have proved that $P$ and $Q$ have a gcd with degree $\geq\ell$
if and only if
$\mbox{\rm PSRC}_0(P,Q)=\ldots=\mbox{\rm PSRC}_{\ell-1}(P,Q)=0$. If the gcd of $P$ and $Q$ has
degree $>\ell$, the problem (*) has an affirmative answer and, therefore,
$\mbox{\rm PSRC}_\ell(P,Q)=0$. Conversely, if $\mbox{\rm PSRC}_0(P,Q)=\ldots=\mbox{\rm
PSRC}_\ell(P,Q)=0$, then (*) has an affirmative answer and there are nonzero polynomials $U$ and
$V$, with $\deg(U)<e-\ell$, $\deg(V)<d-\ell$ and $\deg(UP+VQ)<\ell$. By the inductive assumption,
the gcd of $P$ and $Q$ has degree $\geq\ell$. Since this gcd divides
$UP+VQ$, we have
$UP+VQ=0$. Hence, the gcd of $P$ and $Q$ has degree $>\ell$.
\end{proof}

The preceding proposition will be used in the next chapter. We shall now relate the principal
minors of the Hankel matrix\index{Hankel matrix} ${\cal H}(P)$ with certain principal subresultant
coefficients\index{subresultant!principal -- coefficient}.

\begin{prop} Let $\delta_j$ be the principal minor constructed from the first $j$ rows
and columns of ${\cal H}(P)$. Then 
$$a_0^{2j-1}\delta_j= \hbox{\rm PSRC}_{d-j}(P,P')\;.$$
\end{prop}

\begin{proof}
The expansion  of $P'/P$ as a series in $1/X$ (cf.\ exercise \ref{Newton}) allows one to check that
the matrix of size $2j-1$ whose determinant is $\mbox{\rm PSRC}_{d-j}(P,P')$ is the product of the
matrix
$$\left(\begin{array}{ccccccc}
1     &0     &\ldots&0     &0     &\ldots&0     \\
0     &1     &\ldots&0     &0     &\ldots&0     \\
\vdots&\vdots&\vdots&\vdots&\vdots&\vdots&\vdots\\
0     &0     &\ldots&1     &0     &\ldots&0     \\
0     &0     &\ldots&0     &N_0   &\ldots&N_{j-1}\\
\vdots&\vdots&\vdots&\vdots&\vdots&\vdots&\vdots\\
0   &N_0  &\ldots&N_{j-3}&N_{j-2}&\ldots&N_{2j-3}\\
N_0   &N_1   &\ldots&N_{j-2}&N_{j-1}&\ldots&N_{2j-2}
\end{array}\right)$$
with the matrix
$$
\left(\begin{array}{ccccccccc}
a_0     &a_1   &\ldots&a_{2j-2} \\
0     &a_0   &\ldots&a_{2j-1}\\
\vdots&\vdots&\vdots&\vdots\\
0     &.     &\ldots&a_0\\
\end{array}\right)\;,$$
where $a_\ell=0$ if $\ell>d$. The determinant of the former matrix is
$\delta_j$, and the determinant of the latter is  $a_0^{2j-1}$.
\end{proof}

\begin{cor}\label{racinesdist}
The following properties are equivalent:
\begin{enumerate}\shitem
\item $P$ has $r$ distinct complex roots,
\item The matrix ${\cal H}(P)$ has rank $r$,
\item $\delta_r\not=0$ and $\delta_j=0$ for $r<j\leq d$,
\item $\mbox{\rm PSRC}_{d-r}(P,P')\not=0$ and $\mbox{\rm PSRC}_\ell(P,P')=0$ for $0\leq \ell<d-r$.
\end{enumerate}
\end{cor}

\begin{cor}
The number of distinct real roots of $P$ depends only on the signs ($>0$,
$<0$ or
$=0$) of the principal subresultant coefficients\index{subresultant!principal -- coefficient} of $P$
and
$P'$.
\end{cor}

The rule to compute the number of distinct real roots follows from Proposition \ref{comptesign}.
Let us recall that the method of principal subresultant coefficients behaves better than Sturm's
method when there are parameters: the computation of principal subresultant coefficients is uniform
with respect to the parameters (as long as the degree of $P$ is fixed), whereas the computation of
the Sturm sequence is not uniform. Nevertheless, the best way to compute the principal subresultant
coefficients is by using a variant of the Sturm sequence\index{Sturm sequence}.

\chapter{Semialgebraic sets} 

\section[Stability properties]{Stability properties of the class of semialgebraic sets}

\subsection{Definition and first examples}

A {\em semialgebraic subset} of ${\Bbb R}^n$ is the subset of $(x_1,\ldots,x_n)$ in ${\Bbb R}^n$
satisfying a boolean combination of polynomial equations and inequalities with real coefficients.
In other words, the semialgebraic subsets of ${\Bbb R}^n$ form the smallest class ${\cal SA}_n$ of
subsets of ${\Bbb R}^n$ such that:
\begin{enumerate}
\item If  $P \in {\Bbb R}[X_1,\ldots, X_n]$, then  $\{x \in {\Bbb R}^n\ ;\ P(x)
= 0\} \in {\cal SA}_n$  and  $\{x \in {\Bbb R}^n\ ;\ P(x) > 0\} \in {\cal SA}_n$.
\item
If $A \in {\cal SA}_n$ and $B \in {\cal SA}_n$, then $A \cup B$, $A \cap B$ and ${\Bbb
R}^n \setminus A$ are in ${\cal SA}_n$.
\end{enumerate}
The fact that a subset of ${\Bbb R}^n$ is semialgebraic  does not depend on the choice of affine
coordinates.

\begin{prop}\label{typesa}
Every semialgebraic subset of ${\Bbb R}^n$ is the union of finitely many semialgebraic subsets of
the form
$$\{ x \in {\Bbb R}^n\ ;\ P(x) = 0  \et  Q_1(x) > 0  \et \ldots
\et Q_\ell(x) > 0 \}\;,$$
where $\ell \in {\Bbb N}$ and $P, Q_1 ,\ldots, Q_\ell \in  {\Bbb R}[X_1 ,\ldots,
X_n]$.
\end{prop}

\begin{proof} Check that the class of finite unions of such subsets satisfies the above properties 1
and 2.
\end{proof}

We give now some examples of semialgebraic sets.
\begin{itemize}
\item The semialgebraic subsets of ${\Bbb R}$ are the unions of finitely many points and open
intervals.
\item An algebraic subset of ${\Bbb R}^n$ (defined by polynomial equations) is semialgebraic.
\item Let $F : {\Bbb R}^m \to {\Bbb R}^n$ be a polynomial mapping: $F = (F_1 ,\ldots, F_n)$, where
$F_i\in  {\Bbb R}[X_1 ,\ldots, X_n]$. Let $A$ be a semialgebraic subset of ${\Bbb R}^n$.
Then
$F^{-1}(A)$ is a semialgebraic subset of ${\Bbb R}^m$. 
\item If $A$ is a semialgebraic subset of ${\Bbb R}^n$ and
$L\subset {\Bbb R}^n$ a line, then $L \cap A$ is the union of finitely many points and open
intervals.
\item If $A \subset {\Bbb R}^m$ and $B \subset {\Bbb R}^n$
are semialgebraic,  $A \times B$ is a semialgebraic subset of ${\Bbb R}^m \times {\Bbb
R}^n$.
\end{itemize}

\debexo
Show that the infinite zigzag
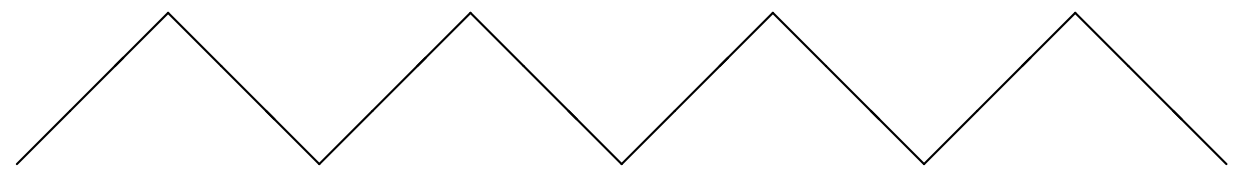
\setlength{\unitlength}{1mm}
$$\begin{picture}(80,10)
\put(0,0){\line(1,1){10}}
\put(10,10){\line(1,-1){10}}
\put(20,0){\line(1,1){10}}
\put(30,10){\line(1,-1){10}}
\put(40,0){\line(1,1){10}}
\put(50,10){\line(1,-1){10}}
\put(60,0){\line(1,1){10}}
\put(70,10){\line(1,-1){10}}
\end{picture}$$
is not semialgebraic. Show that, for every compact semialgebraic subset $K$ of ${\Bbb R}^2$, the
intersection of $K$ with the zigzag is semialgebraic.
\finexo

\subsection{Consequences of Tarski-Seidenberg principle}

We have seen that the class of all semialgebraic subsets is closed under finite unions and
intersections, taking complement, inverse image by a polynomial mapping, cartesian product. It is
also closed under projection.

\begin{thm}[Tarski-Seidenberg -- second form]\label{TS2}\index{Tarski-Seidenberg theorem}
 Let $A$ be a semialgebraic subset of ${\Bbb R}^{n+1}$ and $\pi : {\Bbb R}^{n+1} \to {\Bbb
R}^n$ , the projection on the first $n$ coordinates. Then $\pi(A)$ is a semialgebraic subset of
${\Bbb R}^n$.
\end{thm}

\begin{proof}
Since $A$ is the union of finitely many subsets of the form
$$\{ x=(x_1 ,\ldots, x_{n+1}) \in {\Bbb R}^{n+1}\ ;\ P(x) = 0,\  
Q_1(x) > 0 ,\ldots,\  Q_k(x) > 0 \}\;,$$
we may assume that $A$
itself is of this form. It follows from the Tarski-Seidenberg theorem (first form, \ref{TS1})
that there is a boolean combination\linebreak ${\cal C}(X_1
,\ldots, X_n)$ of polynomial equations and inequalities in $X_1 ,\ldots, X_n$ such that
$$\pi(A)=\{ (x_1
,\ldots, x_n) \in {\Bbb R}^n\ ;\ \exists x_{n+1} \in {\Bbb R}\ (x_1 ,\ldots, x_n ,
x_{n+1}) \in A \}$$
is the set of $(x_1 ,\ldots, x_n)$ which satisfy ${\cal
C}(x_1,\ldots, x_n)$. This means that $\pi(A)$ is semialgebraic.
\end{proof}

We now show some consequences of the Tarski-Seidenberg theorem\index{Tarski-Seidenberg theorem}.

\begin{cor}
\begin{enumerate}\shitem
\item If $A$ is a semialgebraic subset of ${\Bbb
R}^{n+k}$, its image by the projection on the space of the first $n$ coordinates is a semialgebraic
subset of
${\Bbb R}^n$.
\item If $A$ is a semialgebraic subset of ${\Bbb R}^m$ and $F : {\Bbb R}^m \to {\Bbb R}^n$, a
polynomial mapping, then the direct image $F(A)$ is a semialgebraic subset of ${\Bbb
R}^n$.
\end{enumerate}
\end{cor}

\begin{proof} 
The first statement is easily obtained by induction on $k$. For the second statement, note that
$$\{(x, y) \in {\Bbb R}^m \times
{\Bbb R}^n\ ;\ x \in A \et y = F(x)\}$$
is a semialgebraic subset of
${\Bbb R}^m \times {\Bbb R}^n$ and that $F(A)$ is its projection onto ${\Bbb R}^n$.
\end{proof}

\begin{cor}
If
$A$ is a semialgebraic subset of ${\Bbb R}^n$, its closure in
${\Bbb R}^n$ is again semialgebraic.
\end{cor}

\begin{proof}
The closure of $A$ is
$$\clos(A) = \left\{ x
\in {\Bbb R}^n\ ;\ \forall \varepsilon \in {\Bbb R},\;  \varepsilon > 0 \Rightarrow
\exists y \in {\Bbb R}^n,\;  y \in A  \et  \Vert x-y\Vert^2 < \varepsilon^2
\right\}$$ 
and can be written as
$$\clos(A) = {\Bbb R}^n \setminus \left(\pi_1\left(\{(x,\varepsilon)\in {\Bbb R}^n\times{\Bbb
R}\;;\ \varepsilon>0\}
\setminus
\pi_2(B)\right)\right)\;,$$
where
$$B = \left\{ (x, \varepsilon, y) \in {\Bbb R}^n \times {\Bbb R} \times
{\Bbb R}^n\ ;\  y \in A  \et \sum_{i=1}^n (x_i-y_i)^2 < \varepsilon^2)\right\}\;,$$ 
 $\pi_1(x, \varepsilon) = x$ and
 		$\pi_2(x, \varepsilon, y) =  (x, \varepsilon)$. Then observe that $B$ is semialgebraic.
\end{proof}

The example above shows that it is usually boring to write down projections in order to show that a
subset is semialgebraic. We are more used to write down formulas. Let us make precise what is meant
by a {\em first-order formula\index{formula!first order}} (of the language of ordered fields with
parameters in
${\Bbb R}$). A first-order formula is obtained by the following rules.
\begin{enumerate}\shitem
\item If $P \in {\Bbb R}[X_1 ,\ldots, X_n]$, then $P = 0$ and $P > 0$ are first-order formulas.
\item If $\Phi$ and $\Psi$ are first-order formulas, then ``$\Phi \et \Psi$'', ``$\Phi \ou
\Psi$'', ``$\mbox{\rm not }\Phi$'' (often denoted by
$\Phi\wedge\Psi$, $\Phi\vee\Psi$ and  $\neg \Phi$, respectively) are first order formulas.
\item If $\Phi$ is a formula and $X$, a variable ranging over ${\Bbb R}$, then $\exists X \Phi$
and
$\forall X \Phi$ are first-order formulas.
\end{enumerate}

The formulas obtained by using only rules 1 and 2
are called {\em quantifier-free formulas\index{formula!quantifier-free}}. By definition, a subset $A
\subset {\Bbb R}^n$ is semialgebraic if and only if there is a quantifier-free formula $\Phi(X_1
,\ldots, X_n)$ such that
$$(x_1 ,\ldots, x_n)\in A  \Longleftrightarrow  \Phi(x_1 ,\ldots, x_n).$$
The Tarski-Seidenberg theorem has the following useful formulation.

\begin{thm}[Tarski-Seidenberg -- third form]\index{Tarski-Seidenberg theorem} If
$\Phi(X_1,\ldots,X_n)$ is a\linebreak first-order formula, the set of 
$(x_1 ,\ldots, x_n) \in {\Bbb R}^n$ which satisfy $\Phi(x_1 ,\ldots, x_n)$ is semialgebraic.
\end{thm}

\begin{proof}
By induction on the construction of formulas. Rule 1 produces only semialgebraic sets. Rule
2 produces only semialgebrais sets from semialgebraic sets. For rule 3, if
$$\{ (x_1 ,\ldots, x_{n+1}) \in {\Bbb R}^{n+1}\ ;\ \Phi(x_1 ,\ldots, x_{n+1})\}$$
 is semialgebraic, then
$$\{ (x_1 ,\ldots, x_n) \in {\Bbb R}^n\ ;\ \exists x_{n+1}\;
\Phi(x_1 ,\ldots, x_{n+1}) \}$$
is its projection onto ${\Bbb R}^n$ and, hence, it is also semialgebraic. The case of $\forall X
\Phi$ follows by observing that $\forall X \Phi$ is equivalent to $\neg \exists X\neg\Phi$.
\end{proof}

The preceding theorem can be formulated as follows. 
\begin{quote}
\em Every first-order formula is equivalent to a quantifier-free formula,
\end{quote}
or, in other words,
\begin{quote}
\em ${\Bbb R}$ admits the elimination of quantifiers\index{quantifier elimination} in the language of
ordered fields.
\end{quote}
The reader who wants to learn more about model theory and its application to real algebraic
geometry is invited to read the lecture notes \cite{Pr} in the same collection.\par 
\medskip
\noindent{\bf Remark.} One should pay attention to the fact that the quantified variables (or
$n$-tuples of variables) have to range over ${\Bbb R}$, or
${\Bbb R}^n$, or possibly over a {\em semialgebraic} subset of ${\Bbb R}^n$. For instance,
$$\{ (x,y) \in {\Bbb R}^2\ ;\ \exists n \in {\Bbb N}\  
y=nx \}$$
is not semialgebraic (why ?).

\section{Semialgebraic functions}

\subsection{Definition and first properties}

Let $A \subset {\Bbb R}^m$ and	 $B \subset {\Bbb R}^n$ be semialgebraic sets. A mapping $f :
A
\to B$ is said to be {\em semialgebraic} if its graph
$$\Gamma_f = \{(x, y) \in A \times B\ ;\  y = f(x)\}$$
 is a semialgebraic subset of ${\Bbb R}^m \times {\Bbb R}^n$.

For instance:
\begin{itemize}
\item If $f : A \to B$ is a polynomial mapping (all its coordinates are polynomial), it is
semialgebraic.
\item If $f : A \to B$ is a regular rational mapping (all its coordinates are rational
fractions whose denominators do not vanish on
$A$), it is semialgebraic.
\item  If $f : A \to {\Bbb R}$ is a semialgebraic function, then
${\vert} f{\vert}$ is semialgebraic.
\item If $f : A \to {\Bbb R}$ is semialgebraic and $f\geq 0$ on $A$, then $\sqrt{f}$ is a
semialgebraic function.
\end{itemize}
\debexo
 Let $A \subset {\Bbb R}^n$, $A \not =
\emptyset$, be a semialgebraic set. Then the function 
\begin{eqnarray*}
{\Bbb R}^n&\longrightarrow&{\Bbb R}\\
x &\longmapsto& \mbox{\rm dist}(x, A) = \inf \{ \Vert x-y\Vert\ ;\ y \in A \}
\end{eqnarray*}
is continuous semialgebraic.
\finexo

\debexo\label{unbound}
Show that if  $f:(0,1]\to{\Bbb R}$ is a semialgebraic function such that $f(x)$ is not bounded
from above as $x\to 0$, then $\lim_{x\to 0} f(x)=+\infty$. 
\finexo 

Important properties of semialgebraic mappings follow from the 
Tarski-Seidenberg theorem\index{Tarski-Seidenberg theorem}.

\begin{cor}
\begin{enumerate}\shitem
\item The direct image and the inverse image of a semialgebraic set by a semialgebraic
mapping are semialgebraic. For instance, if $P(X_1 ,\ldots,X_n)$ is a polynomial and $f$, a
semialgebraic mapping from $A  \subset {\Bbb R}^m$ to
$B\subset{\Bbb R}^n$, the set $\{ y \in B\ ;\ P(f(y)) > 0 \}$ is semialgebraic.
\item The composition of two semialgebraic mappings is semialgebraic.
\item The semialgebraic functions from $A$ to ${\Bbb R}$ form a ring. \end{enumerate}
\end{cor}

\debexo
Let $U$ be a semialgebraic open subset of ${\Bbb R}^n$, $f:U\to{\Bbb R}$ a semialgebraic function.
Show that, if $f$ admits a partial derivative $\partial f/\partial x_i$ on $U$, then this
derivative is semialgebraic.
\finexo

\subsection{The {\L}ojasiewicz inequality}

The {\L}ojasiewicz inequality\index{Lojasiewicz@{{\L}ojasiewicz inequality}} gives information
concerning the relative rate of growth of two continuous semialgebraic functions. First, we shall
estimate the rate of growth of a semialgebraic function of one variable.

\begin{prop}\label{croispol} Let $f : (A, +\infty ) \to {\Bbb R}$ be a semialgebraic (not necessarily
continuous) function. There exist $B \geq A$ and an integer $N \in {\Bbb N}$ such that
$\vert  f(x)\vert  \leq  x^N $ for all $x \in (B, +\infty)$. \end{prop}

\begin{proof}
Let $\Gamma$ be the graph of $f$. It is a semialgebraic subset of ${\Bbb R}^2$.
By \ref{typesa}, $\Gamma=G_1\cup\ldots\cup G_p$, where each $G_i$ is a nonempty subset of the form
$$G_i=\{(x,y)\in{\Bbb R}^2\ ;\ P_i(x,y)=0,\ Q_{i,1}(x,y)>0,\ldots,\ Q_{i,k_i}(x,y)=0\}\;.$$
All polynomials $P_i$ have degree $>0$ with respect to $y$: otherwise, if
$(x_0,y_0)\in G_i$, $\Gamma$ should contain a nonempty open interval of the vertical line
$\{x_0\}\times{\Bbb R}$, which is impossible since $\Gamma$ is a graph. Let
$$P(x,y)=a_0(x)y^d + a_{1}(x) y^{d-1} +\cdots+ a_d(x)$$ be the product of all $P_i(x,y)$, where
$d>0$ and $a_0\not=0$. Choose
$C\geq A$ big enough so that $a_0(x)$ does not vanish on
$(C, +\infty )$. By Proposition
\ref{borne}, we obtain
$$\vert  f(x)\vert   \leq   \max_{i=1,\ldots,d}\left( d\,\left\vert {a_i(x)\over
a_0(x)}\right|\right)^{1/i}\;.$$ 
As $x$ tends to $+\infty$, the right-hand side of the equality is equivalent to
$\lambda x^\alpha$, where $\lambda >0$ and $\alpha  \in {\Bbb Q}$. Taking $N$ to be a nonnegative
integer $>\alpha$, we obtain 
$B \geq C$, such that $\vert  f(x)\vert  \leq  x^N$ for all $x > B$.
\end{proof}

\begin{thm}[{\L}ojasiewicz inequality]\label{Loja}\index{Lojasiewicz@{{\L}ojasiewicz inequality}}
Let $K \subset  {\Bbb R}^n$ be a compact semialgebraic set, and let $f,g : K \to {\Bbb R}$ be
continuous semialgebraic functions , such that 
$$\forall x \in  K \quad( f(x) = 0 \Rightarrow g(x) =
0 )\;.$$ 
Then there exist an integer $N
\in {\Bbb N}$ and a  constant $C \geq 0$, such that
$$\forall x \in K \quad \vert g(x)\vert ^N \leq  C \vert  f(x)\vert\; .$$
\end{thm}

\begin{proof}
For $t > 0$, set $F_t = \{ x \in K\ ;\ t\vert g(x)\vert  = 1\}$. Since $F_t$ is closed in $K$,
it is compact. Assume $F_t \not= \emptyset$; then $f$ does not vanish on $F_t$ and the continuous
function
$x \mapsto  1/\vert f(x)\vert $  has a maximum on $F_t$, which we denote by
$\theta(t)$. If $F_t = \emptyset$, we set $\theta(t) = 0$. The function $\theta: (0,+\infty)\to
{\Bbb R}$ is semialgebraic (check this fact by writing a formula which
describes its graph). By Proposition \ref{croispol}, there exist
$B > 0$ and
$N
\in {\Bbb N}$ such that
$$\forall t> B\quad  \vert \theta(t)\vert  \leq  t^N\;.$$
This is equivalent to
$$\forall x \in K \quad\left( 0<\vert g(x)\vert  <  {1\over B}   \Rightarrow   {1\over\vert
f(x)\vert } \leq   {1\over\vert g(x)\vert ^N}\right)\;.$$
Let $D$ be the maximum of the continuous function $\vert g(x)\vert
^N / \vert  f(x)\vert$  on the compact set
$$\{ x \in K\ ;\ \vert g(x)\vert  \geq 1/B\} $$
(observe that $f$ does not vanish on this set), and let $C =  \max(1, D)$. We obtain $\vert
g(x)\vert ^N\leq  C\vert f(x)\vert $ for all $x\in K$. 
\end{proof}

\debexo
Let $S$ be a closed semialgebraic subset of the plane which contains the graph of the exponential
function
$y=e^x$. Show that $S$ contains some interval $(-\infty,A)$ of the $x$-axis. Hint: use the function
``distance to $S$''.
\finexo

\section{Decomposition of a semialgebraic set} 

We have seen that a semialgebraic subset of ${\Bbb R}$ can be decomposed as the union of
finitely many points and open intervals. We shall see that every semialgebraic set can be
decomposed as the disjoint union of finitely many pieces which are semialgebraically
homeomorphic to open hypercubes 
$(0,1)^d$ of different dimensions.
A semialgebraic homeomorphism $h : S \to T$ is a bijective continuous semialgebraic mapping from $S$
onto $T$, such that
$h^{-1}: T\to S$ is continuous.

\debexo
Check that  $h^{-1}$ is also semialgebraic.
\finexo

The method of decomposition by using successive codimension $1$ projections
is the main tool for studying semialgebraic sets, and it is used in the foundational paper of S.
{\L}ojasiewicz (1964). We now explain the cylindrical algebraic decomposition of Collins
\cite{Cl}, which makes precise the algorithmic content of this method.

\subsection{Cylindrical algebraic decomposition}

 A {\em cylindrical algebraic decomposition\index{cylindrical algebraic decomposition}}  (abbreviated
to {\em c.a.d.}) of 
${\Bbb R}^n$ is a sequence ${\cal C}_1,\ldots,\ {\cal C}_n$, where, for $1\leq k\leq n$,  ${\cal
C}_k$ is a finite partition of ${\Bbb R}^k$ into semialgebraic subsets (which are called {\em
cells\index{cell}}), satisfying the following properties:
\begin{itemize}
\item[a)] Each cell $C\in {\cal C}_1$ is either a point, or an open interval.
\item[b)] For every $k$, $1\leq k<n$, and for every $C\in {\cal C}_k$, 
there are finitely many continuous semialgebraic functions 
$$\xi_{C,1}<\ldots<\xi_{C,\ell_C}: C\longrightarrow {\Bbb R}\;,$$
and the cylinder $C\times {\Bbb R} \subset  {\Bbb R}^{k+1}$ is the disjoint union of cells of ${\cal
C}_{k+1}$  which are:
\begin{itemize}
\item either the {\em graph} of one of the functions $\xi_{C,j}$, for $j = 1,\ldots,\ell_C$:
$$A_{C,j}=\{ (x', x_{k+1}) \in
C \times {\Bbb R}\ ;\ x_{k+1} = \xi_{C, j}(x') \}\;,$$
\item or a {\em band\index{band}} of the cylinder bounded from below and from above by the graphs of
functions
$\xi_{C,j}$ and $\xi_{C,j+1}$, for $j = 0 ,\ldots, \ell_C$, where we take $\xi_{C,0} = -\infty$
and $\xi_{i,\ell_C+1} = +\infty$:
$$B_{C, j}=\{ (x', x_{k+1}) \in C \times {\Bbb R}\ ;\ 
\xi_{C,j}(x') < x_{k+1} <  \xi_{C, j+1}(x')\}\;.$$
\end{itemize}
\end{itemize}

\begin{prop}\label{paves}
Every cell\index{cell} of a c.a.d.\index{cylindrical algebraic decomposition}\ is semialgebraically
homeomorphic to an open hypercube
$(0,1)^d$ (by convention,
$(0,1)^0$ is a point).
\end{prop}

\begin{proof}
We prove the property of the proposition for cells of ${\cal
C}_k$, by induction on $k$. The key point is to observe that, using the notation above, every
graph $A_{C, j}$ is semialgebraically homeomorphic to $C$ and every band
$B_{C, j}$ is semialgebraically homeomorphic to $C \times (0,1)$. In the case of $A_{C,
j}$, the homeomorphism is simply 
$$C \ni x' \longmapsto  (x',\xi_{C, j} (x')) \in A_{C, j}\;.$$
 For $B_{C, j}$  we take
\begin{eqnarray*}
C\times (0,1) \ni (x', t)&\mapsto&  \left(x',
(1-t) \xi_{C, j} (x') +t \xi_{C, j+1}(x')\right)\mbox{\rm\  if }0 < j < \ell_C\;,\\ 
 &&\left(x',  {t-1\over t}+ \xi_{C, 1}(x')\right)\mbox{\rm\  if }j=0,\ \ell_C\not=0\;,\\
&&\left(x', - {1\over t} + {1\over 1-t}\right)\mbox{\rm\  if }j =\ell_C=0\;,\\
&&\left(x',  {t\over 1-t}+\xi_{C,\ell_C}(x')\right)\mbox{\rm\  if }j=\ell_C\not= 0\;.
\end{eqnarray*}
\end{proof}

\enlargethispage{\baselineskip}
It is time to explain what we want to do with a c.a.d.\index{cylindrical algebraic decomposition}. We
shall use the following terminology: given a finite family $P_1,\ldots ,P_r$ of polynomials in 
${\Bbb R}[X_1,\ldots ,X_n]$, we say that a subset $C$ of ${\Bbb R}^n$ is $(P_1,\ldots ,P_r )$-{\em
invariant\index{invariant@{$(P_1,\ldots,P_r)$-invariant}}} if every polynomial
$P_i$ has a constant sign ($>0$, $<0$, or $=0$) on $C$. We want to construct, from a finite family
$P_1,\ldots ,P_r$  of polynomials in ${\Bbb R}[X_1,\ldots ,X_n]$, a c.a.d.\ of ${\Bbb R}^n$ such
that:
\begin{itemize}
\item[c)] Each cell\index{cell} $C\in {\cal C}_n$ is $(P_1,\ldots ,P_r)$-invariant.
\end{itemize}
A c.a.d. of ${\Bbb R}^n$ satisfying this property will be called 
{\em adapted\index{cylindrical algebraic decomposition!adapted} 
to} $(P_1,\ldots ,P_r)$.
 
  What is a c.a.d.\ adapted to $(P_1,\ldots,P_r)$ good for? 
First, the condition c) shows that every semialgebraic subset of ${\Bbb R}^n$ which is described
by a boolean combination of equations $P_i =0$ and inequalities
$P_j >0$ or $P_j <0$, where $P_i$ and $P_j$ are among $P_1,\ldots ,P_r$, is the union of some
cells of ${\cal C}_n$.  It follows that every semialgebraic set can be decomposed as the disjoint
union of finitely many pieces, each semialgebraically homeomorphic to an open hypercube $(0,1)^d$.
Moreover, the cylindrical arrangement of cells (property b)
allows one to see that every semialgebraic subset of ${\Bbb R}^k$ described by a
formula\index{formula!first order}
$Q_{k+1}x_{k+1}\ldots Q_nx_n \Phi$ , where 
$Q_{k+1},\ldots ,Q_n$ are existential or universal quantifiers and
 $\Phi$, a boolean combination of equations $P_i =0$ and inequalities $P_j >0$ 
or $P_j <0$, is the union of some cells of ${\cal C}_k$. This can be useful, for instance, to
decide whether such a formula is true or false.

\subsection{Construction of an adapted c.a.d.}

The definition of a c.a.d.\index{cylindrical algebraic decomposition}\ shows the importance of the
functions
$\xi_{C,j}$ whose graphs cut the cylinders $C\times{\Bbb R}$. Since we want the bands of a cylinder
contained in
${\Bbb R}^n$ to be
$(P_1,\ldots  ,P_r )$-invariant\index{invariant@{$(P_1,\ldots,P_r)$-invariant}}, the $\xi_{C,j}$ have
to describe the roots of the polynomials
$P_i$, as functions of $(x_1,\ldots ,x_{n-1})\in C$. We give now a
first result in this direction.

\begin{prop}\label{fracines}  Let $P(X_1,\ldots X_n)$ be a polynomial in ${\Bbb 
R}[X_1,\ldots X_n]$. Let $C\subset {\Bbb 
R}^{n-1}$ be a connected semialgebraic subset and $k\leq d$ in ${\Bbb N}$ such that, for every point
$x'=(x_1,\ldots ,x_{n-1})\in C$, the polynomial
$P(x',X_n)$  has degree $d$ and exactly $k$ distinct roots in ${\Bbb C}$. Then there are $\ell\leq k$
continuous semialgebraic functions $\xi_1<\ldots <\xi_\ell : C
\to {\Bbb R}$  such that, for every $x'\in C$, the set of real roots of $P(x',X_n)$ is exactly
$\{\xi_1(x'),\ldots ,\xi_\ell(x')\}$. Moreover, for $i=1,\ldots,\ell$, the multiplicity of the root
$\xi_i(x')$ is constant for $x'\in C$. 
\end{prop}

\begin{proof}  The argument relies on the ``continuity of roots'', in the following form (a proof
is proposed in the next exercise): 
\begin{itemize}
\item[CR] Choose
$a'\in C$ and let
$z_1,\ldots,z_k$ be the distinct roots of $P(a',X_n)$, with multiplicities
$m_1,\ldots,m_k$, respectively. Choose $\varepsilon>0$ so small that
the open disks
$D(z_i,\varepsilon)\subset {\Bbb C}$ with centers $z_i$ and radius $\varepsilon$ are
disjoint. If $b'\in C$ is sufficiently close to $a'$, the polynomial $P(b',X_n)$ has exactly $m_i$
roots, counted with multiplicities, in the disk $D(z_i,\varepsilon)$, for $i=1,\ldots,k$.
\end{itemize} 
Since $P(b',X_n)$ has $k$ distinct roots and $d=m_1+\cdots+m_k$ roots counted with multiplicities, it
follows that each
$D(z_i,\varepsilon)$ contains exactly one root $\zeta_i$ of multiplicity $k$ of $P(b',X_n)$. If
$z_i$ is real, $\zeta_i$ is real (otherwise, its conjugate $\overline\zeta_i$ would be another root
of $P(b',X_n)$ in
$D(z_i,\varepsilon)$). If $z_i$ is nonreal, $\zeta_i$ is nonreal, since $D(z_i,\varepsilon)$ is
disjoint from its image by conjugation. Hence, if $b'\in C$ is sufficiently close to $a'$,
$P(b',X_n)$ has the same number of distinct real roots as $P(a',X_n)$.
Since $C$ is connected, the number of distint real roots of $P(x',X_n)$ is constant for $x'\in C$.
Let $\ell$ be this number. For $1\leq i\leq \ell$, denote by $\xi_i : C \to {\Bbb R}$ the function
which sends
$x'\in C$ to the $i$-th real root (in the increasing order) of $P(x',X_n)$. The argument above, with
$\varepsilon$ as small as we want, shows, moreover, that the functions
$\xi_i$ are continuous. It follows from the connectedness of $C$ that each $\xi_i(x')$ has constant
multiplicity. If 
$C$  is described by the formula
$\Theta(x')$, the graph of $\xi_i$ is described by the formula
$$\displaylines{
\qquad\Theta(x')\et\exists y_1\ldots \exists y_\ell (y_1<\ldots <y_\ell \et\hspace*{\fill}\cr
\hspace*{\fill}P(x',y_1)=0 \et\ldots \et P(x',y_\ell)=0 \et x_n=y_i)\;,\qquad\cr}$$
which shows that $\xi_i$  is semialgebraic.
\end{proof}

\debexo
We identify monic polynomials $X^d+a_1X^{d-1}+\cdots+a_d\in{\Bbb C}[X]$ of degree $d$ with points
$(a_1,\ldots,a_d)\in{\Bbb C}^d$. With this identification, let
\begin{eqnarray*}
\mu : {\Bbb C}^e\times{\Bbb C}^{d-e}&\longrightarrow&{\Bbb C}^d\\
(R,S)&\longmapsto&RS
\end{eqnarray*}
be the mapping defined by the multiplication of monic polynomials.\\
1) Fix $R^0\in{\Bbb C}^e$ and
$S^0\in {\Bbb C}^{d-e}$. Show that the jacobian determinant of $\mu$ at $(R^0,S^0)$ is equal to
$\pm$ the resultant\index{resultant} of $R^0$ and $S^0$.\\
2) Let $Q^0\in{\Bbb C}^d$ and assume that $Q^0=R^0S^0$, where $R^0$ and $S^0$ are relatively prime
monic polynomials of degrees $e$ and $d-e$, respectively. Show that for every $Q$ sufficiently
close to $Q^0$, there is a unique factorization $Q=RS$ with $R$ close to $R^0$ and $S$ close to
$S^0$.\\
3) Assume $Q^0=(X-z_1)^{m_1}\cdots(X-z_k)^{m_k}$, where $z_1,\ldots,z_k$ are the distinct roots
of
$Q^0$. Show that, for every $Q$ close to $Q^0$, there is a unique factorization $Q=R_1\cdots R_k$,
where the $R_i$ are monic polynomials close to $(X-z_i)^{m_i}$.\\
4) Fix $\varepsilon>0$. Show that every monic polynomial sufficently close to $X^m$ has its roots
in  $D(0,\varepsilon)$ (use Proposition \ref{borne}). Deduce that every monic polynomial
sufficiently close to $(X-z)^m$ has its roots in $D(z,\varepsilon)$.\\
5) Let $Q^0$ be a monic polynomial with distinct roots $z_1,\ldots,z_k$ of multiplicities
$m_1$,\ldots, $m_k$, respectively. Choose $\varepsilon>0$ such that all disks $D(z_i,\varepsilon)$
are disjoint. Show that every monic polynomial close to $Q^0$ has exactly $m_i$ roots counted with
multiplicities in $D(z_i,\varepsilon)$, for $i=1,\ldots,k$.\\
6) Prove property CR above (if $P=a_0(x')X_n+\cdots+a_d(x')$, set $Q=P/a_0(x')$).
\finexo

If we have several polynomials $P_i$, we have also to take care that the roots of the different
$P_i$ do not get mixed. 

\begin{prop}
Let $P$ and $Q$ be polynomials of ${\Bbb 
R}[X_1,\ldots X_n]$. Let $C$ be a connected semialgebraic subset of ${\Bbb 
R}^{n-1}$. Assume that the degree and the number of distinct roots of $P(x',X_n)$ (resp.
$Q(x',X_n)$) and the degree of the gcd of $P(x',X_n)$ and $Q(x',X_n)$ are constant for all $x'\in
C$. Let
$\xi,\zeta : C\to{\Bbb R}$ be continuous semialgebraic functions such that $P(x',\xi(x'))=0$ and
$Q(x',\zeta(x'))=0$ for every $x'\in C$. If there is $a'\in C$ such that  $\xi(a')=\zeta(a')$, then
$\xi(x')=\zeta(x')$ for every $x'\in C$.
\end{prop}

\begin{proof}
We use the same method of proof as in the preceding proposition. Let
$z_1=\xi(a')=\zeta(a'),\ldots,z_k$ be the distinct roots in ${\Bbb C}$ of the product
$P(a',X_n)Q(a',X_n)$. Let
$m_i$ (resp.\ $p_i$) be the multiplicity of $z_i$ as a root of $P(a',X_n)$ (resp.\ $Q(a',X_n)$),
where multiplicity zero means ``not a root''. The degree of $\mbox{\rm gcd}(P(a',X_n),Q(a',X_n))$ is
$\sum_{i=1}^k\min(m_i,p_i)$, and each $z_i$ has multiplicity $\min(m_i,p_i)$ as a root of this
gcd. Choose
$\varepsilon>0$ such that all disks
$D(z_i,\varepsilon)$ are disjoint. For every $x'\in C$ sufficiently close to $a'$, each disk
$D(z_i,\varepsilon)$ contains a root of multiplicity $m_i$ of $P(x',X_n)$ and a root of
multiplicity $p_i$ of $Q(x',X_n)$. Since the degree of $\mbox{\rm gcd}(P(x',X_n),Q(x',X_n))$ is
equal to $\sum_{i=1}^k\min(m_i,p_i)$, this gcd must have one root of multiplicity $\min(m_i,p_i)$ in
each disk
$D(z_i,\varepsilon)$ such that $\min(m_i,p_i)>0$. In particular, it follows that
$\xi(x')=\zeta(x')$. Since $C$ is connected, this equality holds for every $x'\in C$. \end{proof}

We have seen in Chapter 1 that the number of distinct complex roots of $P$ and the degree of the
gcd of $P$ and $Q$, can be computed from the fact that the principal subresultant
coefficients\index{subresultant!principal -- coefficient}
$\mbox{\rm PSRC}_i(P,P')$ and
$\mbox{\rm PSRC}_i(P,Q)$ are zero or nonzero, as long as the degrees (with respect to $X_n$) of
$P$ and $Q$ are fixed (cf. Corollary \ref{racinesdist} and Proposition \ref{CSP}).
For the values of the parameters (here, 
$X_1,\ldots,X_{n-1}$) such that some leading coefficients vanish, we have to use the
principal subresultant coefficients for the truncated polynomials. This leads us to the following
definition.
\par
If
$P$ is a polynomial in ${\Bbb R}[X_1,\ldots,X_n]$, we consider it as a polynomial in the
variable $X_n$ with coefficients in ${\Bbb R}[X_1,\ldots,X_{n-1}]$. We denote by $\lc(P)$ its
leading coefficient and by $\mbox{\rm trunc}(P)$ the truncated polynomial obtained by deleting its
leading term. Let $P_1,\ldots,P_r$ be a family of polynomials in ${\Bbb R}[X_1,\ldots,X_n]$. We
define
$\mbox{\rm PROJ}(P_1,\ldots,P_r)$ to be the smallest family of polynomials in ${\Bbb
R}[X_1,\ldots,X_{n-1}]$ satisfying the following rules:
\begin{itemize}
\item If $\deg_{X_n}P_i=d\geq 2$, $\mbox{\rm PROJ}(P_1,\ldots,P_i,\ldots,P_r)$ contains all
nonconstant polynomials among
$\mbox{\rm PSRC}_j (P_i,\partial P_i/\partial X_n)$ for $j=0,\ldots,d-1$.
\item If $1\leq d=\min(\deg_{X_n}(P_i),\deg_{X_n}(P_k))$, $\mbox{\rm
PROJ}(P_1,\ldots,P_i,\ldots,P_k,\ldots,P_r)$ contains all nonconstant
$\mbox{\rm PSRC}_j (P_i,P_k)$ for $j=0,\ldots,d$.
\item If $\deg_{X_n}P_i\geq 1$ and $\lc(P_i)$ is not constant,  $\mbox{\rm
PROJ}(P_1,\ldots,P_i,\ldots,P_r)$ contains $\lc(P_i)$ and $\mbox{\rm
PROJ}(P_1,\ldots,\mbox{\rm trunc}(P_i),\ldots,P_r)$.
\item If $\deg_{X_n}P_i=0$ and $P_i$ is not constant,  $\mbox{\rm
PROJ}(P_1,\ldots,P_i,\ldots,P_r)$ contains $P_i$.
\end{itemize}

The following theorem is a consequence of the results previously proved in this section.

\begin{thm}\label{cylindre}
Let $(P_1,\ldots,P_r)$ be a family of polynomials in ${\Bbb R}[X_1,\ldots,X_n]$, and let $C$ be a
connected, $\mbox{\rm
PROJ}(P_1,\ldots,P_r)$-invariant\index{invariant@{$(P_1,\ldots,P_r)$-invariant}}, semialgebraic
subset of ${\Bbb R}^{n-1}$. Then there are continuous semialgebraic functions
$\xi_1<\ldots<\xi_\ell : C\to {\Bbb R}$, such that, for every $x'\in C$, the set
$\{\xi_1(x'),\ldots,\xi_\ell(x')\}$ is the set of real roots of all nonzero polynomials
$P_1(x',X_n),\ldots,P_r(x',X_n)$. The graph of each
$\xi_i$, and each band of the cylinder $C\times {\Bbb R}$ bounded by these graphs, are
connected semialgebraic sets, semialgebraically homeomorphic to
$C$ or $C\times(0,1)$, respectively, and $(P_1,\ldots,P_r)$-invariant. 
\end{thm}

If we have constructed a c.a.d.\ of ${\Bbb R}^{n-1}$
adapted\index{cylindrical algebraic decomposition!adapted} to $\mbox{\rm PROJ}(P_1,\ldots,P_r)$, the
preceding theorem can be used to extend this c.a.d.\ to a c.a.d.\ of ${\Bbb R}^n$ adapted to
$(P_1,\ldots,P_r)$. On the other hand, by iterating
$(n-1)$ times the operation PROJ, we arrive to a finite family of polynomials in one variable $X_1$.
It is easy to construct a c.a.d.\ of ${\Bbb R}$ adapted to this family: the real roots of the
polynomials in the family cut the line in finitely many points and open intervals.
Finally, we obtain:

\begin{thm}
For every finite family  $P_1,\ldots,P_r$ in ${\Bbb R}[X_1,\ldots,X_n]$,
there is an adapted\index{cylindrical algebraic decomposition!adapted} c.a.d. of ${\Bbb R}^n$.
\end{thm}

We illustrate this result by constructing a 
c.a.d.\index{cylindrical algebraic decomposition!adapted}\ of
${\Bbb R}^3$ adapted to the polynomial
$P=X^2+Y^2+Z^2-1$. The Sylvester matrix\index{Sylvester matrix} of $P$ and $\partial P/\partial Z$ is
$$\left(\begin{array}{ccc}
1&0&X^2+Y^2-1\\
0&2&0\\
2&0&0\end{array}\right)\;.$$
Hence, ${\rm PSRC}_0(P,\partial P/\partial Z)=-4(X^2+Y^2-1)$ and ${\rm PSRC}_1(P,\partial P/\partial
Z)=2$. Getting rid of irrelevant constant factors, we obtain $\mbox{\rm PROJ}(P)=(X^2+Y^2-1)$, then 
$\mbox{\rm PROJ}(\mbox{\rm PROJ}(P))=(X^2-1)$. The c.a.d. obtained is represented on Figure
\ref{sphere}.

\begin{figure}[ht]
\hspace*{20mm}\dessin 120mm by 90mm (sphereeps)
\caption{A c.a.d. adapted\index{cylindrical algebraic decomposition!adapted} to the
sphere\label{sphere}}
\end{figure}

\debexo How many cells\index{cell} of ${\Bbb R}^3$ are there in this c.a.d.? Is it possible to have a
c.a.d. of ${\Bbb R}^3$, such that the sphere is the union of cells,  with less cells?
\finexo

\debexo Let $f:A\to {\Bbb R}$ be a semialgebraic function, which is not supposed to be continuous.
Show that there exists a finite semialgebraic partition $A=\bigcup_{i=1}^s C_i$ of $A$ such that, for
every $i$, the restriction of $f$ to $C_i$ is continuous. Hint: use a c.a.d.\ adapted to the graph
of $f$.
\finexo

\noindent{\bf Remark.} We return to the proof of Proposition \ref{fracines}. For every $x'\in C$,
$\xi_i(x')$ is a root of
$P(x',X_n)$, with constant multiplicity $m_i$. Hence, $\xi_i(x')$ is a simple root of the
$(m_i-1)^{\rm th}$ derivative with respect to $X_n$. Therefore, if $C$ is a ${\cal C}^\infty$
submanifold of
${\Bbb R}^{n-1}$,  the function $\xi_i$ is ${\cal C}^\infty$ on
$C$. The graphs and the bands\index{band} of the cylinder $C\times{\Bbb R}$ are also ${\cal
C}^\infty$ submanifolds, diffeomorphic to $C$ or
$C \times (0,1)$, respectively (cf. the formulas in the proof of Proposition
\ref{paves}). By induction on $n$, one proves in this way that every semialgebraic set is the
disjoint union of finitely many semialgebraic subsets $C_i$, which are ${\cal C}^\infty$
submanifolds each semialgebraically diffeomorphic to an open hypercube $(0,1)^{d_i}$. The
semialgebraic ${\cal C}^\infty$ submanifolds are called {\em Nash manifolds\index{Nash manifold}}.

\subsection{The c.a.d.\ algorithm}

We shall now make more precise the algorithmic aspect of the construction of the
c.a.d.\index{cylindrical algebraic decomposition}. The c.a.d.\ algorithm receives as input a finite
list of polynomials
$P_1,\ldots  ,P_r$ in ${\Bbb Q}[X_1,\ldots ,X_n]$ and produces as output the list of the
cells\index{cell} of a c.a.d.\ adapted\index{cylindrical algebraic decomposition!adapted} to
$P_1,\ldots ,P_r$ (with information on their cylindrical arrangement), together with a ``test
point'' in each cell, whose coordinates are rational or real algebraic numbers. The algorithm works
in the following way:
\begin{itemize}
\item Given a list of polynomials in one variable $X_1$, it counts and isolates in intervals
with rational endpoints all real roots of these polynomials (this can be done by using
Sturm's method\index{Sturm sequence}). The cells of
${\cal C}_1$  are the roots and the intervals between the roots. The roots are characterized by
the polynomial equation they satisfy and the interval with rational endpoints which isolates
them. These endpoints may be taken as test points for the intervals between the roots (but there
may be more convenient choices).
\item Given a list $(P_1,\ldots ,P_r) $ of polynomials in ${\Bbb Q}[X_1,\ldots ,X_n]$, 
where $n>1$, it computes ${\rm PROJ}(P_1,\ldots ,P_r)$ and calls the c.a.d.\ algorithm for this
list of polynomials in
${\Bbb Q}[X_1,\ldots ,X_{n-1}]$. One obtains 
${\cal C}_{n-1}$, which is a partition of ${\Bbb R}^{n-1}$ in ${\rm PROJ}(P_1,\ldots
,P_r)$-invariant cells, and a test point $a'_C$ for each cell  $C\in {\cal C}_{n-1}$. For such a
cell $C$, one can apply Theorem \ref{cylindre} to cut the cylinder $C\times {\Bbb R}$ in 
$(P_1,\ldots ,P_r)$-invariant\index{invariant@{$(P_1,\ldots,P_r)$-invariant}} cells. In order to know
how many cells there are in this cylinder and to produce a test point for each cell, one computes the
real roots of
$P_1(a'_C,X_n),\ldots ,P_r(a'_C,X_n)$. Sturm's method can be used once again, but the coefficients
of the polynomials may be real algebraic numbers.
\end{itemize}
We encounter here the problem of coding  real algebraic numbers and computing with them. One
possible method is to give the coordinates of the test point   $a'_C$ as polynomials in the
primitive element of the extension of  ${\Bbb Q}$ that they generate; this primitive element
is given by its minimal polynomial over ${\Bbb Q}$ and an isolating interval with rational
endpoints. 
\par
The c.a.d.\index{cylindrical algebraic decomposition}\ algorithm can be used to solve the following
problem: decide wether a formula\index{formula!first order} without free variables is true or false
({\em decision problem\index{decision problem}}). More generally, given a first order formula, the
c.a.d.\ algorithm allows one to decide whether the semialgebraic set $S$ defined by this formula is
empty or not, and in case not, produces a point in (every connected component\index{connected
component} of)
$S$.\par
The complexity of the algorithm is {\em doubly exponential} in the number of variables (free
variables or quantified ones) of the formula. This double exponential feature can be explained by
the following remark. The resultant\index{resultant} (with respect to $X_n$) of two polynomials of
total degree
$d$ is, in general, of total degree $d^2$. Hence, taking the $\mbox{\rm PROJ}$ of a family of
polynomials in $n$ variables of maximum degree $d$ will give a family of polynomials in $n-1$
variables of maximum degree $d^2$. Iterating the $\mbox{\rm PROJ}$ operation $n-1$ times will give
a family of polynomials in one variable of maximum degree $d^{2^{n-1}}$.\par
The double exponential explains why the practical applications of the c.a.d.\ algorithm are very
limited. In order to reduce the complexity, it is necessary to avoid the method of successive
codimension $1$ projections. We shall see in Chapter 4 an alternative method: the {\em
critical\index{critical!point} point method}.

\subsection{Connected components of semialgebraic sets}

The c.a.d.\index{cylindrical algebraic decomposition}\ shows that every semialgebraic set $S \subset
{\Bbb R}^n$ is the disjoint union of semialgebraic subsets $C_1 ,\ldots,C_p$, such that each
$C_i$ is semialgebraically homeomorphic to an open hypercube $(0,1)^{d_i}$ (with
$(0,1)^0$ = a point). Each $C_i$ is obviously connected.

\begin{thm}\label{conncomp}
Every semialgebraic set has finitely many connected components\index{connected component} which are
semialgebraic. Every semialgebraic set is locally connected.
\end{thm}

\begin{proof}
Using the notation above, we shall say that $C_i$ is adjacent
to $C_j$ if $C_i  \cap  \clos(C_j) \not= \emptyset$. Let $\sim$ be the equivalence relation
generated by the relation ``is adjacent to'':
$C_i \sim C_j$ if there is a sequence $C_i = C_{i_0} , C_{i_1} ,
\ldots, C_{i_q} = C_j$ such that $C_{i_j}$ is adjacent to $C_{i_{j+1}}$ or
$C_{i_{j+1}}$ is adjacent to $C_{i_j}$. We obtain a partition of $S$ into finitely many
semialgebraic subsets
$S_1 ,\ldots, S_r$ , where each
$S_i$ is the union of all $C_j$ in the same equivalence class for $\sim$. Each
$S_i$ is closed in $S$. Indeed, if $C_i  \cap  \clos(S_j) \not= \emptyset$,
then $C_j$ is adjacent to some $C_k \subset S_i$ and, hence, $C_j \subset S_i$.
Since there are finitely many $S_i$, each one is also open in $S$. We now show that each  $S_i$
is connected. If
$S_i = F_1 \cup F_2$, where $F_1$ and
$F_2$ are disjoint closed subsets of $S_i$, we have:
\begin{itemize}
\item every $C_j \subset S_i$ is contained in $F_1$ or in $F_2$, since $C_j$ is connected;
\item if $C_j$ and $C_k$ are contained in $S_i$ and $C_j$ is adjacent to $C_k$,
then $C_j$ and $C_k$ are both contained in
$F_1$ or both in $F_2$.
\end{itemize}
According to the definition of $S_i$, it follows that $S_i = F_1$ or $S_i = F_2$. The first part of
the theorem is proved. \par
The semialgebraic set $S\subset {\Bbb R}^n$ is locally connected if, for every
$x \in S$, every open ball $B$ with center $x$ contains a connected neighborhood of $x$ in $S$. Since
$B
\cap S$ is semialgebraic, it has a finite number of connected components. The connected component
of $B\cap S$ containing 
$x$ is a connected neighborhood of $x$ in $S$.
\end{proof}

\chapter{Triangulation of semialgebraic sets}

In the preceding chapter, we have decomposed semialgebraic sets into simple pieces (the
cells\index{cell}, which are semialgebraically homeomorphic to open hypercubes). We have also
explained an algorithm which produces this decomposition. But the result obtained is not quite
satisfactory, for the following reasons:
\begin{itemize}
\item We do not have a description of cells of the c.a.d.\index{cylindrical algebraic
decomposition}\ by a boolean combination of polynomial equations and inequalities. In particular,
the c.a.d.\ algorithm described in Chapter 2 does not suffice to eliminate
quantifiers\index{quantifier elimination}.
\item We have no information concerning which cells of a c.a.d.\ are adjacent to others,
except for the cells in a cylinder. We do not know, in general, what happens when we pass from a
cylinder to another. In the case of the c.a.d.\ 
adapted\index{cylindrical algebraic decomposition!adapted} to the sphere, it is not
difficult to determine the topology from the cell decomposition. The two functions on the disk
$x^2+y^2<1$, whose graphs are the two open hemispheres, have an obvious extension by continuity on
the closed disk.  We show an example where this is not so.\par Take
$P=XYZ-X^2-Y^2$. We have 
\begin{eqnarray*}
\mbox{\rm PROJ}(P)&=&(XY,-X^2-Y^2)\;,\\
\mbox{\rm PROJ}(\mbox{\rm PROJ}(P))&=&(X^4,4X^2,X)\;.
\end{eqnarray*}
 The c.a.d.\ of
${\Bbb R}^2$ consists of 9 cells determined by the signs of $X$ and $Y$. The cylinders over each
open quadrant have three cells; the sign of $P$ in these cells is 
$$\pmatrix{+\cr0\cr-\cr}\quad \mbox{\rm where } XY>0,\quad
\pmatrix{-\cr0\cr+\cr}\quad \mbox{\rm where } XY<0\;.$$
The cylinders over each open half-axis have one cell, on which $P$ is negative. The cylinder over the
origin has one cell, on which $P=0$. This information is not sufficient to determine the topology of
the surface $P=0$.\par
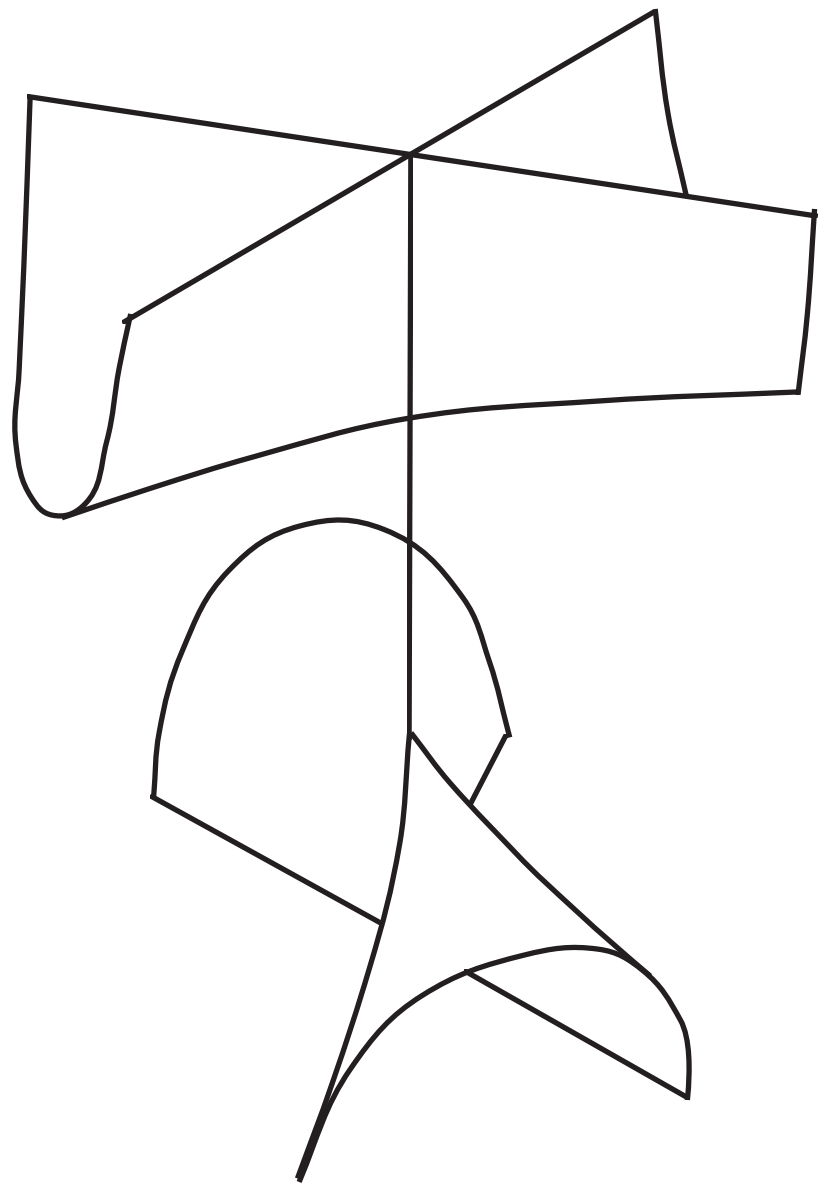
\begin{figure}
\setlength{\unitlength}{1mm}
\begin{picture}(100,80)
\put(35,0){\dessin 65mm by 80mm (doubleumb)}
\end{picture}
\caption{The surface $XYZ-X^2-Y^2=0$.}
\end{figure}

The main difference between the example of the sphere and the example above is the fact the
polynomial $P$ in the latter example is not monic as a polynomial in $Z$: its leading coefficent
$XY$ vanishes, and the polynomial is even identically zero for $X=Y=0$. The functions describing
the zeros of $P$ on each open quadrant have no extension by continuity on the closed quadrants.
\end{itemize}

\debexo
Return to the last example $P=XYZ-X^2-Y^2$, taking the variables in the order $(Z,X,Y)$. Compute
a c.a.d.\  adapted to $P$. Is it possible to recover the topology of $P=0$ from this c.a.d.?
\finexo

We shall see in this chapter how to modify the c.a.d.\index{cylindrical algebraic decomposition}
algorithm in order to solve the two problems listed above.

\section{Thom's lemma}

\subsection{For one variable}

We introduce a notation, which we shall use for relaxing inequalities. If $\varepsilon\in\{-1,0,1\}$
is a sign, we denote 
$$\overline{\varepsilon}=\cases{
\{0\}& if $\varepsilon=0$,\cr
\{0,1\}& if $\varepsilon=1$,\cr
\{0,-1\}& if $\varepsilon=-1$.\cr}$$

\begin{prop}[Thom's lemma]\index{Thom's lemma}
 Let $P_1 ,\ldots, P_s \in {\Bbb R}[X]$ be a finite family of nonzero polynomials, which is
closed under derivation (i.e., if the derivative $P'_i$ is nonzero, there is $j$ such that
$P'_i=P_j$). For $\varepsilon=(\varepsilon_1,\ldots,\varepsilon_s)\in\{-1,0,1\}^s$, let
$A_\varepsilon \subset {\Bbb R}$ be defined by 
$$A_\varepsilon= \{ x \in {\Bbb R}\ ;\ \sign(P_i(x))=\varepsilon_i
\ \mbox{ for }\ i=1,\ldots,s\}\;.$$
 Then
\begin{itemize}\setlength{\itemsep}{-2pt plus1pt minus1pt}
\item either $A_\varepsilon = \emptyset$,
\item or $A_\varepsilon$ is a point (necessarily, at least one of the $\varepsilon_i$ is $0$),
\item or $A_\varepsilon$ is a nonempty open interval (necessarily, all $\varepsilon_i$ are $\pm 1$).
\end{itemize}
Let 
$$A_{\overline \varepsilon} = \{ x \in {\Bbb R}\ ;\ \sign(P_1(x))\in\overline\varepsilon_1 \et
\ldots\et  \sign(P_s(x))\in\overline\varepsilon_s \}\;,$$
which is obtained by relaxing the strict inequalities. Then
$A_{\overline \varepsilon}$ is either empty, or a point, or a closed interval different from a point
(and the interior of this interval is $A_\varepsilon$).
\end{prop}

\begin{proof}
The proof is by induction on the number $s$ of polynomials in the list. If $s = 1$, the only
polynomial must be a nonzero constant, and in this case the statement holds trivially
($A_\varepsilon =
\emptyset$ or ${\Bbb R}$). We prove the induction step from $s$ to $s+1$. We can assume that
$P_{s+1}$ has maximal degree in the list $(P_1 ,\ldots, P_{s+1})$. Then the list
$(P_1 ,\ldots, P_s)$ is also closed under derivation. By the inductive assumption, if
$\varepsilon\in\{-1,0,1\}^s$, the subset
$A_\varepsilon\subset {\Bbb R}$ is either empty, or a nonempty open interval, or a point.
Moreover, if $A_\varepsilon$ is a nonempty open interval,
$P_{s+1}$ is monotone on $A_\varepsilon$, since
$P'_{s+1}$ has constant sign on $A_\varepsilon$. It follows that, for every $\varepsilon_{s+1}\in
\{-1,0,1\}$, the set 
$$A_{\varepsilon,\varepsilon_{s+1}}=A_\varepsilon \cap \{ x\in {\Bbb R}\ ;\
\sign(P_{s+1}(x))=\varepsilon_{s+1}\}$$ 
 satisfies the properties of the proposition.
\end{proof}

\debexo
Let $a$ and $b$ be distinct real roots of $P\in{\Bbb R}[X]$. Show that there exists a derivative
$P^{(i)}$ such that $P^{(i)}(a)P^{(i)}(b)<0$.
\finexo

Thom's lemma\index{Thom's lemma} allows one to answer to the first problem concerning
c.a.d.\index{cylindrical algebraic decomposition}: in order to obtain a description of each
cell\index{cell} by a boolean combination of polynomial equations or inequalities, it is enough to
add the derivatives. More precisely:

\begin{cor}
 Let $(P_{i, j})$,
 where $1 \leq i \leq n$ and $1 \leq j \leq s_i$, be a family of nonzero polynomials such that:
\begin{itemize}\setlength{\itemsep}{-2pt plus1pt minus1pt}
\item For fixed $i$,  $P_{i,1} ,\ldots, P_{i,s_i}$ is a family of polynomials in ${\Bbb R}[X_1
,\ldots, X_i]$ which is closed under derivation with respect to $X_i$.
\item For $i < n$, the family of polynomials $(P_{i,1} ,\ldots, P_{i,s_i})$ contains the
family\linebreak
$\mbox{\rm PROJ}(P_{i+1,1} ,\ldots, P_{i+1,s_{i+1}})$.
\end{itemize}
Then the families ${\cal C}_k$, for $k=1,\ldots,n$, consisting of all nonempty semialgebraic
susbsets of ${\Bbb R}^k$ of the form
 $$\{ x \in
{\Bbb R}^k\ ;\     \sign(P_{i, j}(x))=\varepsilon_{i,j}\ \mbox{\rm for } i = 1,\ldots,k  \et  j =
1,\ldots,s_i\}\;,$$
where $\varepsilon_{i, j}\in\{-1,0,1\}$, constitute a c.a.d.\index{cylindrical algebraic
decomposition}\ of
${\Bbb R}^n$.
\end{cor}

\subsection{For several variables}

In the situation of Thom's lemma\index{Thom's lemma}, the different ``pieces'' (points and open
intervals) are described by sign conditions on polynomials, and their closures are obtained by
relaxing the strict inequalities. We shall extend these nice properties to the case of several
variables. We return to the c.a.d.\index{cylindrical algebraic decomposition}, in order to see that,
in ``good situations'', we can control what happens when we pass from a cylinder $C\times{\Bbb R}$
to another
$C'\times{\Bbb R}$ such that
$C'\subset\clos(C)$.
\par
We shall say that a nonzero polynomial $P
\in {\Bbb R}[X_1,\ldots, X_n]$ is {\em quasi-monic\index{quasi-monic}} with respect to $X_n$ if its
leading coefficient is a constant (we consider $P$ as a polynomial in $X_n$ with coefficients in
${\Bbb R}[X_1,\ldots,X_{n-1}]$).\par
Consider the following situation: \begin{itemize}
\item $P_1 ,\ldots, P_s$ is a list of polynomials in ${\Bbb R}[X_1,\ldots, X_n]$, all quasi-monic
with respect to $X_n$, closed under derivation with respect to  $X_n$.
\item $C$ and $C'$ are both connected, $\mbox{\rm PROJ}(P_1 ,\ldots,
P_s)$-invariant\index{invariant@{$(P_1,\ldots,P_r)$-invariant}}, semialgebraic subsets of ${\Bbb
R}^{n-1}$,
 and $C'$ is contained in the closure of $C$.
 \end{itemize}
It follows from Theorem \ref{cylindre} that there are continuous semialgebraic functions $\xi_1
<\ldots< \xi_\ell : C \to {\Bbb R}$  and  $\xi'_1 <\ldots< \xi'_{\ell'} : C' \to {\Bbb R}$,
which describe, as functions of $x=(x_1,\ldots,x_{n-1})$, the real roots of the polynomials
$P_1(x,X_n),\ldots,P_s(x,X_n)$. Denote by
$A_j$ and $A'_j$ the graphs of $\xi_j$ and $\xi'_j$, respectively. Denote by $B_j$ and $B'_j$ the
bands of the cylinders
$C\times {\Bbb R}$ and $C' \times {\Bbb R}$, respectively, which are cut by these graphs.

\begin{lem}\label{prolongement}
In the situation above:
\begin{enumerate}\setlength{\itemsep}{0pt plus1pt minus1pt}
\item Every function $\xi_j$ can be continuously extended to $C'$, and this extension coincides
with one of the functions $\xi'_{j'}$.
\item For every function $\xi'_{j'}$ , there is a function $\xi_j$ whose extension by continuity to
$C'$ is $\xi'_{j'}$ .
\item For every $\varepsilon=(\varepsilon_1,\ldots,\varepsilon_s)\in\{-1,0,1\}^s$, the set
$$E_\varepsilon = \{ (x, x_n) \in C \times {\Bbb R}\ ;\  \sign(P_i(x,x_n))=\varepsilon_i \ \mbox{
for }\ i=1,\ldots,s
\}$$ is either empty, or one of the $A_j$, or one of the $B_j$. Let $E_{\overline\varepsilon}$
be the subset of $C\times{\Bbb R}$ obtained by relaxing the strict inequalities:  
$$E_{\overline\varepsilon} = \{ (x, x_n) \in C \times {\Bbb R}\
;\  \sign(P_i(x, x_n))\in\overline\varepsilon_i  \ \mbox{ for }\ i=1,\ldots,s\}\;,$$ and let
$$E'_{\overline\varepsilon} = \{ (x, x_n) \in C' \times {\Bbb R}\
;\  \sign(P_i(x, x_n))\in\overline\varepsilon_i  \ \mbox{ for }\ i=1,\ldots,s\}\;.$$
If $E_\varepsilon \not= \emptyset$, we have  $\clos(E_\varepsilon) \cap (C \times {\Bbb R}) = 
E_{\overline\varepsilon}$ and 
$\clos(E_\varepsilon) \cap (C' \times {\Bbb R}) = E'_{\overline\varepsilon}$. Moreover,
$E'_{\overline\varepsilon}$ 
 is either a graph $A'_{j'}$ , or the closure of one of the bands $B'_{j'}$ in
$ C' \times {\Bbb R}$.
\end{enumerate}
\end{lem}

\begin{proof}
Let $x' \in C'$. Choose a function $\xi_j$. There is a polynomial $P_\mu$ of the family
such that, for every $x\in C$, $\xi_j(x)$ is a simple root of
$$P_\mu (x, X_n) = a_0 X_n^d + a_1(x) X_n^{d-1} +\cdots+ a_d(x)\;,$$
where $a_0$ is a nonzero constant. Set
$$M(x') = \max_{i=1,\ldots,d}\left(d\,\left\vert {a_i(x')\over a_0}\right\vert
\right)^{1/i}\;.$$
By Proposition \ref{borne},  there is a neighborhood  $U$  of  $x'$ in ${\Bbb R}^{n-1}$ such that,
for every $x\in U\cap C$, we have $\xi_j(x) \in [-M(x')-1, M(x')+1]$. Choose a sequence
$(x^\nu)$ in
$C$, such that $\lim_{\nu\to\infty} x^\nu=x'$. The sequence $\xi_j(x^\nu)$ is bounded and has,
therefore,  a $\limsup$ $y' \in [-M(x')-1, M(x')+1]$. The point $(x',y')$ belongs to the closure of
the graph of $\xi_j$. Let
$\varphi_1=\sign(P_\mu'(x,\xi_j(x))$,\ldots, $\varphi_d=\sign(P_\mu^{(d)}(x,\xi_j(x))$,
for $x\in C$ (observe that these signs are constant for $x\in C$). Every point  $(x', x'_n)$
in the closure of the graph of $\xi_j$ must satisfy
$$P_\mu(x', x'_n) = 0,\  \sign(P'_\mu(x', x'_n))\in\overline\varphi_1,\ldots,\ 
\sign(P_\mu^{(d)})(x', x'_n)\in\overline\varphi_d\;.$$
By Thom's lemma, there is at most one $x'_n$
satisfying these inequalities. It follows that $\xi_j$ extends continuously at
$x'$. Hence, it extends continuously to $C'$, and this extension coincides with one of the functions
$\xi'_{j'}$. This proves 1.\par
We now prove 2. Choose a function $\xi'_{j'}$.  Since $\xi'_{j'}$ is a simple root of some
polynomial $P_\nu$ in the family, it follows from the implicit function theorem that there is a
function $\xi_j$, also a root of $P_\nu$, whose extension by continuity to $C'$ is
$\xi'_{j'}$ .\par
We now turn to 3. The properties of $E_\varepsilon$ and $E_{\overline\varepsilon}$ are
straightforward consequences of Thom's lemma, since 
$P_1 ,\ldots, P_s$ have constant signs on each graph
$A_j$ and each band $B_j$, and the closure of $B_j$ in $C \times {\Bbb
R}$ is $A_j  \cup B_j \cup A_{j+1}$ (as usual, $A_0 = \emptyset = A_{\ell+1}$). It is obvious that
$\clos(E_\varepsilon)\cap (C'\times {\Bbb R}) \subset E'_{\overline\varepsilon}$. It follows from 1
and 2 that
$\clos(E_\varepsilon) \cap (C' \times {\Bbb R})$ is either a graph $A'_j$ or the closure of one of
the bands
$B'_{j'}$ in $C'\times {\Bbb R}$. By Thom's lemma, this is also the case for
$E'_{\overline\varepsilon}$. It remains to check  that the equality holds if
$E'_{\overline\varepsilon}$ is the closure of a band $B'_{j'}$. In this
case, all
$\varepsilon_i$ must be $\pm 1$, and the sign of $P_i$ is $\varepsilon_i$ on every
sufficiently small neighborhood $V$ of a point $x'$ of $B'_{j'}$. This implies  $V\cap (C\times{\Bbb
R})\subset E_\varepsilon$ and, hence, $x'\in\clos(E_\varepsilon)$. This shows that
$\clos(E_\varepsilon)\cap (C' \times {\Bbb R})$ is also the closure of $B'_{j'}$.
\end{proof}

The following theorem gives an answer to the problem of determining which cells\index{cell} of a
c.a.d.\index{cylindrical algebraic decomposition}\ are adjacent to another. 

\begin{thm}\label{LTG}
 Let $(P_{i, j})$ be a family of polynomials with real coefficients,
$1
\leq i \leq n$, $1 \leq j \leq s_i$, such that:
\begin{itemize}\shitem
\item
for fixed $i$, $(P_{i,1} ,\ldots, P_{i,s_i})$ is a family of polynomials in ${\Bbb R}[X_1 ,\ldots,
X_i]$, all quasi-monic with respect to  $X_i$, closed under derivation with respect to $X_i$,
\item
 for $i < n$, the family of polynomials $(P_{i,1} ,\ldots, P_{i,s_i})$ contains the family\linebreak
$\mbox{\rm PROJ}(P_{i+1,1} ,\ldots, P_{i+1,s_{i+1}})$.
\end{itemize}
For $0<k\leq n$, given a family $\varepsilon=(\varepsilon_{i,j})$ of signs in $\{-1,0,1\}$
indexed by $i=1,\ldots,k$ and $j=1,\ldots,s_i$, set
\begin{eqnarray*}
C_\varepsilon &= &\{ x \in {\Bbb
R}^k\ ;\     \sign(P_{i, j}(x))=\varepsilon_{i, j}\ \mbox{\rm for } i = 1,\ldots,k  \et  j =
1,\ldots,s_i\}\;,\\
C_{\overline\varepsilon}& =& \{ x \in {\Bbb R}^n\ ;\    \sign( P_{i, j}(x))\in\overline{\varepsilon_{i,
j}}\ \hbox{\rm for } i = 1,\ldots,n  \et  j = 1,\ldots,s_i\}\;.
\end{eqnarray*}
Then the non empty $C_\varepsilon$ are the cells of a c.a.d.\ of ${\Bbb R}^n$, and the closure of a
nonempty cell $C_\varepsilon$ is
$C_{\overline\varepsilon}$, which is a union of cells.
\end{thm}

	The proof of the theorem is by induction on $n$, using the preceding lemma for the induction step
and Thom's lemma for $n = 1$. This theorem may be seen as a  {\em generalized Thom
lemma\index{Thom's lemma!generalized}}. Observe that the cells $C_\varepsilon$ are actually Nash
submanifolds\index{Nash manifold}, semialgebraically diffeomorphic to open hypercubes. The theorem
above holds for a family of polynomials with special properties.  Nevertheless, any finite family of
polynomials in ${\Bbb R}[X_1,\ldots,X_n]$ can be, up to a linear change of variables, completed to a
family satisfying these properties.

\begin{prop}\label{chgvar} Let $P_1 ,\ldots, P_\ell \in {\Bbb R}[X_1,\ldots, X_n]$. There is a linear
automorphism $u : {\Bbb R}^n
\to {\Bbb R}^n$ and a family of polynomials $(P_{i, j})$ satisfying the conditions of Theorem
\ref{LTG}, such that
$P_{n, j}(X) = P_j(u(X))$ for $j = 1,\ldots,\ell$ (where $X = (X_1 ,\ldots, X_n)$). \end{prop}

\begin{proof}
First, there is a linear change of variables	
$$v(X_1 ,\ldots, X_n) = (X_1 + a_1 X_n ,\ X_2 + a_2 X_n ,\ldots,\ X_{n-1} + a_{n-1} X_n
,\ X_n)$$
such that all polynomials $P_1(v(X)) ,\ldots, P_\ell(v(X))$ are quasi-monic with respect to $X_n$.
Indeed, if
$P_i(X) =\Pi_i(X) +\cdots$, where $\Pi_i$ is the homogeneous part of highest degree (say
$d_i$) of $P_i$, then $P_i(v(X)) = X_n^{d_i} \Pi_i(a_1 ,\ldots, a_{n-1} , 1) +$ terms of lower
degree with respect to $X_n$. It suffices to choose $a_1 ,\ldots, a_{n-1}$ such that none of the
$\Pi_i(a_1 ,\ldots, a_{n-1} , 1)$ is zero. 
Then we add to the list of polynomials $P_1(v(X))$,\dots, $P_\ell(v(X))$ all their nonzero
derivatives of every order with respect to $X_n$, say
$P_{\ell+1}$,\dots, $P_{s_1}$. Now compute $(Q_1 ,\ldots,
Q_t)=\mbox{\rm PROJ}(P_1(v(X)),\ldots,P_\ell(v(X)), P_{\ell+1} ,\ldots, P_{s_1})$. Using induction,
there is a linear automorphism
$u' : {\Bbb R}^{n-1}
\to {\Bbb R}^{n-1}$ and a family $(P_{i,j} )_{1\leq i\leq n-1\;,\ 1\leq  j\leq s_i}$ of polynomials
satisfying the conditions of the theorem and
such that $P_{n-1, j}(X') = Q_j(u'(X'))$, for
$j = 1,\ldots,t$, where $X' = (X_1 ,\ldots, X_{n-1})$. Finally, set $u  =  (u' \times
\mbox{\rm Id})\circ v$ (where $(u' \times \mbox{\rm Id })(X', X_n) = (u'(X'), X_n)$), $P_{n, j}(X) =
P_j(u(X))$ for $1 \leq j \leq \ell$ and $P_{n,  j}(X) = P_j(u'(X'), X_n)$ for $\ell+1 \leq j \leq
s_1$.
\end{proof}

\begin{cor}\label{stratification}
Let $S \subset {\Bbb R}^n$ be a semialgebraic set and $T_1,\ldots,T_q$, finitely many semialgebraic
subsets of $S$. Then
$S$  can be decomposed as a disjoint finite union $S = \bigcup_{i=1}^p C_i$, where 
\begin{itemize}\shitem
	\item every $C_i$ is semialgebraically homeomorphic (and even diffeomorphic) to an open
hypercube $(0,1)^{d_i}$,
	\item the closure of $C_i$ in $S$ is the union of $C_i$ and some
$C_j$'s, $j \not= i$, with $d_j < d_i$,
\item every $T_k$ is the union of some $S_i$.
\end{itemize}
\end{cor}

\begin{proof}
We start with a list of polynomials $(P_1,\ldots,P_\ell)$ such that $S$ and all $T_k$ are
described by boolean combinations of sign conditions on polynomials of this list. We use
Proposition \ref{chgvar} to be in the conditions of Theorem \ref{LTG}. Then $S$  and all $T_k$ are
the unions of cells
$C_\varepsilon$ of this theorem. If $C_\varepsilon \not= \emptyset$, then
$C_{\overline\varepsilon}$ is the union of $A_\varepsilon$ and some $A_{\varepsilon'}$,
$\varepsilon' \not=
\varepsilon$. We can check by induction on $n$ that $d_{\varepsilon'} <
d_\varepsilon$.
\end{proof}

A decomposition $S = \bigcup_i C_i$ as in the above corollary is called a
{\em stratification\index{stratification}} of $S$, and the $C_i$ are called {\em
strata\index{stratum}} of this stratification.

\debexo
We use the notation of Theorem \ref{LTG}. Assume that $C_\varepsilon$ is non\-empty and
bounded. Show that the semialgebraic diffeomorphism $(0,1)^d\to C_\varepsilon$ induced by the
c.a.d.\ extends to a surjective continuous mapping $[0,1]^d\to \clos(C_\varepsilon)$.
\finexo

\subsection{The finiteness theorem}

The following theorem can be obtained as a consequence of the generalized Thom lemma\index{Thom's
lemma!generalized}.	

\begin{thm}[Finiteness Theorem]\label{finiteness}
 Let $S \subset {\Bbb R}^n$ be a semialgebraic set and $U$, a semialgebraic subset of $S$. Then
$U$ can be written as a finite union of open semialgebraic subsets of $S$ of the form
$$ \{ x \in S\ ;\ P_1(x) > 0  \et\ldots \et  P_s(x) > 0 \}\;,$$
where $P_1,\ldots, P_s \in {\Bbb R}[X_1,\ldots, X_n]$.
\end{thm}

	\begin{proof} Up to a linear change of variables, we can assume that $S$ and $U$ are defined by
boolean combinations of sign conditions on polynomials of a list $(P_{i, j})$ satisfying
the conditions of the generalized Thom lemma \ref{LTG}. It follows that $S$ and $U$ are
finite unions of cells  $C_\varepsilon$ (with the notation of Theorem \ref{LTG}). Define
$C_{\underline\varepsilon}$ to be the set of $x \in {\Bbb R}^n$ such that
$\sign(P_{i, j}(x))=\varepsilon_{i, j}$ for every $i, j$ such that $\varepsilon_{i, j}=\pm1$ (i.e. 
we drop the equations and keep only the strict inequalities). Obviously, $C_{\underline\varepsilon}$
is an open semialgebraic set containing $C_\varepsilon$, and
$C_{\underline\varepsilon}$ is the union of $C_\varepsilon$ and the $C_{\varepsilon'}$ for all
$\varepsilon'$ such that $\varepsilon_{i,j}=\pm1 \Rightarrow \varepsilon_{i,j}=\varepsilon'_{i,j}$.
The open subset $S\cap A_{\underline\varepsilon}$ of $S$ is defined by a conjunction of strict
polynomial inequalities. Hence, it suffices to prove that
$$U =  \bigcup_{C_\varepsilon\subset U} (C_{\underline\varepsilon} \cap S)\;.\leqno{(*)}$$ 
The inclusion of the left-hand side into the right-hand side is clear, since $C_\varepsilon \subset
U$ implies $C_\varepsilon\subset C_{\underline\varepsilon} \cap S$. Let $C_\varepsilon \subset
U$ and $C_{\varepsilon'} \subset C_{\underline\varepsilon} \cap S$. Then
$$C_{\overline{\varepsilon'}} = \clos(C_{\varepsilon'}) \supset C_\varepsilon\;,$$ 
since $\varepsilon'_{i,j}=\pm1$ implies $\varepsilon_{i,j}\in\overline{\varepsilon'_{i,j}}$.
Therefore
$\clos(C_{\varepsilon'}) \cap U \not= \emptyset$ and, since $U$ is open in $S$,
$C_{\varepsilon'} \cap U \not= \emptyset$. It follows that $C_{\varepsilon'} \subset U$. Hence, the
equality $(*)$ is proved.
\end{proof}

\debexo
The semialgebraic set
$$\{ (x, y) \in {\Bbb R}^2\ ;\ (y \not= 0 \et x^2+y^2 < 1)
\ou (y = 0 \et 0 < x < 1) \}$$
is open in ${\Bbb R}^2$. Can you write it in the form described in the Finiteness Theorem?
\finexo

\noindent{\bf Remark.} Since $U$ is open in $S$, it can be written as a union of open balls
intersected with $S$, and each open ball is described by a strict polynomial inequality. The main
point in Theorem \ref{finiteness} is that $U$ is a {\em finite} union of semialgebraic subsets
described by conjunctions of strict polynomial inequalities.

\section{Triangulation}

\subsection{Simplicial complexes}

We first recall some definitions concerning simplicial complexes that we shall need. Let  $a_0
,\ldots, a_d$ be points of ${\Bbb R}^n$ which are affine independent (i.e. not contained in
an affine subspace of dimension $d-1$). The $d$-{\em simplex\/\index{simplex}} with vertices
$a_0 ,\ldots, a_d $ is
$$[a_0 ,\ldots, a_d] = \{ x \in  {\Bbb
R}^n\  ;\ \exists \lambda_0 ,\ldots, \lambda_d \in  [0,1]\  \sum_{i=0}^d \lambda_i = 1 
\et x = \sum_{i=0}^d \lambda_i a_i  \}$$

\begin{figure}[ht]
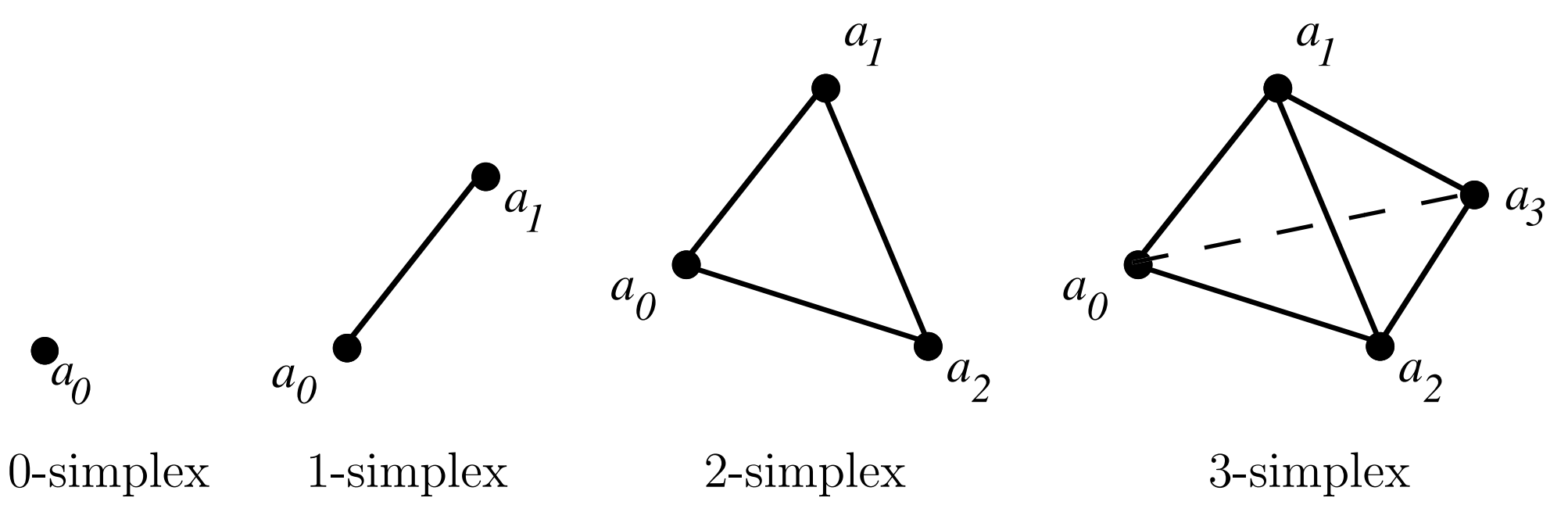

\hspace*{\fill}
\vbox{
\hbox{\hfil\dessin 7.1 mm by 8.2 mm (0-simp)\hfil}
\hbox{0-simplex}}
\hfill
\vbox{
\dessin 25 mm by 23 mm (1-simp)
\hbox to 25 mm{\hfil 1-simplex\hfil}}
\hfill
\vbox{
\dessin 34mm by 36 mm (2-simp)
\hbox to 34mm{\hfil 2-simplex\hfil}}
\hfill
\vbox{
\dessin 43 mm by 36 mm (3-simp)
\hbox to 43mm{\hfil 3-simplex\hfil}}
\hspace*{\fill}
\caption{\label{simplices}Simplices}
\end{figure}

The corresponding {\em open simplex\/} is
$$(a_0 ,\ldots, a_d) = \{ x \in  {\Bbb R}^n\  ; \ \exists \lambda_0 ,\ldots,
\lambda_d \in  (0,1]\  \sum_{i=0}^d \lambda_i = 1  \et x= \sum_{i=0}^d \lambda_i a_i  \}$$ 
We shall denote by $\buildrel\circ\over\sigma$ the open simplex
corresponding to the simplex $\sigma$.  A {\em face\index{face}} of the simplex $\sigma =
[a_0 ,\ldots, a_d]$ is a simplex $\tau  = [b_0 ,\ldots, b_e]$ such that
$$\{b_0 ,\ldots, b_e \} \subset  \{a_0 ,\ldots, a_d \}\;.$$
 
\begin{figure}[ht]
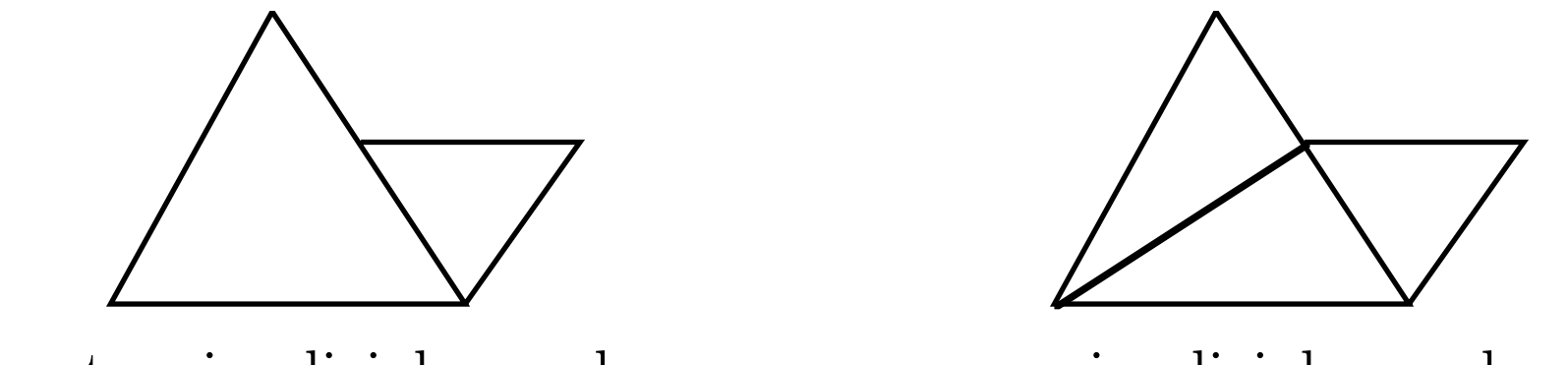

\noindent
\hspace*{\fill}
\vbox{
\hbox to 50mm{\hfil\dessin 33.1mm by 21.1 mm (pas_comp)\hfil}
\hbox to 50mm{\hfil			         not a simplicial complex\hfil}}
\hfill
\vbox{
\hbox to 50mm{\hfil\dessin 33.1mm by 21.1 mm (comp)\hfil}
\hbox to 50mm{\hfil			      a simplicial complex\hfil}}
\hspace*{\fill}
\caption{\label{simpcomp}Simplicial complex}
\end{figure}

A {\em finite simplicial complex\/\index{simplicial complex}} in ${\Bbb R}^n$ is a finite collection
$K =
\{\sigma_1 ,\ldots, \sigma_p \}$ of simplices $\sigma_i \subset  {\Bbb R}^n$ such that, for every
$\sigma_i,
\sigma_j \in  K $, the intersection $\sigma_i \cap \sigma_j$ is a common face of
$\sigma_i$ and $\sigma_j$ (see Figure \ref{simpcomp}).

\begin{figure}[ht]
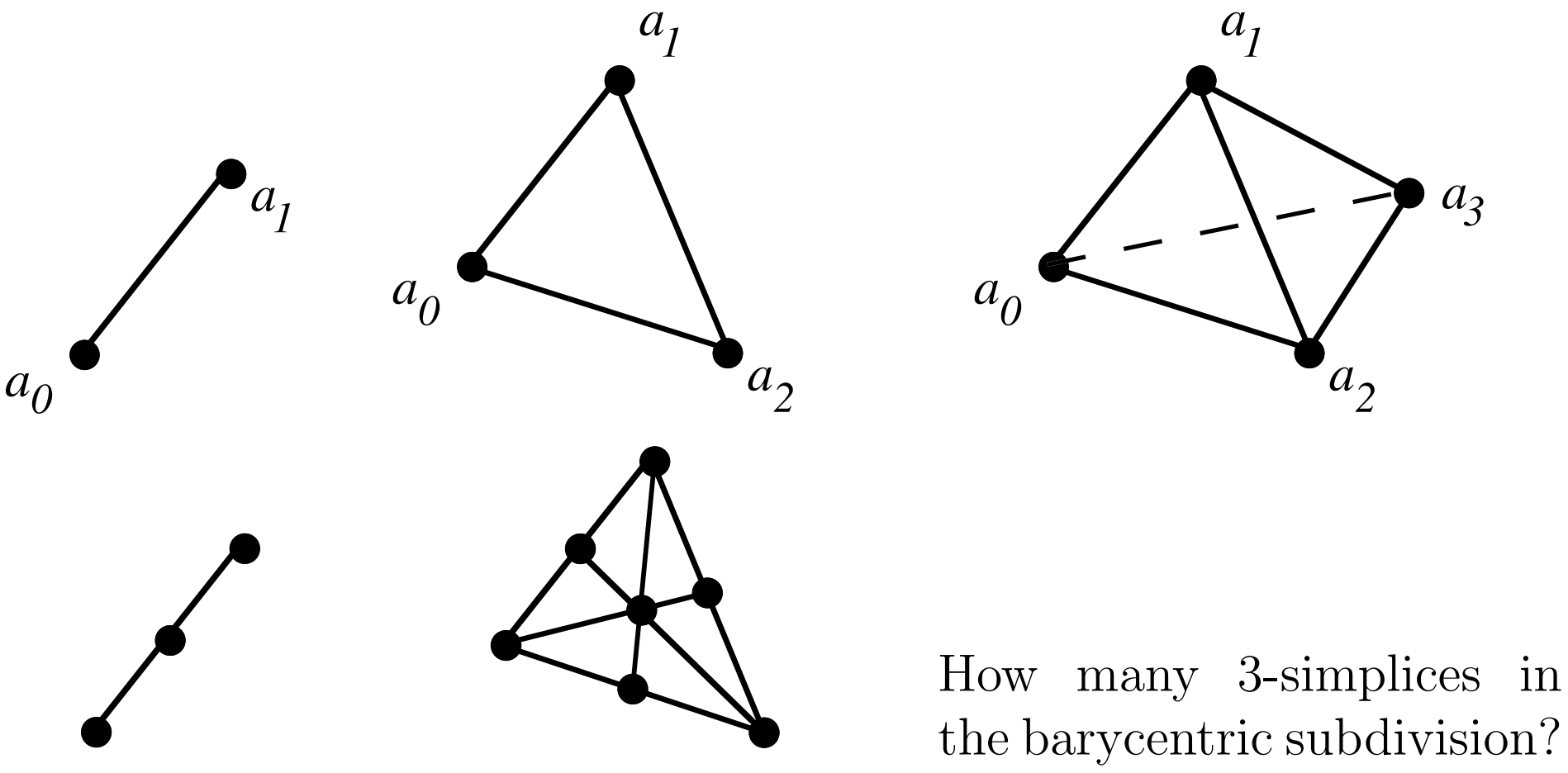

\noindent
\hspace*{\fill}
\dessin 25 mm by 23 mm (1-simp)
\hfill
\dessin 34mm by 36 mm (2-simp)
\hfill
\hbox to 55mm{\hfil\dessin 43mm by 36 mm (3-simp)\hfil}
\hspace*{\fill}\break
\hspace*{\fill}
\hbox to 25mm{\hfil\dessin 15.1 mm by 17.9 mm (1-simp_sub)\hfil}
\hfill
\hbox to 34mm{\hfil\dessin 23.9 mm by 25 mm (2-simp_sub)\hfil}
\hfill
\vbox{\hsize=50mm \noindent How many 3-simplices in the barycentric subdivision?\hfil}
\hspace*{\fill}
\caption{\label{barsub}Barycentric subdivision}
\end{figure}

We set $\vert K\vert  = \bigcup_{\sigma_i \in  K} \sigma_i $; this is a semialgebraic subset of
${\Bbb R}^n$. A {\em polyhedron\/\index{polyhedron}} in ${\Bbb R}^n$ is a subset $P$ of ${\Bbb R}^n
$, such that there exists a finite simplicial complex $K$ in ${\Bbb R}^n$ with $P = \vert K\vert $.
Such a
$K$ will be called a {\em simplicial decomposition} of
$P$. In the following, it will be convenient to agree that if a simplex $\sigma$ belongs to a
finite simplicial complex $K$, then all faces of
$\sigma$ also belong to $K$. With this convention, $\vert K\vert$  is the disjoint union of all
$\buildrel\circ\over\sigma$ for $\sigma \in  K$.
Let $K$ be a finite simplicial complex and, for $\sigma \in  K$, let $\widehat\sigma$ be the
barycenter of $\sigma$. The barycentric subdivision\index{barycentric subdivision} of $K$, denoted by
$K'$, is the finite simplicial complex whose simplices are all 
$[\widehat\sigma_0, \widehat\sigma_1 ,\ldots, \widehat\sigma_d ]$,
such that $\sigma_i$ is a simplex of $K$, for $i = 0 ,\ldots, d$, and $\sigma_i$ is a proper face of
$\sigma_{i+1}$, for $i = 0 ,\ldots, d-1$. This is indeed a finite simplicial complex. The figure
\ref{barsub} shows examples of barycentric subdivisions for complexes which are reduced to a simplex
with its faces.

\subsection{Triangulation of a compact semialgebraic set}  

\begin{thm}\label{triangulation}\index{triangulation}
Let $S \subset  {\Bbb R}^n$ be a compact semialgebraic set, and
$S_1 ,\ldots, S_p$, semialgebraic subsets of $S$. Then there exists a finite simplicial
complex\index{simplicial complex}
$K$ in ${\Bbb R}^n$ and a semialgebraic homeomorphism $h : \vert K\vert  \to S$, such that each
$S_k$ is the image by $h$ of a union of open simplices of $K$.
\end{thm}

\begin{proof}
The proof is by induction on $n$. For $n = 1$, we can take $\vert K\vert  =
S$, and the only open simplices are points and bounded open intervals.
We now prove the theorem for $n>1$, assuming that it holds true for $n-1$. Up to a linear change
of variables (cf. Proposition \ref{chgvar}), we can assume that $S$ and all $S_k$ are unions of
cells of a c.a.d. satisfying the properties of the generalized Thom lemma \ref{LTG},
associated to a family of polynomials $(P_{i,j})$. In particular, we have a partition of ${\Bbb
R}^{n-1}$ into finitely many connected semialgebraic sets $C_i$ and continuous semialgebraic
functions
$\xi_{i,1} <\ldots<\xi_{i,\ell_i} : C_i \to {\Bbb R}$, which describe the real roots of the
polynomials $P_{n,1}(x,X_n),\,\ldots,P_{n,s_n}(x,X_n)$ as functions of $x\in C_i$. We know that
$S$ and the $S_k$ are unions of graphs of $\xi_{i, j}$ and bands of the cylinders $C_i
\times {\Bbb R}$ bounded by these graphs. Denote by $\pi : {\Bbb R}^{n} \to {\Bbb
R}^{n-1}$ the projection onto the space of the first $n-1$ coordinates. The set $\pi(S)$
is compact, semialgebraic, and it is the union of some $C_i$; similarly, each
$\pi(S_K)$ is the union of some $C_i$. By the inductive assumption, there is a triangulation $g :
\vert L\vert  \to \pi(S) $, where $L$ is a finite simplicial complex in ${\Bbb R}^{n-1}$ and
$g$ a semialgebraic homeomorphism, such that each $C_i \subset 
\pi(S)$ is the union of images by $g$ of open simplices of $L$. We can actually subdivide the  $C_i
\subset \pi(S)$ and assume that these $C_i$ are of the form
$g({\buildrel\circ\over\tau_i})$, where $\tau _i$ is a simplex of $L$. 
\par
We shall now construct, for
each $\tau$  of $L$, a finite simplicial complex $K_\tau$ in ${\Bbb R}^n$ and a semialgebraic
homeomorphism
$$h_\tau  :
\vert K_\tau \vert  \longrightarrow \clos(S \cap(g({\buildrel\circ\over\tau})\times{\Bbb R}))\;.$$
We fix $\tau $ in $L$, say $\tau  = [b_0 ,\ldots, b_d]$. Let $\xi : 
g({\buildrel\circ\over\tau}) \to {\Bbb R}$ be one of the functions of the c.a.d., whose graph is
contained in
$S$. By Lemma \ref{prolongement}, $\xi$ has a continuous extension $\overline\xi$ defined on
$\clos( g({\buildrel\circ\over\tau})) = g(\tau)$. Set $a_i = (b_i, \overline\xi(g(b_i)))
\in  {\Bbb R}^n$, for $i = 0,\ldots,d$. Let $\sigma_\xi$ be the simplex $[a_0 ,\ldots, a_d] \subset 
{\Bbb R}^n$.
This $\sigma_\xi$ will be one of the simplices of $K_\tau$, and we define $h_\tau$ on
$\sigma_\xi$ by (see Figure \ref{triang1})
$$h_\tau(\lambda _0 a_0 +\cdots+ \lambda_d a_d) = ( y,
\overline\xi(y))\;,\quad \hbox{\rm where } y=g(\lambda_0b_0+\cdots+\lambda_db_d)\;.$$

\begin{figure}[ht]
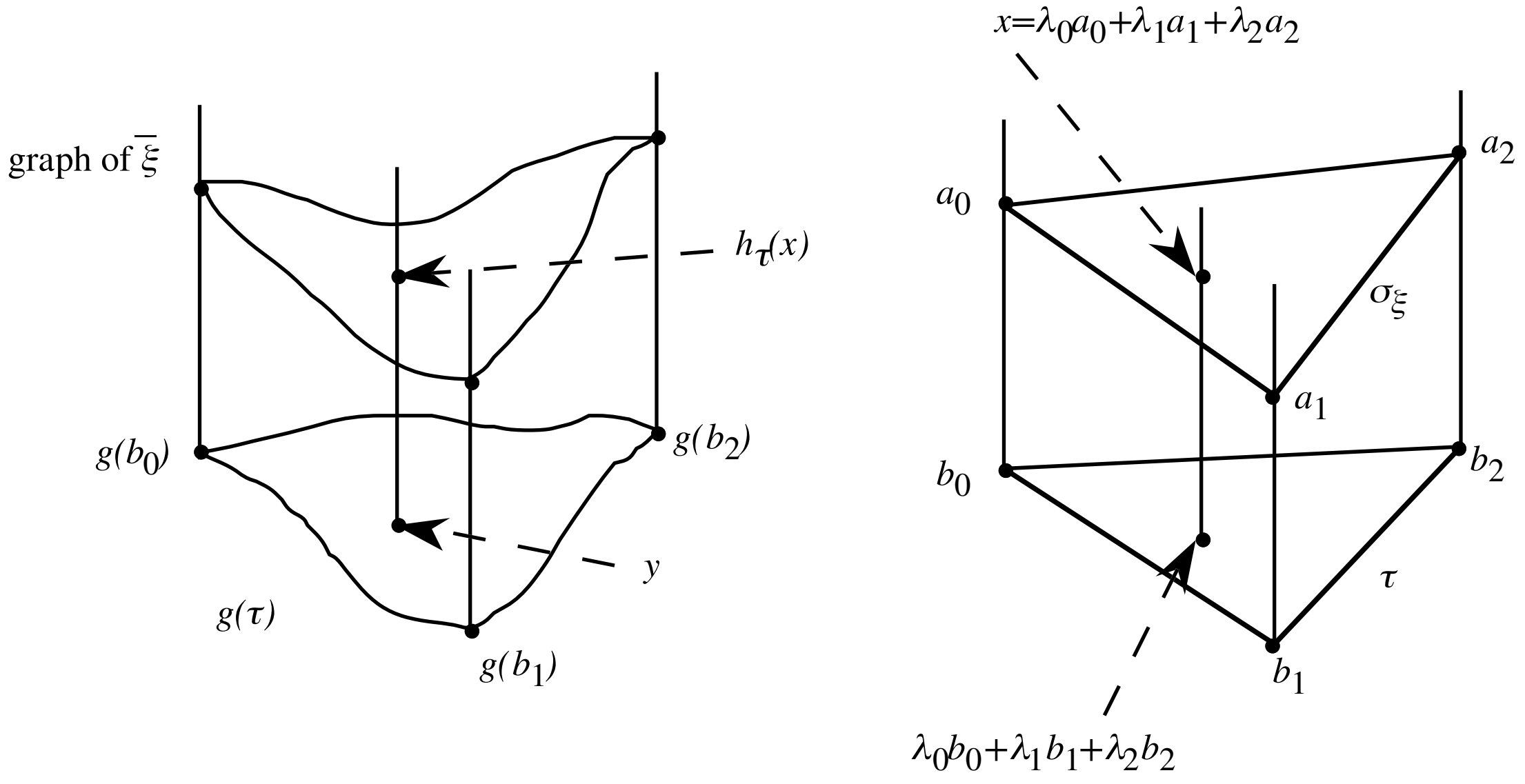

\dessin 148mm by 77mm (triang1eng)
\caption{Construction of $K_\tau$ and $h_\tau$: the case of a graph\label{triang1}}
\end{figure}

If $\xi' = g({\buildrel\circ\over\tau}) \to {\Bbb R}$ is another function of the c.a.d.\ whose
graph is contained in $S$, we define $\sigma_{\xi'} = [a'_0 ,\ldots, a'_d] $ in the same way. We do
not want $\sigma_{\xi'}$ and $\sigma_\xi$ to coincide: at least one of the $a'_i$ must be different
from $a_i$. Also, if the restrictions of $\overline\xi$ and $\overline{\xi'}$ to a face $\rho$ of
$\tau $ are different, the values of $\overline\xi$ and $\overline{\xi'}$ must differ for at least
one of the vertices $b_i$ of $\rho$ (in order that the corresponding $a_i$ and $a'_i$ are distinct).
So we should have the following property:
\begin{itemize}
\item[$\dag$] for every simplex $\tau $ of $L$, if $\xi$ and $\xi'$ are distinct functions
$g({\buildrel\circ\over\tau})\to {\Bbb R}$ of the c.a.d., there is at least one vertex $b$ of
$\tau$  such that $\overline\xi(g(b)) \not=\overline{\xi'}(g(b))$.
\end{itemize}
Observe that if $\tau _1$ is a simplex of the barycentric subdivision $L'$ of $L$, with
${\buildrel\circ\over\tau_1}\subset  {\buildrel\circ\over\tau}$, then the barycenter $\widehat\tau$
of $\tau $ is a vertex of $\tau _1$, and $\xi \not= \xi'$ implies
$\xi(g(\widehat\tau))\not=\xi'(g(\widehat\tau))$. Hence, replacing $L$ with its
barycentric subdivision $L'$, we can assume that property  $\dag$ holds.

\begin{figure}[ht]
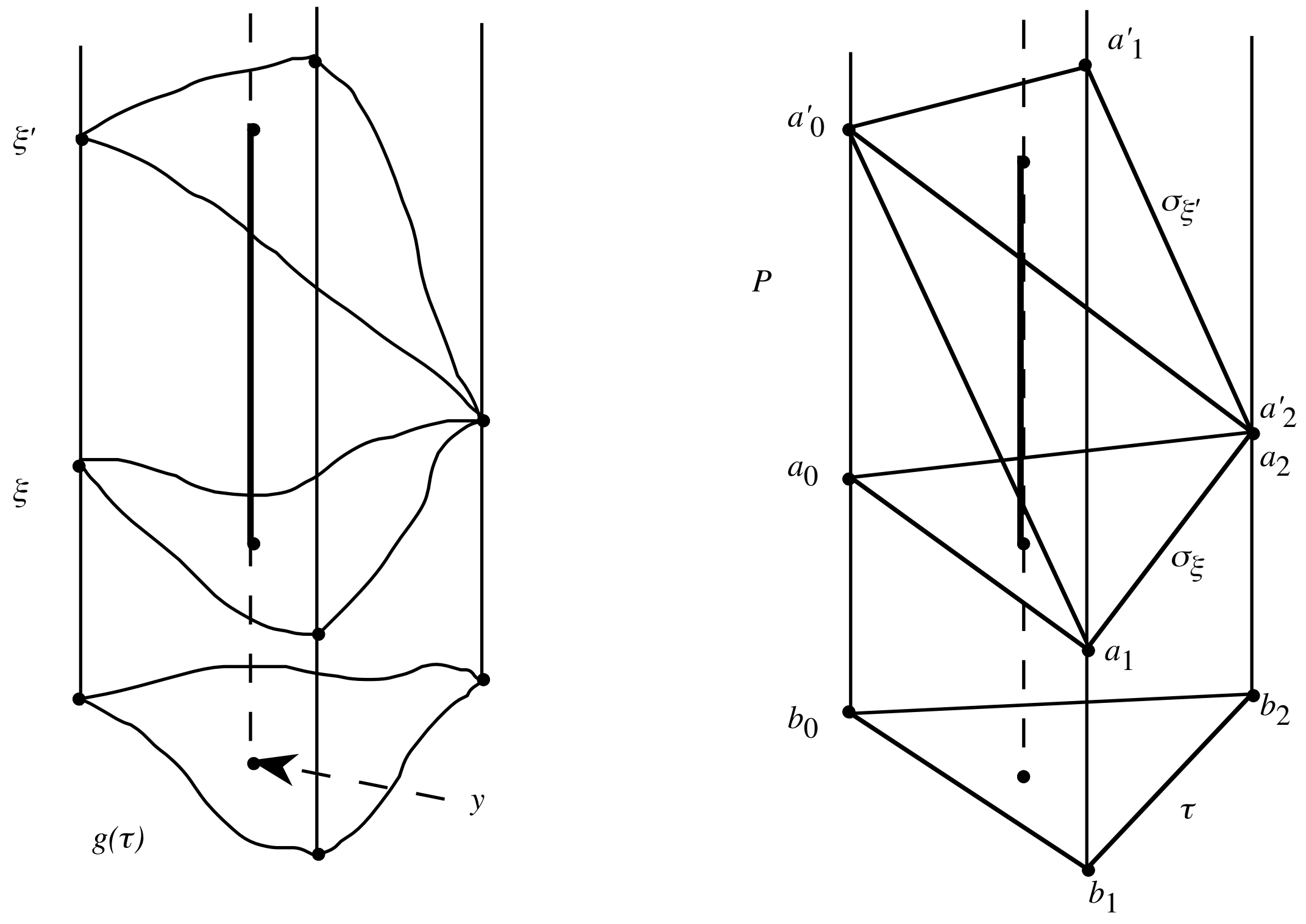
 {\noindent\dessin 141mm by 101mm (triang2eng)}
\caption{Construction of $K_\tau$ and $h_\tau$: the case of a band\label{triang2}}
\end{figure}

We now consider the case of a band $B$ contained in $S$ (see Figure \ref{triang2}). Since
$S$ is compact, the band
$B$ is bounded from below and from above  by the graphs of two consecutive functions
$$\xi <
\xi' : g({\buildrel\circ\over\tau}) \longrightarrow {\Bbb R}\;.$$
Let $P$ be the polyhedron which
is the part of the cylinder
$\tau 
\times  {\Bbb R}$ delimited by the simplices $\sigma_\xi$ and $\sigma_{\xi'}$. The polyhedron
 $P$ has the simplicial decomposition
$$P = \bigcup_{\scriptstyle 0 \leq i
\leq d \atop \scriptstyle a_i \not= a'_i}   [a'_0 ,\ldots, a'_i , a_i ,\ldots, a_d]\;,.$$ 
These
simplices and their faces will belong to $K_\tau$, and we define $h_\tau $ on
$P$ so that it sends linearly the segment $[\lambda_0
a_0 +\cdots+ \lambda_d a_d , \lambda_0 a'_0 +\cdots+ \lambda_d a'_d]$  onto the segment  $[( y,
\overline\xi( y)), ( y, \overline{\xi'}( y))]$, where $y = g(\lambda _0
b_0 +\cdots+\lambda_d b_d)$. By property  $\dag$, the first segment is reduced to a point if and
only if it is the case for the second.

We have constructed $K_\tau$ and $h_\tau$ for every simplex $\tau $ de $L$. We shall now check
that the $K_\tau$  and $h_\tau$  can be glued together to give $K$ and the triangulation
$h : |K|\to S$. It is sufficient to check the gluing property for a simplex $\tau$  and one of its
faces $\rho$. First observe that if $\sigma_\eta $ is a simplex of $K_\rho$ which meets  $\vert
K_\tau \vert $ and is sent by $h_\rho$ onto the closure of the graph of the function $\eta  : g(
{\buildrel\circ\over\rho}) \to {\Bbb R}$ of the c.a.d., then $\sigma_\eta$ is a simplex of
$K_\tau$. Indeed, by property 2 of Lemma \ref{prolongement}, $\eta $ coincides with
$\overline\xi$ on $g({\buildrel\circ\over\rho})$, for some function $\xi
:  g({\buildrel\circ\over\tau}) \to {\Bbb R}$ of the c.a.d. such that
$\sigma_\xi$ is a simplex of $K_\tau$; then $\sigma_\eta$  is a face of $\sigma_\xi$. It follows
also that $h_\tau$  and $h_\rho$ coincide on $\vert K_\tau \vert \cap \vert K_\rho\vert$. It
remains to check that the simplicial decomposition of the polyhedron $P$ in $\tau  \times  {\Bbb
R}$ (cf. above) induces the simplicial decomposition of the polyhedron $P\cap (\rho
\times  {\Bbb R})$. This is the case if we have chosen a total ordering of the set of vertices of
$L$, and if we make the simplicial decomposition
$$P = \bigcup_i [a'_0 ,\ldots, a'_i , a_i ,\ldots, a_d]\;,$$  
where $b_0 <\ldots< b_i<\ldots< b_d$ for the chosen ordering.
\end{proof}

\subsection{The curve selection lemma}

The triangulation\index{triangulation} theorem allows one to give a short proof of the following.

\begin{thm}[Curve selection lemma]\label{petitchemin}
Let\index{curve selection} $S \subset  {\Bbb R}^n$ be a semialgebraic set. Let $x \in 
\clos(S)$, $x \not\in S$. Then there exists a continuous semialgebraic mapping $\gamma : [0,1] \to
{\Bbb R}^n$ such that $\gamma(0) = x$ and $\gamma((0,1])\subset  S$.
\end{thm}

\begin{proof}	Replacing $S$ with its intersection with a ball with center $x$ and radius $1$, we
can assume $S$ bounded. Then $\clos(S)$ is a compact semialgebraic set. By the triangulation
theorem, there is a finite simplicial complex $K$ and a semialgebraic homeomorphism $h : \vert
K\vert  \to \clos(S)$, such that $x = h(a)$ for a vertex $a$ of $K$ and $S$ is the union of some
open simplices of $K$. In particular, since $x$ is in the closure of $S$ and not in $S$, there is a
simplex $\sigma$ of $K$ whose $a$ is a vertex, and such that $h( {\buildrel\circ\over\sigma})
\subset  S$. Taking a linear parametrization of the segment joining $a$ to the barycenter of
$\sigma$, we obtain $\delta : [0,1]\to \sigma$ such that $\delta(0) = a$ and $\delta((0,1]) \subset 
\buildrel\circ\over\sigma$. Then $\gamma = h\circ\delta$ satisfies the property of the theorem.
\end{proof}

\debexo\label{arcconnect}\index{connected component} Show that a connected semialgebraic set $A$ is
{\em semialgebraically arcwise connected}: for every points $a,b$ in $A$, there exists a continuous
semialgebraic mapping
$\gamma : [0,1]\to A$ such that $\gamma(0)=a$ and $\gamma(1)=b$.
\finexo

\section{Dimension}

\subsection{Dimension via c.a.d.}

The c.a.d.\index{cylindrical algebraic decomposition}\ allows one to decompose a semialgebraic set
$S$ as a finite union of cells\index{cell} semialgebraically homeomorphic (and even diffeomorphic) to
open hypercubes
$(0,1)^d$. This leads to the definition of the dimension of $S$ as the maximum of these $d$.

\begin{prop}\label{dimension} Let $S \subset  {\Bbb R}^n$ 
be a semialgebraic set, and let $S = \bigcup_{i = 1}^p C_i$ be a decomposition of $S$ into a
disjoint union of semialgebraic subsets $C_i$, each semialgebraically diffeomorphic
to $(0,1)^{d_i}$. The {\em  dimension\index{dimension}} of $S$ is, by definition, $d = \max \{ d_i \
;\  i = 1,\ldots , p
\}$. This dimension is independent of the decomposition.
\end{prop}

\begin{proof} Assume that $S = \bigcup_{i = 1}^p C_i$ and $S = 
\bigcup_{j = 1}^q D_j$  are both decompositions of $S$ into a disjoint union of semialgebraic
subsets, each semialgebraically diffeomorphic to an open hypercube. 
By corollary \ref{stratification}, there is a semialgebraic stratification
$S=\bigcup_{k = 1}^r\Sigma_k$, which is a common refinement of the two decompositions. Each $C_i$
and each $D_j$  is a finite union of strata $\Sigma_k$. Let us compare the dimensions
of
$C_i$,
$D_j$ and $\Sigma_k$, as submanifolds of ${\Bbb R}^n$.  If
$\Sigma_k$ is contained in $C_i$, the dimension of $\Sigma_k$ is not greater than the dimension of
$C_i$. On the other hand, if $\Sigma_k$ is of maximal dimension among the strata contained in $C_i$,
then $\Sigma_k$ is open in $C_i$ and, hence, has the same dimension as $C_i$. Indeed, the closure of
$C_i\setminus\Sigma_k$ contains only $C_i\setminus\Sigma_k$ and strata of smaller dimension than
$\Sigma_k$; therefore, it is disjoint from $\Sigma_k$. It follows that the maximum of the
dimensions of the $C_i$ is equal to the maximum of the dimensions of the $\Sigma_k$, which is,
for the same reason, equal to the maximum of the dimensions of the $D_j$.\end{proof}

The following properties of the dimension of a semialgebraic set are obvious: the dimension of
the union of finitely many semialgebraic sets is the maximum of the dimensions of these
semialgebraic sets; the dimension of the cartesian product of semialgebraic sets is the sum of their
dimensions. The dimension of a semialgebraic set behaves well with respect to the topological
closure.

\begin{prop} Let $A \subset  {\Bbb R}^n$ be a semialgebraic set. Then
\begin{enumerate}\shitem
\item $ \dim(\clos (A)) = \dim A$,
\item	$\dim(\clos (A) \setminus  A) < \dim A$.
\end{enumerate}
\end{prop}

\begin{proof} Both properties follow from the definition of dimension and the fact that the
closure of a stratum is the union of this stratum and strata of smaller dimensions 
(cf. \ref{stratification}).
\end{proof}

The dimension is invariant by a semialgebraic homeomorphism.

\begin{lem}
Let $A \subset  {\Bbb R}^{n+k}$ be a semialgebraic set, $\pi :{\Bbb R}^{n+k}\to{\Bbb R}^n$ the
projection on the space of the first $n$ coordinates. Then
$\dim(\pi (A)) \leq  \dim(A)$. Moreover, if the restriction of $\pi$ to $A$ is one-to-one, 
then $\dim(\pi (A)) = \dim  A$.
\end{lem}

\begin{proof}
If $k = 1$ and  $A$ is either the graph of a function or a band of a c.a.d. of ${\Bbb
R}^{n+1}$, the lemma is obvious. If $A$ is any semialgebraic subset of ${\Bbb R}^{n+1}$, it is the
union of cells of a c.a.d.; hence, the lemma still holds true in this case. We prove the case of
$k>1$ by an easy induction.
\end{proof}

\begin{thm} Let $S$ be a semialgebraic subset of ${\Bbb R}^n$, and $f : S 
\to  {\Bbb R}^k$ a semialgebraic mapping (not necessarily continuous). Then $\dim f(S) 
\leq  \dim S$. If $f$ is one-to-one, then $\dim f(S) = \dim S$.
\end{thm}

\begin{proof}
Let $A\subset  {\Bbb R}^{n+k}$  be the graph of $f$. From the preceding lemma, it follows that
$\dim(S) = \dim(A)$ and $\dim( f(S)) \leq  \dim(A)$,  with moreover $\dim( f(S)) =  \dim(A)$ if $f$
is one-to-one.
\end{proof}

\debexo
Let $S\subset{\Bbb R}^n$ be a semialgebraic set, $x\in S$. Show that there exist a neighborhood
$V$ of $x$ in ${\Bbb R}^n$ and a nonnegative integer $d$ such that, for every semialgebraic
neighborhood $W\subset V$ of $x$ in ${\Bbb R}^n$, $\dim(W\cap S)=d$. The integer $d$ is called the
{\em dimension of $S$ at $x$} and denoted by $\dim_x S$. Show that 
$$\dim S=\max\{\dim_x S\;;\ x\in S\}\;.$$
Show that $\{x\in S\;;\ \dim_x S=\dim S\}$ is a closed semialgebraic subset of $S$.
\finexo

\subsection{Dimension of algebraic sets}

We shall compare the dimension of semialgebraic sets, defined via decomposition, with the
dimension of algebraic sets. We recall the results concerning the dimension of algebraic sets
that we shall need (cf.\ \cite{S}, Chapter 1).

For $A \subset  {\Bbb R}[ X_1 ,\ldots , X_n]$, we denote by ${\cal Z}(A) 
\subset  {\Bbb R}^n$ the common zeroset of all polynomials of $A$:
$${\cal Z}(A) = \{ x \in  {\Bbb R}^n\ ; \ \forall  P \in  A \ P(x) = 0 \}\;.$$
An algebraic subset of ${\Bbb R}^n$ is a subset of the form ${\cal Z}(A)$, for some subset
$A$ of ${\Bbb R}[ X_1 ,\ldots , X_n]$. If $I$ is the ideal of
${\Bbb  R}[ X_1 ,\ldots , X_n]$ generated by $A$, then ${\cal Z}(I) = {\cal 
Z}(A)$. Since every ideal of  ${\Bbb R}[ X_1 ,\ldots , X_n]$ is generated by finitely many
polynomials, an algebraic set is the common zeroset of finitely many polynomials. Since, for
$x\in{\Bbb R}^n$, $P_1(x) =\ldots = P_s(x) = 0$ 
is equivalent to $(P_1^2+\cdots + P_s^2)(x) = 0$, an algebraic subset of 
${\Bbb R}^n$ can always be described by {\em one} equation (this is not the case for complex
algebraic sets).

For $S \subset  {\Bbb R}^n$, we denote by ${\cal I}(S) \subset  {\Bbb 
R}[ X_1 ,\ldots , X_n]$ the subset of polynomials which vanish on $S$:
$${\cal I}(S) = \{ P\in  {\Bbb R}[ X_1 ,\ldots , X_n]\  ; \  \forall  x 
\in  S \  P(x) = 0 \}\;.$$
${\cal I}(S)$ is an ideal of ${\Bbb R}[ X_1 ,\ldots , X_n]$. A subset $V \subset  {\Bbb R}^n$ is
algebraic if and only if
$V = {\cal Z}({\cal I}(V))$.
The quotient ring ${\cal P}(S) =  {\Bbb R}[ X_1 ,\ldots , X_n] 
/ {\cal I}(S)$ is called {\em the ring of polynomial functions on $S$}. Indeed, it can be
identified with the ring of functions $S \to   {\Bbb R}$ which are the restriction of a polynomial.
A nonempty algebraic set $V$ is said to be {\em irreducible\index{irreducible}} if it cannot be
written as the union of two algebraic sets strictly contained in $V$. If $V$ is irreducible, ${\cal
P}(V)$ is an integral domain, and we shall denote by ${\cal K}(V)$ its field of fractions (the field
of rational fractions on $V$).\par

The {\em dimension\index{dimension}} of an algebraic set $V$ is, by definition, the Krull dimension
of the ring
${\cal P}(V)$, i.e. the maximal length of chains of prime ideals in 
${\cal P}(V)$~:
$\dim({\cal P}(V))$ is the maximum of $d$ such that there exist  prime ideals ${\frakt p}_0, 
{\frakt p}_1 ,\ldots , {\frakt p}_d$ of  
${\cal P}(V)$, with ${\frakt p}_0 \strictsubset {\frakt 
p}_1 \strictsubset\ldots  \strictsubset {\frakt p}_d$.
If $V$ is irreducible, this dimension is equal to the transcendence degree of the field 
${\cal K}(V)$ over ${\Bbb R}$ (i.e. ${\cal K}(V)$ is an algebraic extension of the field of rational
fractions ${\Bbb  R}(T_1 ,\ldots , T_d)$, where $d = \dim(V)$). An algebraic set  
$V$ has a unique decomposition as a union of finitely many irreducible algebraic subsets $V_1 ,\ldots
, V_p$, where $V_i \not\subset V_j$ for  $i \not= j$. The $V_i$ are called the {\em 
irreducible components} of $V$.
If $W_1 ,\ldots , W_k$ are algebraic sets, so is $W_1 \cup  \ldots  
\cup  W_k$, and $\dim(W_1 \cup  \ldots  \cup  W_k) = \max( 
\dim(W_i))$. In particular, the dimension of an algebraic set is the maximum of the dimensions of
its irreducible components.\par

 If $S \subset  {\Bbb R}^n$ is any subset,
${\cal  Z}({\cal I}(S))$ is the smallest algebraic subset of ${\Bbb 
R}^n$ containing $S$. It is called the {\em Zariski closure\index{Zariski closure} of $S$}, and it
will be denoted by
$\overline{S}^Z$. The algebraic subsets of ${\Bbb R}^n$ are the closed sets of a topology on ${\Bbb
R}^n$, which is called the Zariski topology, and $\overline{S}^Z$ is the closure of
$S$ for this topology. The Zariski topology is coarser than the usual topology, and it is not
separated.

\begin{thm} Let $S \subset  {\Bbb R}^n$ be a semialgebraic set. Its dimension\index{dimension} as a
semialgebraic set (cf. \ref{dimension}) is equal to the dimension, as an algebraic set, of its
Zariski closure
$\overline{S}^Z$. In particular, if $V \subset  {\Bbb R}^n$ is an algebraic set, its dimension as a
semialgebraic set is equal to its dimension as an algebraic set (i.e. the Krull dimension of ${\cal
P}(V)$).
\end{thm}

\begin{proof}
If $S=\bigcup_{i=1}^p C_i$, then  $\overline{S}^Z= 
\bigcup_{i=1}^p \overline{C_i}^Z$. Hence, it is sufficient to prove the theorem for a cell $C\subset
{\Bbb R}^n$ of a c.a.d.. The proof is by induction on $n$.\par
If $n = 1$, either $C$ is a point and $\overline{C}^Z$ is equal to this point, or $C$ is a nonempty
open interval and $\overline{C}^Z= {\Bbb R}$. We have algebraic dimension $0$ in the first case and
$1$ in the second case.
Let $n>1$, and assume the theorem proved for $n-1$. Let $\pi  : {\Bbb R}^n \to   {\Bbb R}^{n-1}$ be
the projection on the space of the first $n-1$ coordinates. Then
$\pi (C) = D$ is a cell of the c.a.d.,  semialgebraically homeomorphic to $(0,1)^d$.
By the inductive assumption, $\dim \overline{D}^Z= d $ (dimension as algebraic set). 
We have to consider two cases.
\begin{enumerate}\shitem
\item $C$ is the graph of a function $\xi  : D \to  {\Bbb R}$ of the c.a.d., and, hence, $C$ is
semialgebraically diffeomorphic to $(0,1)^d$. There is a polynomial
$P\in  {\Bbb R}[ X_1 ,\ldots , X_n] $ such that, for every $x \in  D,  P( x, X_n)$ is not
identically zero and $P( x, \xi (x)) = 0$. Let $Z = {\cal Z}(P) \subset  {\Bbb R}^n$. Let 
$\overline{D}^Z= V_1 \cup \ldots \cup  V_p$ be the decomposition into irreducible components
(actually, it can be shown that $\overline{D}^Z$ is irreducible). We have
$$\overline{C}^Z\subset  Z \cap  ( (V_1 \times  {\Bbb R}) \cup  (V_2 
\times  {\Bbb R}) \cup \ldots \cup  (V_p \times  {\Bbb R} )),$$
and $Z \cap  (V_i \times  {\Bbb R}) \strictsubset  V_i \times  {\Bbb 
R}$, for $i = 1,\ldots , p$.
We use the following facts concerning algebraic sets and their dimensions:
\begin{itemize}\shitem
		\item $V_i \times  {\Bbb R}$ is irreducible, since both $V_i $ 
and ${\Bbb R}$ are irreducible.
		\item $\dim(V_i \times  {\Bbb R}) = \dim V_i + 1$ 
(algebraic dimension).
		\item If $W$ is an irreducible algebraic set and $V \strictsubset W$ a proper algebraic subset,
then $\dim V < \dim W$.
\end{itemize}
It follows that $$\dim Z \cap  (V_i \times  {\Bbb R}) < \dim V_i 
\times  {\Bbb R} = \dim V_i + 1\;,$$ 
therefore $\dim  \overline{C}^Z \leq  \dim \overline{D}^Z= d$.
It remains to prove the reverse inequality. Let $W$
be an irreducible component of $\overline{C}^Z$. Then $Y = 
\overline{\pi(W)}^Z$ is irreducible. Indeed, if $Y = F_1 \cup  F_2$, 
where $F_1$ and $F_2$ are algebraic sets, we have $W \subset  \pi ^{-
1}(F_1) $ or $W \subset  \pi ^{-1}(F_2)$ since $W$ is irreducible; hence, $Y \subset  F_1$ or $Y 
\subset  F_2$.
The projection $\pi$  induces an injective homomorphism from ${\cal 
P}(Y) = {\cal P}(\pi (W))$ into ${\cal P}(W)$ and, hence, a field homomorphism ${\cal K}(Y) \to 
{\cal K}(W)$. It follows that $\dim(Y) \leq  \dim(W)$. Since $\overline{D}^Z$ is the union of these
$Y$, we have $\dim \overline{D}^Z \leq  \dim  \overline{C}^Z$.
	\item $C$ is a band of the cylinder $D \times  {\Bbb R} $. Then $C$ is semialgebraically
diffeomorphic to $(0,1)^{d+1}$. In this case,
$\overline{C}^Z\supset D\times  {\Bbb R}$, therefore 
$\dim\overline{C}^Z = \dim\overline{D}^Z+ 1 = d + 1$.
\end{enumerate}
The proof is completed.\end{proof}

\chapter[Families. Uniform bounds]{Families of semialgebraic sets.
				Uniform bounds}

\section{Semialgebraic triviality of families}
\subsection{Hardt's theorem}

	Let $A \subset  {\Bbb R}^n$ be a semialgebraic set, 
defined by a boolean combination of sign conditions on polynomials $P_1 ,\ldots , P_q$. Construct
a c.a.d.\index{cylindrical algebraic decomposition}\ of ${\Bbb R}^n$ adapted to $P_1,\ldots,P_q$.
The set
$A$ is a union of graphs and bands in cylinders  
$C_i \times  {\Bbb R}$, where ${\Bbb R}^{n-1}=C_1 \cup\ldots\cup C_r$  is a finite semialgebraic
partition.  Each $A\cap  (C_i\times  {\Bbb R})$ is semialgebraically homeomorphic to a product $C_i
\times  F_i$, where
$F_i$ is a semialgebraic subset of ${\Bbb R}$: one can take for instance $F_i = p^{-1}(b_i)$, where
$p : A  \to  {\Bbb R}^{n-1}$ is the restriction of the projection onto the space of the $n-1$ first
coordinates, and $b_i$, a point chosen in $C_i$. Hence,
we have decomposed the target space ${\Bbb R}^{n-1}$ as the disjoint union of finitely many
semialgebraic subsets $C_i$, such that $p$ is semialgebraically trivial over each 
$C_i$ in the following sense.\par
A continuous semi-algebraic mapping $p:A\to{\Bbb R}^k$ is said to be {\em semialgebraically
trivial\index{trivial} over a semialgebraic subset $C\subset {\Bbb R}^k$} is there is a
semialgebraic set
$F$ and a semialgebraic homeomorphism
$h:p^{-1}(C)\to C\times F$, such that the composition of $h$ with the projection $C\times F\to C$
is equal to the restriction of $p$ to $p^{-1}(C)$.
\setlength{\unitlength}{1mm}
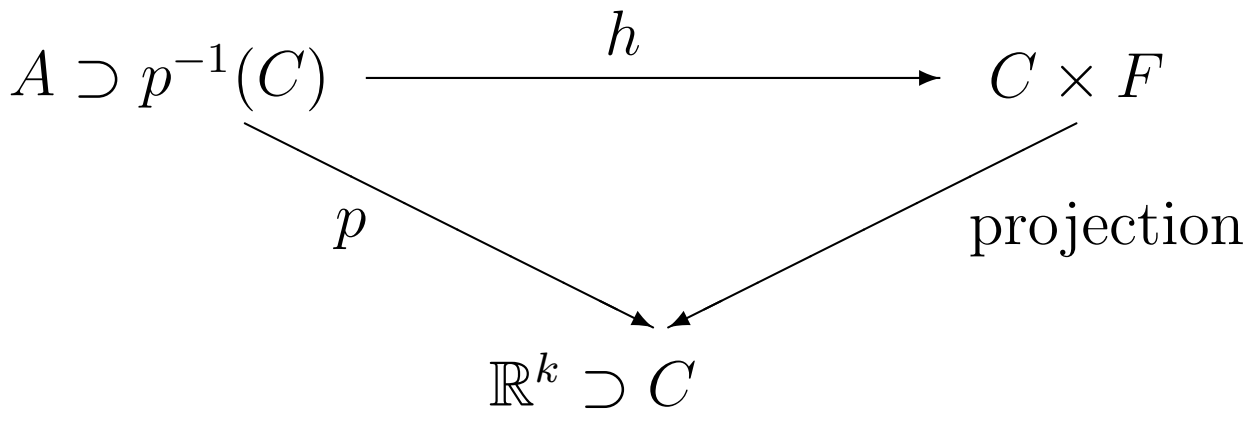
$$\begin{picture}(60,30)
\put(0,25){\makebox(0,0){$A \supset p^{-1}(C)$}}
\put(13,25){\vector(1,0){38}}
\put(30,28){\makebox(0,0){$h$}}
\put(60,25){\makebox(0,0){$C \times  F$}}
\put(5,22){\vector(2,-1){27}}
\put(12,15){\makebox(0,0){$p$}}
\put(60,22){\vector(-2,-1){27}}
\put(62,15){\makebox(0,0){\rm projection}}
\put(28,5){\makebox(0,0){${\Bbb R}^k \supset C$}}
\end{picture}$$
 The homeomorphism $h$ is called a semi-algebraic
trivialization of $p$ over $C$. We say that the trivialization $h$ is compatible with a
semialgebraic subset $B\subset A$ if there is a semialgebraic subset $G\subset F$ such that
$h(B\cap p^{-1}(C))=C\times G$.\par
The above property  of the projection mapping
$p:A\to {\Bbb R}^{n-1}$ holds actually for every continuous semialgebraic mapping.

\begin{thm} [Hardt's semialgebraic triviality]\label{Hardt}\index{Hardt's theorem}
Let $A \subset  {\Bbb R}^n$ be a semialgebraic set and $p : A 
\to  {\Bbb R}^k$, a continuous semi-algebraic mapping. There is a finite semialgebraic partition 
of ${\Bbb R}^k$ into $C_1 ,\ldots, C_m$ such that  $p$ is
semialgebraically trivial\index{trivial} over each $C_i$. Moreover, if $B_1 ,\ldots , B_q$ are
finitely many semialgebraic subsets of $A$, we can ask that each trivialization $h_i:p^{-1}(C_i)\to
C_i\times F_i$ is compatible with all $B_j$.
\end{thm}
In particular, if $b$ and $b'$ are in the same $C_i$, then 
$p^{-1}(b) $ and $p^{-1}(b')$ are semialgebraically homeomorphic, since they are both
semialgebraically homeomorphic to $F_i$.
Actually we can take for $F_i$ a fiber 
$p^{-1}(b_i)$, where $b_i$ is a chosen point in $C_i$, and we ask in this case that 
 $h_i(x) = (x, b_i)$ for all $x\in p^{-1}(b_i)$.
\par

For the proof of Hardt's theorem, we refer to \cite{BR} or \cite{BCR}. Hardt's theorem for
o-minimal structures is proved in \cite{D} and in the lecture notes \cite{Co} in the same
collection.\par

We can easily derive from Hardt's theorem\index{Hardt's theorem} a useful information about the
dimensions\index{dimension} of the fibers of a continuous semialgebraic mapping. We keep the
notation of the theorem. For every
$b \in C_i$  :
$$\dim p^{-1}(b)  =  \dim F_i  =  \dim p^{-1}(C_i) - \dim C_i  
\leq  \dim A - \dim C_i\; .$$
From this observation follows: 

\begin{cor} Let $A \subset  {\Bbb R}^n$ be a semialgebraic set and $f : A \to  {\Bbb R}^k$, a
continuous semialgebraic mapping. For
$d \in {\Bbb N}$, the set
$$\{ b \in {\Bbb R}^k\ ;\  \dim( p^{-1}(b)) = d \}$$
is a semialgebraic subset of ${\Bbb R}^k$ of dimension not greater than  $\dim A - d$.
\end{cor}

\debexo
Let $p:{\Bbb R}^{n+1}\to{\Bbb R}^n$ be the projection on the space of the first $n$ coordinates.
Let $A$ be a semialgebraic subset of ${\Bbb R}^{n+1}$, of dimension $n$. Show that there is an
integer $N\geq 0$ such that
$$\{t\in{\Bbb R}^n\ ;\ p^{-1}(t)\cap A\ \mbox{\rm has exactly}\ N\ \mbox{\rm elements }\}$$ 
is a semialgebraic set of dimension $n$ and
$$\{t\in{\Bbb R}^n\ ;\ p^{-1}(t)\cap A\ \mbox{\rm has} > N\ \mbox{\rm elements }\}$$ 
is a semialgebraic set of dimension $< n$ or empty.
\finexo

\subsection{Local conic structure of semialgebraic sets}

Let $A$ be a semialgebraic subset of ${\Bbb R}^n$ and $a$, a nonisolated point of $A$: for every
$\varepsilon  > 0$ there is $x \in  A$, $x \not= a$,  such that $\Vert x-a\Vert  <\varepsilon $.
Let
$\overline{B}(a, \varepsilon )$ (resp. $ S(a,\varepsilon )$) be the closed ball (resp. the sphere)
with center $a$ and radius $\varepsilon$. We denote by $a \ast (S(a, \varepsilon ) \cap  A)$ the 
cone with vertex $a$ and basis $S(a, \varepsilon ) \cap  A$, i.e. the set of points in ${\Bbb R}^n$
of the form $\lambda a + (1-\lambda)x$, where $\lambda \in  [0,1]$ and $x \in  S(a, \varepsilon )
\cap  A$.

\begin{thm} For $\varepsilon  > 0$ sufficiently small, there is a semialgebraic homeomorphism
$h : \overline{B}(a, \varepsilon ) \cap  A \to  a \ast (S(a, \varepsilon ) \cap  A)$ such that
$\Vert h(x) -
 a\Vert  = \Vert  x - a\Vert $ and $h_{\vert S(a,\varepsilon ) \cap  A }= {\rm Id}$.
\end{thm}

\begin{proof}
We apply Hardt's theorem \ref{Hardt} to the mapping $p : A \to  {\Bbb R}$ defined by $p(x) = \Vert x
- a\Vert$. We obtain semialgebraic trivializations of $p$ over a finite semialgebraic partition 
of
${\Bbb R}$. We can assume that this partition has as member an interval $( 0,
\varepsilon _0 )$. Choose $\varepsilon$ such that $0 < \varepsilon  < \varepsilon  _0$. Since
$p^{-1}(\varepsilon)= (A \cap  S(a, \varepsilon  ))$, we have a semialgebraic homeomorphism
$$g :  p^{-1}(( 0, \varepsilon _0 )) \to  ( 0, \varepsilon _0) \times  (A \cap  S(a, \varepsilon 
))$$
such that  $g(x) = (\Vert  x - a\Vert , g_1(x))$, where the restriction of  $g_1$
to $S(a,\varepsilon ) \cap  A$ is the identity. Now define  $h :\overline{B}(a,
\varepsilon ) \cap  A \to  C_\varepsilon$  by
$$\left\{ \begin{array}{rcl}
h(x) &= &\left( 1 - \displaystyle{\Vert  x - a\Vert \over \varepsilon }\right) a +
\displaystyle{\Vert  x - a\Vert 
\over \varepsilon }\, g_1(x)\quad\mbox{\rm  if   }x \not= a\;,\\
h(a) &= &a\;.
\end{array}\right.$$
We can check that $h$ has the properties of the theorem. The inverse mapping of $h$ 
is defined by  
$$\left\{ \begin{array}{rcl}
h^{-1}( \lambda a + (1-\lambda) x) &= &g^{-1}((1-\lambda)\varepsilon , x)\ \mbox{\rm  for  }
\lambda 
\in  [0,1),\ x 
\in  S(a, \varepsilon ) \cap  A\;,\\
h^{-1}(a) &= &a\;.
\end{array}\right.$$
\end{proof}

\debexo
Let $A\subset {\Bbb R}^n$ be a semialgebraic set. For $r>0$, we denote by
$B(r)$ (resp. $\overline{B}(r)$, $S(r)$) the open ball (resp. the closed ball, the sphere) with
center $0$ and radius $r$ in ${\Bbb R}^n$.\\ 
1) Show that there exist $r>0$ and a continuous semialgebraic mapping
$$h_1 : A\cap\big({\Bbb R}^n\setminus B(r)\big)
\longrightarrow \big(A\cap S(r)\big)$$ 
such that 
$$ \forall x \in A\cap S(r), \quad h_1(x)=x$$
and
\begin{eqnarray*}
 h : A\cap\big({\Bbb R}^n\setminus B(r)\big) &\longrightarrow &\big(A\cap
S(r)\big)\times [r,+\infty)\\ 
x&\longmapsto& \big(h_1(x),\Vert x\Vert\big)
\end{eqnarray*} 
is a homeomorphism.
\\
2) Using 1), construct a semialgebraic mapping
$$H : A\times [0,1]\longrightarrow A$$ such that
\begin{eqnarray*}
 &\forall x\in A, \quad H(x,0)=x \et H(x,1)\in A\cap\overline{B}(r)\;,&\\
&\forall x\in A\cap\overline{B}(r),\ \forall t\in[0,1], \quad H(x,t)=x\;.&
\end{eqnarray*} 
($A\cap\overline{B}(r)$ is a semialgebraic deformation retract of $A$.)
\finexo

\subsection{Finiteness of the number of topological types}

\begin{thm}\label{finitudetyptop}
For every positive integers $d$ and $n$, there exist a positive integer $p = p(n, d)$ and
algebraic subsets $V_1 ,\ldots , V_p \subset  {\Bbb R}^n$, defined by polynomial equations of
degrees at most $d$, such that, for every algebraic subset
$W \subset  {\Bbb R}^n$ defined by polynomial equations of degrees at most $d$, 
there exist $i \in \{1,\ldots , p\}$ and a semialgebraic homeomorphism  $h :
{\Bbb R}^n \to  {\Bbb R}^n$  such that  $h(W) = V_i$.
\end{thm}

In other words, there are finitely many semialgebraic topological types of inclusions $V\subset{\Bbb
R}^n$, where $V$ is an algebraic subset defined by equations of degrees at most $d$. By
{\em (semialgebraic) topological type\index{topological type}} of inclusion, we mean an equivalence
class of subsets
$V
\subset   {\Bbb R}^n$ for the equivalence relation $(V \subset  {\Bbb R}^n) \sim (W 
\subset  {\Bbb R}^n)$ if there exists a (semialgebraic) homeomorphism $h : {\Bbb R}^n \to  
{\Bbb R}^n$ such that $h(V) = W$.\par
\vskip11pt
\begin{proof}
An algebraic subset $V \subset  {\Bbb R}^n$ 
defined by equations $P_1 =\ldots = P_q = 0$ of degrees $\leq d$ can always be regarded as
defined by only one equation of degree $\leq  2d$, that is $P_1^2+\ldots + P_q^2= 0$.
	A polynomial of degree $2d$ in $n$ variables has ${2d + n\choose 
n}= N(n, d)$ coefficients (check by induction on $n$). 
	We identify the space of coefficients with ${\Bbb R}^N$, and we denote by  
 $P_a \in  {\Bbb R}[X_1 ,\ldots , X_n]$ the polynomial of degree $\leq  2d$ corresponding to $a
\in  {\Bbb R}^N$. The set
$$ A = \{ (a, x) \in  {\Bbb R}^N \times  {\Bbb R}^n\  ; \  P_a(x) 
= 0 \}$$
is an algebraic subset of $ {\Bbb R}^N \times  {\Bbb R}^n$. Let $p : {\Bbb R}^N \times  {\Bbb R}^n
\to  {\Bbb  R}^N$ be the projection. Hardt's theorem implies that there is a finite
semialgebraic partition ${\Bbb R}^N = C_1\cup\ldots\cup  C_q$  such that, for every $i = 1,\ldots
,q$, there is a semialgebraic trivialization $h_i : C_i\times{\Bbb R}^n\to C_i\times{\Bbb R}^n$ of
$p$ over
$C_i$, compatible with $A$. Choose a point $a_i$ such that $P_{a_i}$ is a sum of squares of
polynomials in every $C_i$ containing such points, say $C_1,\ldots,C_p$. For $i=1,\ldots,p$, set
$$ V_i = \{ x \in  {\Bbb R}^n \ ;\  P_{a_i}(x) = 0 \}\; .$$
Any algebraic set $W$ as in the statement of the theorem is of the form
$$W = \{ x \in  {\Bbb R}^n \ ;\  P_a(x) = 0 \}\;,$$
where $a\in C_i$ for some $i\in\{1,\ldots,p\}$. The semialgebraic trivialization $h_i$
induces a semialgebraic homeomorphism $p^{-1}(a)\to p^{-1}(a_i)$ sending $A\cap p^{-1}(a)$ to
$A\cap p^{-1}(a_i)$ which, in turn, gives a semialgebraic homeomorphism $h:{\Bbb R}^n\to{\Bbb
R}^n$ such that $h(W)=V_i$.
\end{proof}

\section[Number of connected components]{Uniform bound on the number of connected components}

It follows from Theorem \ref{finitudetyptop} that there exists a positive integer $\varphi (n, d)$,
such that every algebraic subset of ${\Bbb R}^n$ defined by equations of degrees
$\leq  d$ has at most $\varphi (n, d)$ connected components\index{connected component}. Indeed, we
can take for
$\varphi (n, d)$ the maximum number of connected components of $V_1 ,\ldots , V_p$ (with the
notation of Theorem
\ref{finitudetyptop}). The aim of this section is to give an explicit upper bound for $\varphi (n,
d)$. First note that we have the lower bound $d^n\leq\varphi (n, d)$. Indeed, the system of
equations
 $(X_i-1)(X_i-2)\cdots(X_i-d)$, for
$i=1,\ldots,n$,  has a set of solutions consisting of $d^n$ real points. We assume everywhere in
this section that $d$ is a positive integer.

\begin{thm}\label{bornecc} Let $V \subset  {\Bbb R}^n$ be an algebraic subset defined by equations
of degrees $\leq  d$. The number of connected components\index{connected component} of 
$V$ is not greater than $d(2d - 1)^{n-1}$.
\end{thm}

This result is related to the Thom-Milnor bound\index{Thom-Milnor bound} on the sum of the Betti
numbers of a real algebraic set (the number of connected components is the Betti number $b_0$). For
this result, see for instance \cite{BR,BCR}. The proof of Theorem
\ref{bornecc} is essentially the proof of the Thom-Milnor bound, without its homological part.
\par

The proof proceeds by reducing to the case of a compact smooth hypersurface, and then looking for
extremal points of a function on this hypersurface. This is an example of the ``critical
point\index{critical!point} method''. This method can be applied to design algorithms (for deciding
whether a semialgebraic set is empty, for instance) with better complexity than the
c.a.d.\index{cylindrical algebraic decomposition}\ algorithm. See
\cite{Madrid} and the references cited there.
\par 

We state here a result which will be useful in this section. Let $M$ (resp. $N$) be a smooth
submanifold of ${\Bbb R}^n$ (resp. ${\Bbb R}^p$), and let $f:M\to N$ be a smooth map. For $x \in 
M$, denote by 
 $df_x : T_xM \to  T_{f(x)} N$ the tangent linear mapping from the tangent space to  $M$ at $x$ to
the tangent space to $N$ at 
$f(x)$. We say that $x$ is a {\em critical point\index{critical!point}} of $f$, and $f(x)$ a  {\em
critical value\index{critical!value}} of
$f$, if $df_x$ is not surjective.

\begin{thm}[Sard's theorem -- semialgebraic version]\label{Sard} Let $M\subset{\Bbb R}^n$ and
$N\subset{\Bbb R}^p$ be semialgebraic smooth submanifolds, and let
$f : M\to  N$ be a semialgebraic ${\cal C}^\infty $ mapping.  Then the set of critical values of
$f$ is a semialgebraic subset of $N$, of dimension\index{dimension} $<\dim N$.
\end{thm}

 For a proof, we refer to \cite{BCR} chap.\ 9 section 5, or
\cite{BR}  2.5.12. We propose a proof of a particular case in the following exercise.\par

\debexo\label{petitSard} Let $U$ be an open semialgebraic subset of ${\Bbb R}^n$. Let $Q$ and $f$
be polynomials on $R^n$. Assume that, for every $x\in U\cap Q^{-1}(0)$, $\grad Q(x)\not =0$ (hence,
$M=U\cap Q^{-1}(0)$ is a smooth hypersurface). Show that the set of critical values of
$f|_M$ is finite. Hints:\\ 1) Show that $x\in M$ is a critical point of $f|_M$ if and only $\grad
f(x)$ is colinear to $\grad Q(x)$. Deduce that the set $Z$ of critical points of $f|_M$ is
semialgebraic.\\ 2) Show that $f$ is constant along a smooth path in $Z$. Deduce that $f$ is
constant on each connected component of $Z$. Conclude.
\finexo

\subsection{Reduction to the case of a compact smooth hypersurface}

Assume that
$$V = \{ x \in  {\Bbb R}^n \ ;\ P_1(x) =\ldots = P_q(x) = 0 \}\;,$$  where $\deg P_i \leq  d$, and
let $P = P_1^2 +\ldots +  P_q^2$. We also assume that $V$ has at least two connected components.
Choose
$R>0$ larger than the maximum of the distances from the origin $0\in {\Bbb R}^n$ to the connected
components of $V$. The closed ball $\overline B$ with center $0$ and radius $R$ has a nonempty
intersection with every connected component of $V$. Hence, the number of connected components of
$\overline B \cap V$ is greater than or equal to the number of connected components of $V$.\par Let
${\cal F}$ be the finite set of connected components of $\overline B\cap V$. For
$C\in {\cal F}$, let
$K_C$ be the set of $x\in \overline B$ such that ${\rm dist}(x,C)={\rm dist}(x, (\overline B\cap
V)\setminus C)$. Let $K$ be the union of the sets $K_C$, for all $C\in {\cal F}$. The set $K$ is a
closed semialgebraic subset of $\overline B$, disjoint from $V$.
\par If $C_1$ and $C_2$ are different connected components of $\overline B\cap V$, every continuous
path in $\overline B$ joining a point of $C_1$ to a point of $C_2$ must intersect
$K_{C_1}$. Hence, each connected component of
$\overline B\setminus K$ contains at most one element of ${\cal F}$ (actually, exactly one).\par
 For
$0<\varepsilon$, set
$$Q_\varepsilon(x)=P(x)+\varepsilon(\Vert x\Vert^2 - R^2)\;.$$ Note that $Q_\varepsilon$ is a
polynomial of degree $2d$ in $x$. Since the polynomial
$P(x)$ is nonnegative on
${\Bbb R}^n$, the zeroset
$W_\varepsilon$ of
$Q_\varepsilon$ is contained in $\overline B$.  Let $\varepsilon_0>0$ be the minimum of
$P(x)/R^2$ for $x\in K$. Now assume $\varepsilon<\varepsilon_0$. It follows that
$W_\varepsilon$ is disjoint from $K$. Let $A$ be a connected component of
$\overline B\setminus K$ containing a $C\in {\cal F}$. Since $Q_\varepsilon$ takes nonpositive
values on $C\subset A$ and positive values on $K\cap \clos(A)$, its zeroset
$W_\varepsilon$ must have a nonempty intersection with $A$. This shows that the number of connected
components of $W_\varepsilon$ is greater than or equal to the number of connected components of
$\overline B\cap V$, which is  greater than or equal to the number of connected components of
$V$.\par Now we show that we can choose $\varepsilon$ such that $W_\varepsilon$ is a smooth
hypersurface. This means that for every point $x\in W_\varepsilon$,  the partial derivatives
$\partial Q_\varepsilon/\partial x_i$ do not all vanish.
\par
 First, we claim that the set
$$X=\{(x,\varepsilon)\in {\Bbb R}^n\times{\Bbb R}\;;\ \varepsilon>0 \et Q_\varepsilon(x)=0\}\;.$$
is a smooth hypersurface of ${\Bbb R}^{n+1}$. Since
$\partial Q_\varepsilon/\partial \varepsilon=\Vert x\Vert^2 - R^2$, this partial derivative can
vanish at $(x,\varepsilon)\in X$ only if $\Vert x\Vert=R$ and $P(x)=0$. Since $\partial P/\partial
x_i$ vanishes when $P=0$, we have moreover $(\partial Q_\varepsilon/\partial x_i)
(x,\varepsilon)=2\varepsilon x_i$. Hence, one of the partial derivatives $(\partial
Q_\varepsilon/\partial x_i)$ is nonzero at $(x,\varepsilon)$. This proves the claim.\par Second, we
claim that the function $f:X\to R$ defined by $f(x,\varepsilon)=\varepsilon$ has finitely many
critical values. This follows from the semialgebraic version of Sard's theorem \ref{Sard}.
Actually, it is sufficient to use here the result of Exercise
\ref{petitSard}. Hence, we can choose $\varepsilon$, with $0<\varepsilon<\varepsilon_0$ and
$\varepsilon$ not a critical value of $f$. This implies that
$\grad Q_\varepsilon=(\partial Q_\varepsilon/\partial x_1,\ldots,
\partial Q_\varepsilon/\partial x_n)$ is never $0$ on $W_\varepsilon$.\par Now it suffices to show
that the number of connected components of $W_\varepsilon$ is not greater than $d(2d-1)^{n-1}$.
This will be done in the next lemma.

\subsection{The case of a compact smooth hypersurface}

\begin{lem}\label{compacthyper} Let $Q : {\Bbb R}^n\to {\Bbb R}$ be a polynomial of degree $2d$.
Assume that
$W=Q^{-1}(0)$ is compact and $\grad Q$ has no zero on $W$ (hence, $W$ is a smooth compact
hypersurface). Then the number of connected components\index{connected component} of $W$ is at most
$d(2d-1)$.
\end{lem}

\begin{proof}
 The connected components of $W$ are compact smooth hypersurfacess.  The coordinate function
$x_n$ has a maximum and a minimum on each of these components. At a point where
$x_n$ reaches its maximum or minimum,  the tangent hyperplane to
$W$ is parallel to the hyperplane $x_n = 0$, which means that the first $n-1$ coordinates of $\grad
Q$ vanish. Hence, the points of $W$ where $x_n$ reaches an extremum are solutions of the system of
equations
$$ (S)\quad  \left\{\begin{array}{rcl} Q(x) &= &0\\
\displaystyle{\partial Q\over \partial x_1}\, (x) &= &0\\ &\vdots &\\
\displaystyle{\partial Q\over \partial x_{n-1}}\, (x) &= &0
\end{array}\right. \qquad  (n\ \mbox{\rm equations, } n\ \mbox{\rm variables})\;.$$

\begin{lem}\label{Morse} We can choose the coordinates in ${\Bbb R}^n$ such that all real solutions
of the above system 
$(S)$ are nondegenerate, i.e. solutions where the jacobian determinant
$$
\det  \left(\begin{array}{ccc}
\displaystyle{\partial Q\over\partial x_1}&\ldots&\displaystyle{\partial Q\over\partial x_n}
\\
\vbox to 20pt{}\displaystyle{\partial ^2Q\over\partial x_1\partial
x_1}&\ldots&\displaystyle{\partial ^2Q\over\partial  x_1\partial x_n}\\
\vdots&\vdots&\vdots\\
\vbox to 20pt{}\displaystyle{\partial ^2Q\over\partial x_{n-1}\partial x_1}&\ldots
&\displaystyle{\partial  ^2Q\over\partial x_{n-1}\partial x_n}
\end{array}\right)$$ does not vanish.
\end{lem}

	Assume for the moment that Lemma \ref{Morse} is proved. Then, by Bezout's theorem (cf.\ for
instance \cite{BR}, appendix B, or \cite{BCR}, chap.\ 9), the number of nondegenerate real
solutions of the system
$(S)$ is
$\leq 2d(2d-1)^{n-1}$, which is the product of the degrees of the equations. Since each connected
component of
$W$ has at least two points where $x_n$ has an extremum, $W$ has at most
$d(2d-1)^{n-1}$ connected components. This concludes the proof of Lemma
\ref{compacthyper} and also the proof of Theorem \ref{bornecc}.
\end{proof}

{\par\noindent{\it Proof of Lemma \ref{Morse}\/: }\nopagebreak\normalsize}
Apply the semialgebraic
version of Sard's theorem \ref{Sard} to the map
$$\varphi  =  \grad Q / \Vert  \grad Q\Vert   :  W   \longrightarrow   S^{n-1}\;, $$ where
$S^{n-1}$ is the unit sphere in ${\Bbb R}^n$. Since the set of critical values of $\varphi $  is of
dimension $< n-1$, we can find a pair of antipodal points $b$ and $-b$ in $S^{n-1}$, which are not
critical values. After rotating the coordinate axes, we can assume that $b = (0,\ldots ,0,1)$. 
Observe that the real solutions $a$ of the system $(S)$ are exactly the points $a\in W$ such that
$\varphi (a) = (0,\ldots  ,0,\pm1)$. For such an $a$, the tangent hyperplanes $T_aW$
 and $T_{f(a)}S^{n-1}$ are both $x_n = 0$, and we have 
$\partial Q/\partial x_i (a) = 0$, for $i = 1,\ldots ,n-1$. The matrix of 
$d\varphi _a$, in the coordinates $( x_1 ,\ldots , x_{n-1})$, is
$${1\over \Vert  \grad Q(a)\Vert }  \left(  {\partial ^2Q\over\partial x_i \partial x_j}  (a)
\right)_{i=1,\ldots ,n-1\ ;\ j=1,\ldots  ,n-1}$$ and, therefore,
$$\Delta = \det \left( {\partial ^2Q\over\partial x_i \partial x_j}  (a) 
\right)_{i=1,\ldots ,n-1\ ;\ j=1,\ldots ,n-1}\not= 0\;. $$ Since the value of the jacobian
determinant of the system $(S)$ at $a$ is equal to
$\pm\Delta \,\displaystyle{\partial Q\over \partial x_n} (a)$, it is nonzero.
 {\linebreak[2]\hspace*{\fill}$\square$\ifdim\lastskip<12pt
    \removelastskip \penalty-200 \vskip12pt \fi}

\subsection{Bound for the number of connected components of a semialgebraic set}

We have just seen that, for an algebraic set, the bound on the number of connected components
depends only on the degree of the equations and not on the number of these equations. In the case
of a semialgebraic set defined by polynomial inequalities, we have to take into account the number
of these inequalities.

\debexo  Prove that any finite union of open intervals in ${\Bbb R}$ can be defined by a system of
inequalities of degree $2$.
\finexo

\begin{prop}\label{borneccsa}
 Let $(S)$ be a system of $s$ polynomial equations and inequalities in
$k$ variables, of degrees at most $d \geq  2$. The number of connected components\index{connected
component} of the set of solutions of
$(S)$ in ${\Bbb R}^n$ is not greater than $d(2d - 1)^{k+s-1}$.
\end{prop}

\begin{proof} Assume $P>0$ is a strict polynomial inequality in $(S)$. Choose
$\varepsilon>0$ so small that there is a point $x$ in each connected component of the set of
solutions of $(S)$, with
$P(x)>\varepsilon$. Replacing $P > 0$ with $P- \varepsilon\geq 0$ in $(S)$ can only increase the
number of connected components of the set of solutions. Hence, we can assume that $(S)$ contains
only lax inequalities $Q_1\geq 0$,\ldots, $Q_t\geq 0$ and equations. Next, we replace every lax
inequality
$Q_i\geq0$ with the equation $Q_i- T_i^2=0$, introducing a new variable $T_i$ for each inequality.
This replacement can only increase the number of connected components of the set of solutions.
Finally, we have a system of $s$ equations in $k+t\leq k+s$ variables of degrees $\leq d$. By
Theorem \ref{bornecc}, the number of connected components of the set of solutions is not greater
than $d(2d - 1)^{k+s -1}$.
\end{proof}

The bound of Proposition \ref{borneccsa} is very coarse. Using more sophisticated arguments, one
can obtain a bound  with a polynomial dependence on the number of equations and inequalities, of
the form $s^kO(d)^k$ (cf. \cite{Madrid}).

\section{An application to lower bounds}

This section is devoted to an application of the bound on the number of connected
components\index{connected component} of a semialgebraic set, due to Ben-Or (Proc. 15th ACM ann.
Symp. on Theory of Comp., 80-86 (1983)).
First, we describe a model of algorithm to decide whether an element  $x=( x_1 ,\ldots , x_n) \in 
{\Bbb R}^n $ satisfies a boolean combination of sign conditions on polynomials in $ n$
variables  (in other words, to decide whether $x$ belongs to a given semialgebraic subset $W
\subset  {\Bbb R}^n$). This model is an {\em algebraic computation tree\index{algebraic computation
tree}}.  Such a tree has one root and several leaves. The vertices different from the root have one
father. The vertices different from the leaves have one or two sons. A vertex $v$ with one son is
labelled with a variable $U_v$ and an instruction $U_v := a * b$, 
where $*$ is an arithmetic operation $(+, -, \times)$ and 
$a $ and $b$ are either $x_i$,  
$i\in\{1,\ldots ,n\}$, or a real constant, or a variable 
$U_{v'}$ with $v'$ ancestor of $v$. A vertex $v$ with two sons is labeled with a test $a\, 
?\, 0$, where $a$ is as above and $?$ is $=$, $>$, or 
$\geq $. A leaf $v$ is labeled with a boolean constant $b_v$ (true or false). 
The algorithm modelled by such a tree has as inputs  
$n$-tuples $x_1 ,\ldots , x_n$ of real numbers and goes down in the tree starting from the root. At
a vertex $v$ with one son, it computes the value of $U_v$ following the instruction of the label.
At a vertex with two sons, it chooses the left (resp. right) son if the answer to the test of
the label is yes (resp. no). When the algorithm arrives to a leaf, it returns the boolean
constant in the label of this leaf. The {\em cost\index{cost}} of an algorithm (in this model) is
the maximal length of a path (from root to leaf) taken by an input $( x_1 ,\ldots , x_n) \in  {\Bbb
R}^n$.

\begin{figure}
\setlength{\unitlength}{1mm}

\begin{picture}(120,85)
\put(35,75){\framebox(30,10){$U_1:=y-x$}}
				\put(50,75){\line(0,-1){4}}
\put(35,59){\test{$U_1>0$ ?}}
				\put(35,65){\line(-1,0){15}}
	   \put(27,62){\makebox(0,0){yes}}
				\put(20,65){\line(0,-1){10}}
\put(10,45){\framebox(20,10){false}}
				\put(65,65){\line(1,0){15}}
	   \put(72,62){\makebox(0,0){no}}
				\put(80,65){\line(0,-1){10}}
\put(65,45){\framebox(30,10){$U_2:=y\times U_1$}}
				\put(80,45){\line(0,-1){5}}
\put(65,30){\framebox(30,10){$U_3:=U_2-1$}}
	 		\put(80,30){\line(0,-1){4}}
\put(65,14){\test{$U_3=0$ ?}}
				\put(65,20){\line(-1,0){15}}
	   \put(57,17){\makebox(0,0){yes}}
				\put(50,20){\line(0,-1){10}}
\put(40,0){\framebox(20,10){true}}
				\put(95,20){\line(1,0){15}}
	   \put(102,17){\makebox(0,0){no}}
				\put(110,20){\line(0,-1){10}}
\put(100,0){\framebox(20,10){false}}
\end{picture}
\caption{An algebraic computation tree deciding whether $y^2-xy=1\et y\leq x$.}
\end{figure}
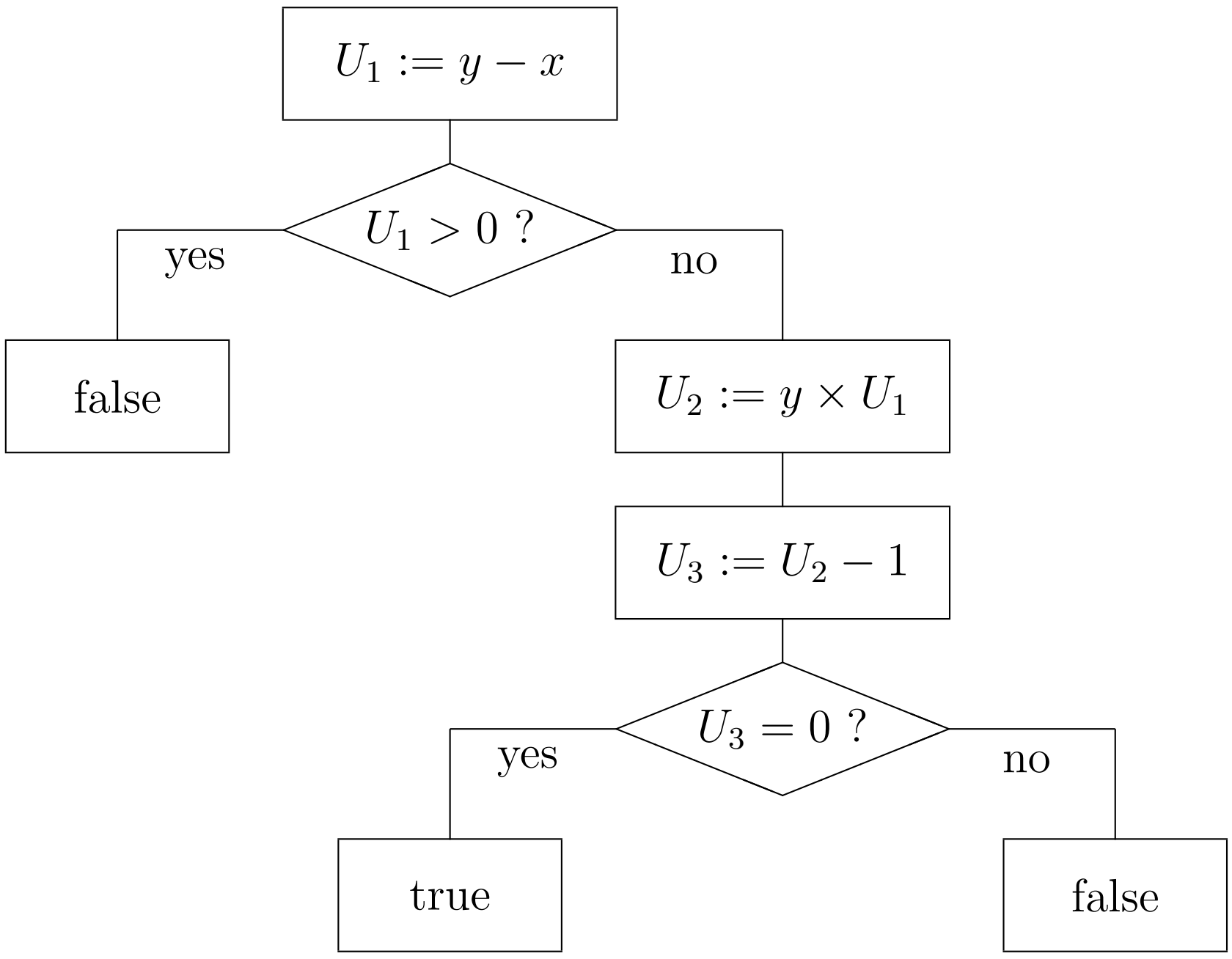

\begin{thm} If an algorithm with cost $c$ decides whether $x \in  
W$, where  $W$ is a semialgebraic subset of ${\Bbb R}^n$, then the number of connected
components\index{connected component} of $W$ is not greater than $2^{2n+5c}$.
\end{thm}

\begin{proof}
Let $v$ be a leaf labelled with $b_v = \mbox{\rm true}$. Denote by $W_v$ the semialgebraic subset
of $W$ consisting of those inputs  $x \in  {\Bbb R}^n$ for which the algorithm arrives to the leaf
$v$. Consider the system $(S_v)$  obtained in the following way.
\begin{itemize}
\item For each ancestor $v'$ of $v$ with one son, take the equation 
 $U_{v'}= a * b$ which is in the label of $v'$. 
\item For each ancestor $v'$ of $v$ with two sons, take the equation or inequality $a\, ?\, 0$
in the label of $v'$ if $v$ is a heir of $v'$ on the left side, and its negation if $v$ is a heir
of $v'$ on the right side.
\end{itemize} 
The system $(S_v)$ 
has $s$ equations and inequalities in the variables 
$$X_1,\ldots , X_n, U_{v_1},\ldots , U_{v_m},$$ 
where $m \leq  s$. Assume $W_v 
\not= \emptyset$. There are inputs for which the algorithm arrives to $v$. Hence $s \leq  c$.
Finally,
$( x_1 ,\ldots , x_n)
\in  W_v$ if and only if there exists $(u_{v_1},\ldots , u_{v_m}) \in  {\Bbb R}^m$ such that 
$( x_1,\ldots , x_n, u_{v_1},\ldots , u_{v_m}) $ is a solution of $(S_v)$. Hence, the number of
connected components of $W_v$ is not greater than the number of connected components of the set of
solutions of $(S_v)$. Since all equations and inequalities in
$(S_v)$  have degree $\leq  2$, Proposition \ref{borneccsa} implies that the number of connected
components of $W_v$ is not greater than $2 \times  3^{n+m+ s
-1}
\leq  2
\times 
 3^{n+2c-1}$. The number of leaves $v$ with $b_v = \mbox{\rm true}$ and $W_v 
\not= \emptyset$ is at most $2^c$, since the paths from the root to the leaves which are taken for
some input have length
$c$ at most, and each vertex has at most two sons. Since $W$ is the union of these $W_v$, the
number of connected components of
$W$ is not greater than $2^c \times  2 \times  3^{n+2c-1} \leq  2^{2n+5c}$.
\end{proof}

\begin{cor}
The cost (for the algebraic computation tree model) of an algorithm deciding whether $n$ real
numbers $( x_1 ,\ldots , x_n)$ are all distinct is at least $\Omega(n \log n)$ (recall that
$f=\Omega (g)$ means $g=O(f)$).
\end{cor}

Remark that there are algorithms solving this problem with cost $O(n \log n)$ (sorting by divide
and conquer, for instance).\par
\medskip

\begin{proof}
The set $$W = \{ ( x_1 ,\ldots , x_n) \in  {\Bbb R}^n \ ;\ \forall i \not= j, 
x_i \not= x_j \}$$ has $n!$ connected components (there is a faithful and transitive action of
the group of permutation of $\{1,\ldots ,n\}$ on the set of connected components of $W$). The cost
$c$ of an algorithm deciding whether
$x \in  W$ must satisfy $n !  \leq  2^{2n+5c}$. Taking logarithm and using $ n
\log n=O(\log(n !))$ we obtain
$n \log n= O(2n+5c)$.
Hence $c =\Omega (n \log n)$.
\end{proof}

\debexo
(This exercise is taken from a paper by J.L. Monta{\~n}a, L.M. Pardo and T. Recio in Effective
Methods in Algebraic Geometry,\linebreak Birkh{\"a}user 1991).\\
1) We denote by ${\cal F}_d$ the family of all algebraic subsets of ${\Bbb R}^n$
defined by an equation $P=0$, where $P$ is a nonzero polynomial of degree $\leq d$.
Let $W$ be a semialgebraic subset of
${\Bbb R}^n$. Show that there exists  $i\in {\Bbb N}$ such that, for every 
$H$ in ${\cal F}_d$, the number of connected components of
$H\cap W$ is $\leq i$. Let $I_d(W)$ be the smallest such integer. \\
2) We now assume $W$ to be defined by a system of
$\ell$ polynomial equations and inequalities in $n$ variables, of degrees at most
$d\geq 2$. Give an upper bound for $I_d(W)$. \smallskip\\ 
We shall now find a lower bound for the cost of an algorithm deciding the following problem
(``big hole'' problem): given 
$n\geq 2$ real numbers $x_1,\ldots,x_n$, is there a closed interval of length $1$ containing
no $x_i$ and contained in the convex hull of the $x_i$'s in ${\Bbb R}$?\\
3) Let $W_n\subset {\Bbb R}^n$ be the subset of all
$(x_1,\ldots,x_n)$ for which there is no big hole. Show that $W_n$ is a connected
semialgebraic set.\\
4) We assume that $(x_1,\ldots,x_n)\in W_n$ and
$x_1\leq x_2\leq\ldots\leq x_n$. Show that 
$$\sum_{i<j} (x_i-x_j)^2\leq \sum_{k=1}^{n-1} k(n-k)^2\;.$$ 
Show that the equality holds if and only if $x_{i+1}=x_i+1$ for
$i=1,\ldots,n-1$. Deduce that the intersection of $W_n$ with the algebraic set defined by the
equation
$$\sum_{i<j} (x_i-x_j)^2= \sum_{k=1}^{n-1} k(n-k)^2$$ is the union of $n!$ disjoint lines.
Hence $I_2(W_n)\geq n!$.\\
5) Show that the cost (in the algebraic computation tree model) of an algorithm deciding the big
hole problem is at least
$\Omega(n\,\log n)$.
\finexo

\printindex
\end{document}